\documentclass[]{amsart}
\usepackage{amssymb}
\usepackage{amsmath}
\usepackage{amsthm}
\usepackage{amscd}
\usepackage{epsfig}
\usepackage{graphicx,psfrag}
\usepackage[all]{xy}

\input epsf
\def\relabelbox{%
  \hbox\bgroup%
}%
\def\endrelabelbox{%
\egroup%
}%
\def\relabel #1#2 {%
  \special{ps:/a {} def}%
  \smash{\rlap{#2}}%
}%
\def\adjustrelabel <#1,#2> #3#4 {%
  \special{ps:/a {} def}%
  \smash{\rlap{\kern #1 \raise #2\hbox{#4}}}%
}%
\def\extralabel <#1,#2> #3 {\smash{\rlap{\kern #1 \raise #2\hbox{#3}}}}%

\newcommand{\etalchar}[1]{$^{#1}$}
\def\MCG{{\mathcal MCG}}

\def\C{{\mathcal C}}

\newtheorem{thm}{Theorem}[section]
\newtheorem{thm*}{Theorem*}
\newtheorem{prop}[thm]{Proposition}
\newtheorem{cor}[thm]{Corollary}

\newtheorem{qn}[thm]{Question}
\newtheorem{lem}[thm]{Lemma}

\newtheorem{cvn}[thm]{Convention}

\theoremstyle{definition}
\newtheorem{defn}[thm]{Definition}
\newtheorem{exa}[thm]{Example}

\theoremstyle{remark}
\newtheorem{rmk}[thm]{Remark}
\newtheorem{n}[thm]{Notation}

\def\square{\hfill${\vcenter{\vbox{\hrule height.4pt \hbox{\vrule
width.4pt
height7pt \kern7pt \vrule width.4pt} \hrule height.4pt}}}$}
\newenvironment{pf}{\medskip\noindent {\it Proof:}\quad}{\square
\vskip 12pt}

\def\R{{\mathbb R}}
\def\free{F}
\def\Z{{\mathbb Z}}
\def\N{{\mathbb N}}
\def\Q{{\mathbb Q}}

\def\C{{\mathbb C}}

\def\ga{{\mathcal{G}}_A}
\def\C{{\mathcal C}}
\def\calc{{\mathcal C}}

\def\G{{\mathcal G}}
\def\g{{\mathcal G}}

\def\H{{\mathcal H}}
\def\lll{\mathcal{L}}
\def\P{{\mathcal P}}
\def\pp{{\mathcal P}}
\def\calp{{\mathcal P}}

\def\Rcal{{\mathcal R}}
\def\oo{{\mathcal O}}
\def\nn{{\mathcal N}}
\def\bb{{\mathcal B}}
\def\fg{{\mathfrak g}}
\def\fc{{\mathfrak c}}
\def\fp{{\mathfrak p}}
\def\MCG{\mathcal {MCG}}
\def\MCGS{\mathcal {MCG}(S)}
\def\mcg{\MCG}

\def\pmcg{\mathcal{PMCG}}

\def\nbhd{\mathcal{N}}
\def\diam{{\rm{diam}}}
\def\bg{\mathbf{G}}
\def\bs{\mathbf{S}}
\def\bn{\mathbf{N}}
\def\ban{\mathbf{AN}}
\def\coneomega{{\rm{Cone}}_{\omega}}
\def\ocone{{\rm{lim}}_{\omega}}
\def\con{{\rm{Cone}}_{\omega}}

\def\dist{{\rm{dist}}}
\def\cone{\coneomega}
\def\into{\hookrightarrow}
\def\ulim{\lim_\omega}
\def\ulimi{\lim_\omega}

\def\<{\langle}
\def\>{\rangle}
\def\SLNZ{SL_{n}(\Z)}
\long\def\Restate#1#2#3#4{
\medskip\par\noindent
{\bf #1 \ref{#2} #3}
{\it #4}\par\medskip
}
\def\autfn{{\rm{Aut}}(F_{n})}
\def\sautfn{{\rm{SAut}}(F_{n})}
\def\outfn{{\rm{Out}}(F_{n})}

\def\pants{\C_{P}(S)}
\def\pant{\C_{P}}
\DeclareMathOperator{\CAT}{CAT}

\newcommand {\me}{\medskip}

\newcommand\suchthat{\bigm|}
\newcommand\infinity\infty
\newcommand\wt\widetilde
\newcommand\inject\hookrightarrow
\newcommand\union\cup
\newcommand\reals{{\mathbf R}}
\newcommand\abs[1]{\left| #1 \right|}
\newcommand{\from}{\colon\thinspace}
\newcommand{\co}{\colon\thinspace}
\newcommand\join\Lambda
\newcommand\cross\times
\newcommand {\inv}{^{-1}}
\newcommand\lub\vee

\newcommand\glb\wedge
\newcommand\B{\mathcal{B}}

\renewcommand\paragraph[1]{\bigskip\noindent\textbf{#1} }
\newcommand\Id{{\mathrm Id}}
\newcommand {\gfrak}{\mathfrak g}
\newcommand {\p}{\mathfrak p} 
\newcommand {\q}{\mathfrak q} 
\newcommand{\aaa}{{\mathcal A}}
\newcommand{\atgs}{asymptotically tree-graded }

\newcommand{\atgrs}{asymptotically tree-graded with respect to }

\newcommand{\semidirect}{\rtimes}
\newcommand{\subgroup}{<}

\newcommand{\intersect}{\cap}
\newcommand{\even}{\text{even}}
\newcommand{\odd}{\text{odd}}
\newcommand{\la}{\langle}
\newcommand{\ra}{\rangle}

\begin{document}
\title{Thick metric spaces, relative hyperbolicity, \\ and
quasi-isometric rigidity}
\author{Jason Behrstock}
\address{Department of Mathematics\\ Columbia University\\ 
2990 Broadway\\ New York, NY 10027}
\email{jason@math.columbia.edu}
\author{Cornelia Dru\c{t}u}
\address{Mathematical Institute\\
24-29 St Giles, Oxford OX1 3LB\\
United Kingdom.}
\email{drutu@maths.ox.ac.uk}
\author{Lee Mosher}
\address{Department of Mathematics and Computer Science\\
Rutgers University at Newark\\
Newark, NJ 07102}
\email{mosher@andromeda.rutgers.edu}
\date{\today}
\thanks{To appear in \emph{Mathematische Annalen}.}

\begin{abstract}
    We study the geometry of nonrelatively hyperbolic groups. Generalizing a
    result of Schwartz, any quasi-isometric image of a non-relatively hyperbolic
    space in a relatively hyperbolic space is contained in a bounded neighborhood
    of a single peripheral subgroup. This implies that a group being relatively
    hyperbolic with nonrelatively hyperbolic peripheral subgroups is a
    quasi-isometry invariant. As an application, Artin groups are relatively
    hyperbolic if and only if freely decomposable.

    We also introduce a new quasi-isometry invariant of metric spaces called
    \emph{metrically thick}, which is sufficient for a metric space to be
    nonhyperbolic relative to any nontrivial collection of subsets. Thick
    finitely generated groups include: mapping class groups of most surfaces;
    outer automorphism groups of most free groups; certain Artin groups; and
    others. Nonuniform lattices in higher rank semisimple Lie groups are thick and
    hence nonrelatively hyperbolic, in contrast with rank one which provided the
    motivating examples of relatively hyperbolic groups.  Mapping class groups
    are the first examples of nonrelatively hyperbolic groups having cut points
    in any asymptotic cone, resolving several questions of Drutu and Sapir about
    the structure of relatively hyperbolic groups. Outside of group theory,
    Teichm\"{u}ller spaces for surfaces of sufficiently large complexity are
    thick with respect to the Weil-Peterson metric, in contrast with
    Brock--Farb's hyperbolicity result in low complexity.
\end{abstract}

\maketitle
\pagestyle{myheadings}
\markboth{JASON BEHRSTOCK, CORNELIA DRU\c{T}U, AND LEE MOSHER}
{THICK METRIC SPACES, RELATIVE HYPERBOLICITY, \& QUASI-ISOMETRIC
RIGIDITY}

\setcounter{tocdepth}{1}

\tableofcontents

\section{Introduction}

Three of the most studied families of groups in geometric group theory
are the mapping class group of a surface of finite type, $\MCG(S)$;
the outer automorphism group of a finite rank free group, $\outfn$;
and the special linear group, $\SLNZ$. Despite the active interest in
these groups, much of their quasi-isometric structure remains
unknown, particularly for the first two families. We introduce the
notion of a \emph{thick} group (or more generally, metric space), a
property which is enjoyed by all groups in each of the families
$\MCG(S)$, $\outfn$, and $SL_{n}(\Z)$ except in the lowest complexity
cases where the groups are actually hyperbolic. The notion of
thickness helps unify the study of these groups and casts light on
some of their geometric properties.

Before proceeding, we recall some relevant developments. In
\cite{Gromov:hyperbolic}, M. Gromov introduced the notion of a
relatively hyperbolic group. The theory of relatively hyperbolic
groups was developed by Farb in \cite{Farb:RelHyp}, then further
developed in \cite{Bowditch:RelHyp}, \cite{Dahmani:thesis},
\cite{Osin:RelHyp}, \cite{Yaman:RelHyp}, and
\cite{DrutuSapir:TreeGraded}. Several alternate characterizations
of relative hyperbolicity have been formulated, all of them more
or less equivalent to each other. We recall the definition due to
Farb. In the sequel $G$ denotes a finitely generated group endowed
with a word metric, $\mathcal{H}= \{ H_1,...,H_n\}$ is a finite
family of subgroups of $G$ and $\mathcal{L}\mathcal{H}$ denotes
the collection of left cosets of $\{ H_1,..., H_n \}$ in $G$. The
group $G$ is \emph{weakly hyperbolic relative to} $\mathcal{H}$ if
collapsing the left cosets in $\mathcal{L}\mathcal{H}$ to finite
diameter sets, in a Cayley graph of $G$, yields a
$\delta$--hyperbolic space. The subgroups $H_1,..., H_n$ are
called \emph{peripheral subgroups}.

The group $G$ is \emph{(strongly) hyperbolic relative to}
$\mathcal{H}$ if it is weakly hyperbolic relative to $\mathcal{H}$
and if it has the \emph{bounded coset property}. This latter
property, roughly speaking, requires that in a Cayley graph of $G$
with the sets in $\mathcal{LH}$ collapsed to bounded diameter
sets, a pair of quasigeodesics with the same endpoints travels
through the collapsed $\mathcal{LH}$ in approximately the same
manner.

In \cite[$\S 8$ and Appendix]{DrutuSapir:TreeGraded}, Dru\c{t}u,
Osin and Sapir provide a geometric condition which characterizes
relative hyperbolicity of a group. They show that $G$ is
hyperbolic relative to $\mathcal{H}$ if and only if any
\emph{asymptotic cone} of $G$ is \emph{tree-graded} with respect
to the collection of \emph{pieces} given by ultralimits of
elements in $\mathcal{L}\mathcal{H}$ (see Section \ref{prel} for
definitions). In particular any asymptotic cone of $G$ has
(global) cut-points.

\medskip

The asymptotic characterization of relative hyperbolicity
mentioned above is in turn equivalent to three metric
properties in the Cayley graph of $G$ (formulated without
asymptotic cones), which are approximately as follows:

\begin{itemize}
\item[$(\alpha_1)$] Finite radius neighborhoods of distinct elements
in
$\mathcal{L}\mathcal{H}$ are either disjoint or intersect in sets
of uniformly bounded diameter;
\item[$(\alpha_2)$] geodesics diverging slower than linearly
from a set $gH_i$ in $\mathcal{L}\mathcal{H}$ must intersect a
finite radius neighborhood of $gH_i$;
\item[$(\alpha_3)$] fat geodesic polygons must stay close to
a set in $\mathcal{L}\mathcal{H}$ (``fat'' here is the contrary of
``thin'' in its metric hyperbolic sense; see Definition
\ref{deffat}).
\end{itemize}

This definition of relative hyperbolicity also makes sense in a
general metric setting: a geodesic metric space $X$ is said to be
\emph{asymptotically tree-graded} (ATG in short) with respect to a
collection $\mathcal{A}$ of subsets of $X$ (called
\textit{peripheral subsets}) if the three conditions above hold
with $G$ replaced by $X$ and $\mathcal{L}\mathcal{H}$ replaced by
$\mathcal{A}$ (see also \cite{BrockFarb:curvature} for another
metric version of the notion of relative hyperbolicity). For
instance, the complementary set in $\mathbb{H}^3$ of any family of
pairwise disjoint open horoballs is asymptotically tree-graded
with respect to the collection of boundary horospheres. It was
recently proven by Dru\c{t}u that if a group is asymptotically
tree-graded in a metric sense, that is with respect to a
collection $\mathcal{A}$ of subsets, then it is relatively
hyperbolic with respect to some family of subgroups
\cite{Drutu:RelHyp} (see Theorem \ref{t3} in this paper). The
converse of the above statement was shown in
\cite{DrutuSapir:TreeGraded} (see Theorem
\ref{DrutuSapir:TreeGraded:relhypclassification}).

\medskip

\begin{cvn}\label{convnrh}
Throughout the paper, we exclude the trivial case of a metric
space $X$ asymptotically tree-graded with respect to a
collection $\mathcal{A}$ where some finite radius neighborhood of
some subset $A\in \mathcal{A}$ equals $X$. In the case of an
infinite group, $G$, hyperbolic relative to a collection of
subgroups, the trivial case we are excluding is where one of the
subgroups is $G$. (Note that a group is never hyperbolic relative to a finite index subgroup, so we need not exclude such cases in our convention.)

When a group contains no collection of proper subgroups with
respect to which it is relatively hyperbolic, we say the group is
\emph{not relatively hyperbolic (NRH)}.
\end{cvn}

\medskip

Thickness is, in many respects, opposite to relative
hyperbolicity. The notion of thickness is built up inductively. A
geodesic metric space is \textit{thick of order zero} if it is
\textit{unconstricted}, in the terminology of
\cite{DrutuSapir:TreeGraded}, that is: for at least one sequence
of scaling constants $d=(d_n)$ and one ultrafilter, all asymptotic
cones constructed by means of $d$ and $\omega$ are without
(global) cut-points. If the metric space is a group then this is
equivalent to the condition that \emph{at least one} asymptotic cone is
without cut-points. See Section \ref{sunc} for details, and for a
list of examples of groups thick of order zero
 (unconstricted).
A metric space is \emph{thick of order} $n$ if, roughly speaking,
it can be expressed as a coarse union of a network of subspaces
thick of order $n-1$, each quasi-isometrically embedded, so that
two adjacent subspaces in this network have infinite coarse
intersection. The exact definition of thickness can be found in
Section \ref{section:thick}. Because thickness is a quasi-isometry
invariant, thickness of a finitely generated group $G$ is
well-defined by requiring that the Cayley graph of a finite
generating set of $G$ be a thick metric space. Thick metric spaces
behave very rigidly when embedded into asymptotically tree-graded
metric spaces in particular we obtain (see
Theorem~\ref{thick:peripheral} for a generalization of this
result):

\Restate{Corollary}{crit}{(Thick spaces are not asymptotically
tree-graded).}{If $X$ is a thick metric space, then $X$ is not
asymptotically tree-graded. In particular, if $X$ is a finitely
generated thick group, then $X$ is not relatively hyperbolic.}

The following result puts strong restrictions on
how NRH groups can be quasi-isometrically embedded in ATG spaces.

\Restate{Theorem}{qinonrelhyp}{(NRH subgroups are peripheral).}{
Let $(X,\dist_X)$ be a metric space \atgrs a collection $\aaa$ of
subsets. For every $L\geq 1$ and $C\geq 0$ there exists $R=R(L,C,
X,\aaa)$ such that the following holds. If $G$ is a finitely
generated group endowed with a word metric $\dist$ and $G$ is not
relatively hyperbolic, then for any $(L,C)$-quasi-isometric
embedding $\q \co (G,\dist ) \to (X,\dist_X)$  the image $\q (G)$
is contained in the radius $R$ neighborhood of some $A\in \aaa$. }

Note that in the theorem above the constant $R$ \emph{does not}
depend on the group~$G$.

This theorem shows that the presence of NRH (in particular thick)
peripheral subgroups in a relatively hyperbolic group
``rigidifies'' the structure. A similar rigidity result, with
additional hypotheses on both the domain and the range plays a key
role in Schwartz's quasi-isometric classification of rank one
non-uniform lattices in semisimple Lie groups
\cite{Schwartz:RankOne}. Dru\c{t}u--Sapir proved a similar
rigidity result under the assumption that the domain is
unconstricted \cite{DrutuSapir:TreeGraded}; using work of
\cite{DrutuSapir:TreeGraded} allows one to obtain the following
theorem. (For special cases of this result see also
Theorem~\ref{cutp1} and Theorem~\ref{thick:peripheral} in this
paper or other results in \cite{DrutuSapir:TreeGraded}.)

\Restate{Theorem}{sqir}{(Quasi-isometric rigidity of hyperbolicity
relative to NRH subgroups).} {If $\Gamma$ is a finitely generated
group hyperbolic relative to a finite collection of finitely
generated subgroups $\G$ for which each $G\in\G$ is not relatively
hyperbolic, then any finitely generated group $\Gamma'$ which is
quasi-isometric to $\Gamma$ is hyperbolic relative to a finite
collection of finitely generated subgroups $\G'$ where each
subgroup in $\G'$ is quasi-isometric to one of the subgroups in
$\G$.}

In \cite{Drutu:RelHyp} is proved the quasi-isometry invariance of
relative hyperbolicity (see Theorem \ref{t4} in this paper), but
without establishing any relation between the peripheral subgroups
(which is impossible to do in full generality, see the discussion
following Theorem \ref{t4}). Theorem \ref{sqir} resolves this
question. Moreover, it advances towards a classification of
relatively hyperbolic groups. By results in
\cite{PapasogluWhyte:ends}, the classification of relatively
hyperbolic groups reduces to the classification of one-ended
relatively hyperbolic groups. Theorem \ref{sqir} points out a
fundamental necessary condition for the quasi-isometry of two
one-ended relatively hyperbolic groups with NRH peripheral
subgroups: that the peripheral subgroups define the same
collection of quasi-isometry classes. Nevertheless the condition
is not sufficient, as can be seen in \cite{Schwartz:RankOne},
where it is proved for instance that two fundamental groups of
finite volume hyperbolic three-manifolds are quasi-isometric if
and only if they are commensurable (while all their peripheral
subgroups are isomorphic to $\Z^2$, when there is no torsion).
This raises the question on what finer invariants of
quasi-isometry may exist for relatively hyperbolic groups (besides
the q.i. classes of peripherals) which would allow advancing
further in the classification.

\medskip

Theorems \ref{qinonrelhyp} and \ref{sqir} motivate the study of
non-relative hyperbolicity and, in particular, thickness. In order to verify
thickness of a finitely generated group, we formulate an algebraic
form of thickness in the setting of groups endowed with word
metrics and their undistorted subgroups (see Definition
\ref{dgthick}). Many important groups turn out to have this
property, and therefore are NRH:

\begin{thm}\label{list}
The following finitely generated groups (keyed to section numbers)
are algebraically thick with respect to the word metric:

\begin{itemize}
\item[\S\ref{SectionMCG}.] $\MCG(S)$, when $S$ is an orientable
finite type surface with
\\
\hphantom\qquad $3\cdot \, \mathrm{genus}(S) + \#\,
\mathrm{punctures} \geq 5;$
\item[\S\ref{section:outfn}.] $\autfn$ and $\outfn$, when $n\geq 3$;
\item[\S\ref{section:artin}.] A freely indecomposable Artin group
with any of the following properties: the integer labels on the Artin
presentation graph are all even; the Artin presentation graph is a
tree; the Artin presentation graph has no triangles; the associated
Coxeter group is finite or affine of type $\wt A_n$.
\item[\S\ref{Sectiongrm}.] Fundamental groups of 3-dimensional graph
manifolds;
\item[\S \ref{lattices}.] Non-uniform lattices in semisimple groups
of rank at least two.
\end{itemize}
\end{thm}
The failure of strong relative hyperbolicity for $\SLNZ$ when $n
\ge 3$ was first proved in \cite{KarlssonNoskov}. For the case of
mapping class groups, the failure of strong relative hyperbolicity
is also proved in \cite{AASh:RelHypMCG},
\cite{Bowditch:3manifolds}, and \cite{KarlssonNoskov}; see the
discussion after Corollary~\ref{mcgnobcp}. If one is solely
interested in disproving strong relative hyperbolicity, there are
more direct approaches which avoid asymptotic cones, such as the
one taken in \cite{AASh:RelHypMCG}. In Propositions~\ref{conjb}
and \ref{conjb2} we also give such results, generalizing the main
theorem of \cite{AASh:RelHypMCG}.

\medskip

In the particular case of Artin groups, more can be proved
concerning relative hyperbolicity. The following is an immediate
consequence of Proposition~\ref{conjb2} and
Example~\ref{twogeneratorartin}:

\begin{prop}\label{allArtin} Except for the integers, any
Artin group with connected Artin presentation graph is not
relatively hyperbolic.
\end{prop}

Note that this gives a complete classification of which Artin
groups are relatively hyperbolic, since any group with a
disconnected presentation graph is freely decomposable and hence
relatively hyperbolic with respect to the factors in the free
decomposition.

\begin{rmk}
For the Artin groups which are not in the list of Theorem
\ref{list} we do not know whether they are thick or not. Possibly
some of them might turn out to be examples of NRH groups that are
not thick. This would provide a nice class of examples, as the
groups we know which are NRH, but not thick, are fairly
pathological, cf. the end of Section~\ref{section:thick}.
\end{rmk}

\medskip

Theorem \ref{list} and Proposition \ref{allArtin} are interesting
also because some of the listed groups are known to be weakly
relatively hyperbolic. Examples include: mapping class groups
\cite{MasurMinsky:complex1}, certain Artin groups
\cite{KapovichSchupp:Artin}, and fundamental groups of graph
manifolds. Thus the study we begin in this paper, of thick groups
from the point of view of quasi-isometric rigidity, may also be
perceived as a first attempt to study quasi-isometric rigidity of
weakly relatively hyperbolic groups. Note that up to now there is
no general result on the quasi-isometric behavior of weakly
relatively hyperbolic groups. In \cite{KapovichLeeb:Haken},
\cite{Papasoglu:Zsplittings}, \cite{DunwoodySageev:JSJ},
\cite{MSW:QTOne} and \cite{MSW:QTTwo} strong quasi-isometric
rigidity results are proved for some particular cases of weakly
relatively hyperbolic groups---in fact all of them are fundamental
groups of some graphs of groups (fundamental groups of Haken
manifolds, groups with a JSJ decomposition, fundamental groups of
finite graphs of groups with Bass-Serre tree of finite depth).

Some of the groups mentioned in Theorem~\ref{list} present even
further similarities with (strongly) relatively hyperbolic groups,
in that all their asymptotic cones are tree-graded metric spaces.
This is the case for the mapping class groups, where it was proved
by Behrstock \cite{Behrstock:asymptotic}; and for the fundamental
groups of 3-dimensional graph manifolds, where it follows from
results in \cite{KapovichLeeb:3manifolds} and \cite{KKL:QI}; the
latter class includes right angled Artin groups whose Artin
presentation graph is a tree of diameter at least three (see
Proposition~\ref{raagcutpoints}).

In particular these examples answer in the negative two questions
of Dru\c{t}u and Sapir (see
\cite[Problem~1.18]{DrutuSapir:TreeGraded}) regarding a finitely
generated group $G$ for which every asymptotic cone is
tree-graded: Is $G$ relatively hyperbolic? And is $G$
asymptotically tree-graded with respect to some collection of
subsets of $G$? The negative answers to these questions indicate
that a supplementary condition on the pieces in the asymptotic
cones is indeed necessary.

Another question resolved by the example of mapping class groups
is whether every relatively hyperbolic group is in fact hyperbolic
relative to subgroups that are \emph{unconstricted} (see
\cite[Problem~1.21]{DrutuSapir:TreeGraded}). Indeed, consider the
finitely presented
relatively hyperbolic group $\Gamma=\MCG(S)*\MCG(S)$. Suppose that it is
hyperbolic relative to a finite collection of unconstricted
peripheral subgroups $\H$. Corollary~\ref{c2} implies that each
$H\in\H$ must be contained in a conjugate $\gamma \MCG(S)
\gamma^{-1}$ of one of the two free factors isomorphic to
$\MCG(S)$ in $\Gamma$. Applying Corollary~\ref{c2} again to
$\Gamma$ seen as hyperbolic relative to the subgroups in $\H$ we
obtain that $\gamma \MCG(S) \gamma^{-1}$ is contained in a
conjugate of a subgroup $H_1\in \H$. This implies that $H$ is
contained in a conjugate of $H_1$, a situation which can occur
only if $H$ coincides with the conjugate of $H_1$. Thus the two
inclusions above are equalities, in particular $H=\gamma \MCG(S)
\gamma^{-1}$. On the other hand, all asymptotic cones of $\MCG(S)$
have (global) cut-points, and hence the same holds for $\gamma
\MCG(S) \gamma^{-1}$ (see \cite{Behrstock:asymptotic}); this
contradicts the hypothesis that $H$ is unconstricted. Note that in
the previous argument $\MCG(S)$ can be replaced by any group which
is thick (or more generally not relatively hyperbolic) and with
(global) cut-points in any asymptotic cone (i.e.,
\textit{constricted}, in the terminology of
\cite{DrutuSapir:TreeGraded}).

In Section~\ref{section:Dunwoody}, we answer a related weaker
question, namely, does any relatively hyperbolic group admit a
family of peripheral subgroups which are not relatively
hyperbolic? The answer is no, with Dunwoody's inaccessible group
providing a counterexample. Since finitely presented groups are
accessible, this raises the following natural question.

\begin{qn}
Is there any example of a finitely presented relatively hyperbolic
group such that every list of peripheral subgroups contains
a relatively hyperbolic group?

A similar question can be asked for groups without torsion, as
these groups are likewise accessible.
\end{qn}

Thickness can be studied for spaces other than groups. As an example
of this we prove the following:

\Restate{Theorem}{WPnotatg}{}{For any surface $S$ with $3\cdot \,
\mathrm{genus}(S) + \#\, \mathrm{punctures} \geq 9$, the
Teichm\"{u}ller space with the Weil-Petersson metric is thick.} In
particular the Teichm\"{u}ller space is not asymptotically
tree-graded. An interesting aspect of this theorem is that
although these higher complexity Teichm\"{u}ller spaces are not
asymptotically tree-graded, they do have tree-graded asymptotic
cones as proven in \cite{Behrstock:asymptotic}. We also note that
the lack of relative hyperbolicity contrasts with the cases with
$3\cdot \, \mathrm{genus}(S) + \#\, \mathrm{punctures} \le 5$
where it has been shown that Teichm\"{u}ller space is
$\delta$--hyperbolic with the Weil-Petersson metric (see
\cite{BrockFarb:curvature}, and also \cite{Aramayona:thesis},
\cite{Behrstock:asymptotic}). It also contrasts with the relative
hyperbolicity of Teichm\"{u}ller space in the cases where $3\cdot
\, \mathrm{genus}(S) + \#\, \mathrm{punctures} = 6$, as recently
shown in \cite{BrockMasur:WPrelhyp}.

\bigskip

The paper is organized as follows. Section~\ref{prel} provides
background on asymptotic cones and various tools developed in
\cite{DrutuSapir:TreeGraded} for studying relatively hyperbolic
groups. In Section~\ref{sunc} we discuss the property of (not)
having cut-points in asymptotic cones.

Section~\ref{snrh} contains some general results regarding
quasi-isometric embeddings of NRH groups into relatively
hyperbolic groups and our main theorem of rigidity of relatively
hyperbolic groups. Motivated by these results we provide examples
of NRH groups, and in Section~\ref{section:network} we describe a
way to build NRH groups. In Section~\ref{section:Dunwoody} we
discuss an example of a relatively hyperbolic group such that any
list of peripheral subgroups contains a relatively hyperbolic
group.

In Section~\ref{section:thick} we define metric and algebraic
thickness, we provide results on the structure and rigidity of
thick spaces and groups and we discuss an example of an NRH group
which is not thick.

The remaining sections of this work establish thickness for
various groups and metric spaces. For the mapping class groups,
the automorphism group of a free group, and the outer automorphism
group of a free group we prove thickness in all cases except when
these groups are virtually free (and hence are not thick), this is
done in Sections~\ref{SectionMCG} and~\ref{section:outfn}. Artin
groups are studied in Section~\ref{section:artin}. Graph manifolds
and Teichm\"{u}ller space are shown to be thick in
Sections~\ref{Sectiongrm} and~\ref{Section:teich}. Finally in
Section~\ref{lattices}, we establish thickness for non-uniform
lattices (thickness in the uniform case follows from
\cite{KleinerLeeb:buildings}).

\subsection*{Acknowledgements}

The authors would like to thank Ruth Charney for sharing her
knowledge of Artin groups with us by answering our many questions.
We would also like to thank Michael Handel, Ilya Kapovich, Bruce
Kleiner, Walter Neumann, Mark Sapir and Kevin Wortman for helpful
conversations about aspects of this work. Also, we thank the 
anonymous referee for several useful comments.

The first named author
would like to acknowledge the support of the Barnard
College/Columbia University Mathematics Department. 
He would also like to thank the June
2005 Asymptotic and Probabilistic Methods in Geometric Group
Theory Conference, where he presented results from this paper. The
second named author thanks IHES and MPIM in Bonn for their
hospitality during the summer of 2005, when part of this work was
done. The work of the second author has been supported by the ANR 
grant GGPG.

\section{Preliminaries}\label{prel}

A \emph{non-principal ultrafilter} on the positive integers,
denoted by $\omega$, is a nonempty collection of sets of positive
integers with the following properties:

\begin{enumerate}
\item If $S_{1}\in\omega$ and $S_{2}\in\omega$, then
$S_{1}\cap S_{2}\in\omega$.
\item If $S_{1}\subset S_{2}$ and $S_{1}\in\omega$, then
$S_{2}\in\omega$.
\item For each $S\subset \N$ exactly one of the following must
occur:
$S\in\omega$ or $\N\setminus S\in\omega$.
\item \label{ItemCofinite}
$\omega$ does not contain any finite set.
\end{enumerate}

\medskip

\noindent\textit{Convention:} The adjective ``non-principal''
refers to item~(\ref{ItemCofinite}). Since we work only with
non-principal ultrafilters, we shall tacitly drop this adjective
throughout the sequel.

\medskip

\noindent For an ultrafilter $\omega$, a topological space $X$,
and a sequence of points $(x_{i})_{i\in\N}$ in $X$, we define $x$
to be the \emph{ultralimit of $(x_{i})_{i\in\N}$ with respect to
$\omega$}, and we write $x=\ulimi x_{i}$, if and only if for any
neighborhood $\nn$ of $x$ in $X$ the set $\{i\in\N:x_i \in \nn\}$
is in $\omega$. Note that when $X$ is compact any sequence in $X$
has an ultralimit \cite{Bourbaki:topologie}. If moreover $X$ is
Hausdorff
then the ultralimit of any sequence is unique.
Fix an ultrafilter $\omega$ and a family of based metric spaces
$(X_{i}, x_{i}, \dist_{i})$. Using the ultrafilter, a
pseudo-distance on $\prod_{i\in\N}X_{i}$ is provided by:
$$\dist_{\omega}((a_{i}),(b_{i}))=\ulimi
\dist_{i}(a_{i},b_{i})\in [0,\infty]\, .$$
One can eliminate the possibility of the previous pseudo-distance
taking the value $+\infty$ by restricting to sequences $y=(y_i)$
such that $\dist_{\omega}(y,x)<\infty$, where $x=(x_{i})$. A
metric space can be then defined, called the \emph{ultralimit of
$(X_{i},x_{i}, \dist_{i})$}\index{$\ulim$}, by:
$$\ulimi (X_{i},x_{i}, \dist_{i}) =\left\{ y\in\prod_{i\in\N}X_{i}:
\dist_{\omega}(y,x)<\infty \right\} / \sim\, ,$$
\noindent where
for two points $y,z\in \prod_{i\in\N}X_{i}$ we define $y\sim z$ if
and only if $\dist_{\omega}(y,z)=0$. The pseudo-distance on
$\prod_{i \in \N} X_i$ induces a complete metric on $\ulimi
(X_i,x_i,\dist_i)$.

Let now $(X,\dist )$ be a metric space. Consider $x=(x_n)$ a
sequence of points in $X$, called \textit{observation points}, and
$d=(d_n)$ a sequence of positive numbers such that $\lim_\omega
d_n=+\infty$, called \textit{scaling constants}. First defined in
\cite{Gromov:PolynomialGrowth} and \cite{DriesWilkie}, the
\emph{asymptotic cone of $(X,\dist)$ relative to the ultrafilter
$\omega$ and the sequences $x$ and $d$} is given by:
$$\coneomega (X,x,d)=\ulimi
\left( X,x_{n},\frac{1}{d_n}\dist \right).$$

When the group of isometries of $X$ acts on $X$ so that all orbits
intersect a fixed bounded set, the asymptotic cone is independent
of the choice of observation points. An important example of this
is when $X$ is a finitely generated group with a word metric;
thus, when $X$ is a finitely generated group we always take the
observation points to be the constant sequence $(1)$ and we drop
the observation point from our notation.

Every sequence $(A_n)_{n\in \N}$ of non-empty subsets of $X$ has a
\textit{limit set} in the asymptotic cone $\coneomega (X,x,d)$,
denoted by $\ulim A_n$ and defined as the set of images in the
asymptotic cone of sequences $(a_n)_{n\in \N}$ with $a_n \in A_n$
for every~$n$. The set $\ulim A_n$ is empty when $\ulim
\frac{\dist (x_n , A_n)}{d_n}=\infty $, otherwise it is a closed
subset of~$\coneomega (X,x,d)$. In the latter case, $\ulim A_n$ is
isometric to the ultralimit of $(A_n,y_n,\frac{\dist}{d_n})_{n\in
\N}$ with the metric $\dist$ on $A_n$ induced from $X$, and with
basepoints $y_n \in A_n$ such that
$\ulim\frac{\dist(x_n,y_n)}{d_n}<\infty\, $.

Given a collection $\P$ of subsets in $X$ and an asymptotic cone
$\coneomega (X,x,d)$ of~$X$, we denote by $\ocone(\P)$ the
collection of non-empty limit sets $\ulim P_n$ where
$(P_{n})_{n\in\N}$ is a sequence of subsets $P_n\in\P$. We will
often consider the case where $X=G$ is a group and $\H$ is a fixed
collection of subgroups of $G$, in this case we take $\P$ to be
the collection of left cosets $gH$, with $g\in G$ and $H\in\H$. We
denote the latter collection also by $\mathcal{L}\H$. We now
recall a notion introduced in \cite[\S 2]{DrutuSapir:TreeGraded}.

\begin{defn} Let ${\mathbb F}$ be a complete geodesic metric space and
let $\mathcal{P}$ be a collection of closed geodesic subsets
(called \emph{pieces}).
The space ${\mathbb F}$ is said to be \emph{tree-graded
with respect to $\mathcal{P}$}
when the following two properties are satisfied:

\begin{itemize}
\item[$(T_1)$] The intersection of each pair of distinct pieces has at
most one point.
\item[$(T_2)$] Every simple non-trivial geodesic triangle in ${\mathbb F}$
is contained in one piece.
\end{itemize}
When the collection of pieces $\mathcal{P}$ is understood then we say
simply that ${\mathbb F}$ is \emph{tree-graded}.
\end{defn}

\begin{lem}[Dru\c{t}u--Sapir  \cite{DrutuSapir:TreeGraded}]\label{tx}
Let ${\mathbb F}$ be a complete geodesic metric
space which is tree-graded with respect to a collection of pieces
$\mathcal{P}$.

\begin{itemize}
\item[(1)] For every point $x\in \free$, the set $T_x$ of topological
arcs originating at $x$ and intersecting any piece in at most one point
is a complete real tree (possibly reduced to a point). Moreover if
$y\in T_x$ then $T_y=T_x$.
\item[(2)] Any topological arc joining two points in a piece is
contained in the same piece. Any topological arc joining two points
in a tree $T_x$ is contained in the same tree $T_x$.
\end{itemize}
\end{lem}
A tree as in Lemma \ref{tx}~(1) is called a \textit{transversal
tree}, and a geodesic in it is called a \textit{transversal
geodesic}. Both of these notions are defined relative to the
collection of pieces $\mathcal{P}$, which when understood is
suppressed.

The notion of tree-graded metric space is related to the existence
of cut-points.

\medskip

\textit{Convention:} By cut-points we always mean \emph{global}
cut-points. We consider a singleton to have a cut-point.

\medskip

\begin{lem}[Dru\c{t}u--Sapir \cite{DrutuSapir:TreeGraded},
    Lemma 2.31]\label{cutting}
Let $X$ be a complete geodesic metric space containing at least
two points and let $\calc$ be a non-empty set of cut-points in
$X$. There exists a uniquely defined (maximal in an appropriate
sense) collection $\calp$ of subsets of $X$ such that
\begin{itemize}
\item $X$ is tree-graded with respect to $\pp$;
\item any piece in $\pp$ is either a singleton or a set
with no cut-point in~$\calc$.
\end{itemize}
Moreover the intersection of any two distinct pieces from $\calp$
is either empty or a point from $\calc$.
\end{lem}

\begin{defn}
Let $X$ be a metric space and let $\aaa$ be a collection of
subsets in $X$. We say that $X$ is \emph{asymptotically
tree-graded (ATG)} with respect to $\aaa$ if
\begin{itemize}
\item[(I)] every asymptotic cone $\coneomega(X)$ of $X$ is
tree-graded with respect to $\ocone(\aaa)$;
\item[(II)] $X$ is not contained in a finite radius neighborhood
of any of the subsets in~$\aaa$.
\end{itemize}
\end{defn}

The subsets in $\aaa$ are called \emph{peripheral subsets}.

The second condition does not appear in
\cite{DrutuSapir:TreeGraded}. It is added here to avoid the
trivial cases, like that of $X$ asymptotically tree-graded with
respect to $\aaa =\{X\}$. For emphasis, one could refer to an ATG
structure satisfying (II)  as being a \emph{proper asymptotically
tree-graded structure}. Since we always assume that the tubular
neighborhoods of peripheral subsets are proper subsets (see
Convention \ref{convnrh}), we suppress the use of the adjective
``proper.'' Similarly, we assume that relative hyperbolicity is
always with respect to a collection of proper peripheral
subgroups.

\medskip

As mentioned in the introduction, Dru\c{t}u--Sapir provide a
characterization of ATG metric spaces, further simplified by
Dru\c{t}u in \cite{Drutu:RelHyp}, in terms of three metric
properties involving elements of
$\aaa$, geodesics, and geodesic polygons.
There are several versions of the list of three
properties, we recall here those that we shall use most, keeping
the notation in \cite{Drutu:RelHyp}.

First we recall the notion of fat polygon introduced in
\cite{DrutuSapir:TreeGraded}. This notion is in some sense the
opposite of the notion of ``thin'' polygon (i.e., a polygon behaving
metrically like a polygon in a tree, up to bounded perturbation).

\medskip

Throughout the paper $\nn_r(A)$ denotes the set of points $x$
satisfying $\dist(x,A)<r$ and $\overline{\nn}_r(A)$ the set of
points $x$ with $\dist(x,A)\leq r$.

\begin{n}\label{breve}
For every quasi-geodesic $\fp$ in a metric space $X$, we denote
the origin of $\fp$ by $\fp_-$ and the endpoint of $\fp $ by
$\fp_+$.

Given $r>0$ we denote by $\breve{\fp}_r$ the set $\fp \setminus
\nn_r \left( \left\{ \fp_-\, ,\, \fp_+ \right\} \right)$.
\end{n}

\medskip

A \emph{geodesic (quasi-geodesic) $k$-gonal line} is a set $P$ which is the
 union of $k$ geodesics (quasi-geodesics) $\q_1,...,\q_k$
such that $(\q_i)_+=(\q_{i+1})_-$ for $i=1,...,k-1$. If moreover
$(\q_k)_+=(\q_{1})_-$ then we say that $P$ is a \emph{geodesic
(quasi-geodesic) $k$-gon}.

\medskip

\begin{n}\label{ox}
Given a vertex $x\in \mathcal{V}$ and $\q , \q'$ the consecutive
edges of $P$ such that $x=\q_+ =\q'_-$, we denote the polygonal
line $P\setminus ( \q\cup \q')$ by $\oo_{x}(P)$. When there is no
possibility of confusion we simply denote it by $\oo_{x}$.
\end{n}

\bigskip

\unitlength .3mm 
\linethickness{0.4pt}
\ifx\plotpoint\undefined\newsavebox{\plotpoint}\fi 
\begin{picture}(248.81,123.52)(0,0)
\multiput(47.32,42.92)(21.75,-.03125){8}{\line(1,0){21.75}}
\put(74.33,42.92){\line(0,1){1.138}}
\put(74.31,44.06){\line(0,1){1.136}}
\put(74.24,45.19){\line(0,1){1.132}}
\multiput(74.12,46.32)(-.03313,.22514){5}{\line(0,1){.22514}}
\multiput(73.95,47.45)(-.03035,.15969){7}{\line(0,1){.15969}}
\multiput(73.74,48.57)(-.03237,.1385){8}{\line(0,1){.1385}}
\multiput(73.48,49.68)(-.0305,.10962){10}{\line(0,1){.10962}}
\multiput(73.17,50.77)(-.03186,.09841){11}{\line(0,1){.09841}}
\multiput(72.82,51.86)(-.03294,.08891){12}{\line(0,1){.08891}}
\multiput(72.43,52.92)(-.031392,.074967){14}{\line(0,1){.074967}}
\multiput(71.99,53.97)(-.032193,.068685){15}{\line(0,1){.068685}}
\multiput(71.51,55)(-.032842,.063077){16}{\line(0,1){.063077}}
\multiput(70.98,56.01)(-.033361,.058025){17}{\line(0,1){.058025}}
\multiput(70.41,57)(-.03199,.050626){19}{\line(0,1){.050626}}
\multiput(69.81,57.96)(-.032371,.046785){20}{\line(0,1){.046785}}
\multiput(69.16,58.9)(-.032662,.043232){21}{\line(0,1){.043232}}
\multiput(68.47,59.8)(-.032872,.039929){22}{\line(0,1){.039929}}
\multiput(67.75,60.68)(-.033009,.036848){23}{\line(0,1){.036848}}
\multiput(66.99,61.53)(-.03308,.033962){24}{\line(0,1){.033962}}
\multiput(66.2,62.34)(-.034469,.032552){24}{\line(-1,0){.034469}}
\multiput(65.37,63.13)(-.037353,.032436){23}{\line(-1,0){.037353}}
\multiput(64.51,63.87)(-.040432,.032252){22}{\line(-1,0){.040432}}
\multiput(63.62,64.58)(-.045917,.03359){20}{\line(-1,0){.045917}}
\multiput(62.7,65.25)(-.049767,.03331){19}{\line(-1,0){.049767}}
\multiput(61.76,65.89)(-.053954,.032938){18}{\line(-1,0){.053954}}
\multiput(60.79,66.48)(-.058533,.032461){17}{\line(-1,0){.058533}}
\multiput(59.79,67.03)(-.063577,.031864){16}{\line(-1,0){.063577}}
\multiput(58.77,67.54)(-.074115,.033352){14}{\line(-1,0){.074115}}
\multiput(57.74,68.01)(-.081246,.032556){13}{\line(-1,0){.081246}}
\multiput(56.68,68.43)(-.08941,.03157){12}{\line(-1,0){.08941}}
\multiput(55.61,68.81)(-.10878,.03337){10}{\line(-1,0){.10878}}
\multiput(54.52,69.14)(-.12231,.032){9}{\line(-1,0){.12231}}
\multiput(53.42,69.43)(-.13898,.03023){8}{\line(-1,0){.13898}}
\multiput(52.31,69.67)(-.18683,.03253){6}{\line(-1,0){.18683}}
\multiput(51.18,69.87)(-.22563,.02965){5}{\line(-1,0){.22563}}
\put(50.06,70.02){\line(-1,0){1.133}}
\put(48.92,70.12){\line(-1,0){1.137}}
\put(47.79,70.17){\line(-1,0){1.138}}
\put(46.65,70.18){\line(-1,0){1.137}}
\put(45.51,70.14){\line(-1,0){1.134}}
\multiput(44.38,70.05)(-.22594,-.0272){5}{\line(-1,0){.22594}}
\multiput(43.25,69.91)(-.18717,-.0305){6}{\line(-1,0){.18717}}
\multiput(42.12,69.73)(-.1592,-.03282){7}{\line(-1,0){.1592}}
\multiput(41.01,69.5)(-.12265,-.03067){9}{\line(-1,0){.12265}}
\multiput(39.91,69.22)(-.10914,-.03218){10}{\line(-1,0){.10914}}
\multiput(38.82,68.9)(-.09791,-.03337){11}{\line(-1,0){.09791}}
\multiput(37.74,68.53)(-.081594,-.031672){13}{\line(-1,0){.081594}}
\multiput(36.68,68.12)(-.074473,-.032546){14}{\line(-1,0){.074473}}
\multiput(35.63,67.67)(-.06818,-.03325){15}{\line(-1,0){.06818}}
\multiput(34.61,67.17)(-.058882,-.031824){17}{\line(-1,0){.058882}}
\multiput(33.61,66.63)(-.054308,-.03235){18}{\line(-1,0){.054308}}
\multiput(32.63,66.04)(-.050126,-.032768){19}{\line(-1,0){.050126}}
\multiput(31.68,65.42)(-.046279,-.03309){20}{\line(-1,0){.046279}}
\multiput(30.76,64.76)(-.042722,-.033326){21}{\line(-1,0){.042722}}
\multiput(29.86,64.06)(-.039417,-.033485){22}{\line(-1,0){.039417}}
\multiput(28.99,63.32)(-.036334,-.033575){23}{\line(-1,0){.036334}}
\multiput(28.16,62.55)(-.033447,-.033601){24}{\line(0,-1){.033601}}
\multiput(27.35,61.74)(-.033407,-.036488){23}{\line(0,-1){.036488}}
\multiput(26.58,60.9)(-.033304,-.03957){22}{\line(0,-1){.03957}}
\multiput(25.85,60.03)(-.033129,-.042874){21}{\line(0,-1){.042874}}
\multiput(25.16,59.13)(-.032877,-.046431){20}{\line(0,-1){.046431}}
\multiput(24.5,58.21)(-.032537,-.050276){19}{\line(0,-1){.050276}}
\multiput(23.88,57.25)(-.0321,-.054456){18}{\line(0,-1){.054456}}
\multiput(23.3,56.27)(-.033525,-.062717){16}{\line(0,-1){.062717}}
\multiput(22.77,55.27)(-.032937,-.068332){15}{\line(0,-1){.068332}}
\multiput(22.27,54.24)(-.032204,-.074622){14}{\line(0,-1){.074622}}
\multiput(21.82,53.2)(-.031297,-.081739){13}{\line(0,-1){.081739}}
\multiput(21.41,52.13)(-.03292,-.09806){11}{\line(0,-1){.09806}}
\multiput(21.05,51.06)(-.03168,-.10928){10}{\line(0,-1){.10928}}
\multiput(20.74,49.96)(-.03011,-.12279){9}{\line(0,-1){.12279}}
\multiput(20.46,48.86)(-.03209,-.15935){7}{\line(0,-1){.15935}}
\multiput(20.24,47.74)(-.02964,-.18731){6}{\line(0,-1){.18731}}
\multiput(20.06,46.62)(-.0327,-.2826){4}{\line(0,-1){.2826}}
\put(19.93,45.49){\line(0,-1){1.135}}
\put(19.85,44.35){\line(0,-1){1.137}}
\put(19.81,43.22){\line(0,-1){1.138}}
\put(19.82,42.08){\line(0,-1){1.136}}
\put(19.88,40.94){\line(0,-1){1.133}}
\multiput(19.99,39.81)(.03068,-.22549){5}{\line(0,-1){.22549}}
\multiput(20.14,38.68)(.03339,-.18668){6}{\line(0,-1){.18668}}
\multiput(20.34,37.56)(.03086,-.13884){8}{\line(0,-1){.13884}}
\multiput(20.59,36.45)(.03256,-.12216){9}{\line(0,-1){.12216}}
\multiput(20.88,35.35)(.03079,-.09875){11}{\line(0,-1){.09875}}
\multiput(21.22,34.26)(.03198,-.08927){12}{\line(0,-1){.08927}}
\multiput(21.6,33.19)(.032929,-.081096){13}{\line(0,-1){.081096}}
\multiput(22.03,32.14)(.033692,-.073961){14}{\line(0,-1){.073961}}
\multiput(22.5,31.1)(.032156,-.06343){16}{\line(0,-1){.06343}}
\multiput(23.02,30.09)(.032729,-.058384){17}{\line(0,-1){.058384}}
\multiput(23.57,29.1)(.033185,-.053802){18}{\line(0,-1){.053802}}
\multiput(24.17,28.13)(.033538,-.049614){19}{\line(0,-1){.049614}}
\multiput(24.81,27.19)(.032191,-.043583){21}{\line(0,-1){.043583}}
\multiput(25.49,26.27)(.032437,-.040284){22}{\line(0,-1){.040284}}
\multiput(26.2,25.38)(.032608,-.037204){23}{\line(0,-1){.037204}}
\multiput(26.95,24.53)(.03271,-.034319){24}{\line(0,-1){.034319}}
\multiput(27.73,23.7)(.034113,-.032924){24}{\line(1,0){.034113}}
\multiput(28.55,22.91)(.036999,-.03284){23}{\line(1,0){.036999}}
\multiput(29.4,22.16)(.04008,-.032689){22}{\line(1,0){.04008}}
\multiput(30.29,21.44)(.043381,-.032463){21}{\line(1,0){.043381}}
\multiput(31.2,20.76)(.046933,-.032156){20}{\line(1,0){.046933}}
\multiput(32.14,20.12)(.053593,-.033521){18}{\line(1,0){.053593}}
\multiput(33.1,19.51)(.058177,-.033094){17}{\line(1,0){.058177}}
\multiput(34.09,18.95)(.063227,-.032552){16}{\line(1,0){.063227}}
\multiput(35.1,18.43)(.068832,-.031878){15}{\line(1,0){.068832}}
\multiput(36.13,17.95)(.080888,-.033436){13}{\line(1,0){.080888}}
\multiput(37.18,17.52)(.08906,-.03253){12}{\line(1,0){.08906}}
\multiput(38.25,17.13)(.09856,-.03141){11}{\line(1,0){.09856}}
\multiput(39.34,16.78)(.12196,-.03332){9}{\line(1,0){.12196}}
\multiput(40.43,16.48)(.13864,-.03173){8}{\line(1,0){.13864}}
\multiput(41.54,16.23)(.15983,-.02962){7}{\line(1,0){.15983}}
\multiput(42.66,16.02)(.22529,-.0321){5}{\line(1,0){.22529}}
\put(43.79,15.86){\line(1,0){1.132}}
\put(44.92,15.74){\line(1,0){1.136}}
\put(46.06,15.68){\line(1,0){1.138}}
\put(47.2,15.66){\line(1,0){1.137}}
\put(48.33,15.69){\line(1,0){1.135}}
\multiput(49.47,15.77)(.2828,.0309){4}{\line(1,0){.2828}}
\multiput(50.6,15.89)(.18749,.02847){6}{\line(1,0){.18749}}
\multiput(51.72,16.06)(.15955,.03109){7}{\line(1,0){.15955}}
\multiput(52.84,16.28)(.13835,.033){8}{\line(1,0){.13835}}
\multiput(53.95,16.54)(.10948,.031){10}{\line(1,0){.10948}}
\multiput(55.04,16.85)(.09826,.03231){11}{\line(1,0){.09826}}
\multiput(56.12,17.21)(.08876,.03335){12}{\line(1,0){.08876}}
\multiput(57.19,17.61)(.074822,.031736){14}{\line(1,0){.074822}}
\multiput(58.24,18.05)(.068537,.032509){15}{\line(1,0){.068537}}
\multiput(59.26,18.54)(.062925,.033132){16}{\line(1,0){.062925}}
\multiput(60.27,19.07)(.057871,.033627){17}{\line(1,0){.057871}}
\multiput(61.25,19.64)(.050479,.032222){19}{\line(1,0){.050479}}
\multiput(62.21,20.25)(.046636,.032585){20}{\line(1,0){.046636}}
\multiput(63.15,20.91)(.043081,.03286){21}{\line(1,0){.043081}}
\multiput(64.05,21.6)(.039778,.033055){22}{\line(1,0){.039778}}
\multiput(64.93,22.32)(.036696,.033178){23}{\line(1,0){.036696}}
\multiput(65.77,23.09)(.033809,.033236){24}{\line(1,0){.033809}}
\multiput(66.58,23.88)(.032393,.034618){24}{\line(0,1){.034618}}
\multiput(67.36,24.71)(.033731,.039207){22}{\line(0,1){.039207}}
\multiput(68.1,25.58)(.033593,.042512){21}{\line(0,1){.042512}}
\multiput(68.81,26.47)(.033379,.046071){20}{\line(0,1){.046071}}
\multiput(69.47,27.39)(.033081,.04992){19}{\line(0,1){.04992}}
\multiput(70.1,28.34)(.03269,.054105){18}{\line(0,1){.054105}}
\multiput(70.69,29.31)(.032192,.058682){17}{\line(0,1){.058682}}
\multiput(71.24,30.31)(.033677,.06797){15}{\line(0,1){.06797}}
\multiput(71.74,31.33)(.033012,.074268){14}{\line(0,1){.074268}}
\multiput(72.21,32.37)(.032182,.081395){13}{\line(0,1){.081395}}
\multiput(72.62,33.43)(.03115,.08956){12}{\line(0,1){.08956}}
\multiput(73,34.5)(.03287,.10893){10}{\line(0,1){.10893}}
\multiput(73.33,35.59)(.03144,.12246){9}{\line(0,1){.12246}}
\multiput(73.61,36.69)(.02959,.13912){8}{\line(0,1){.13912}}
\multiput(73.85,37.81)(.03167,.18698){6}{\line(0,1){.18698}}
\multiput(74.04,38.93)(.02861,.22576){5}{\line(0,1){.22576}}
\put(74.18,40.06){\line(0,1){1.134}}
\put(74.28,41.19){\line(0,1){1.728}}
\put(248.81,42.57){\line(0,1){1.138}}
\put(248.79,43.71){\line(0,1){1.136}}
\put(248.72,44.84){\line(0,1){1.132}}
\multiput(248.6,45.97)(-.03313,.22514){5}{\line(0,1){.22514}}
\multiput(248.43,47.1)(-.03035,.15969){7}{\line(0,1){.15969}}
\multiput(248.22,48.22)(-.03237,.1385){8}{\line(0,1){.1385}}
\multiput(247.96,49.33)(-.0305,.10962){10}{\line(0,1){.10962}}
\multiput(247.65,50.42)(-.03186,.09841){11}{\line(0,1){.09841}}
\multiput(247.3,51.51)(-.03294,.08891){12}{\line(0,1){.08891}}
\multiput(246.91,52.57)(-.031392,.074967){14}{\line(0,1){.074967}}
\multiput(246.47,53.62)(-.032193,.068685){15}{\line(0,1){.068685}}
\multiput(245.99,54.65)(-.032842,.063077){16}{\line(0,1){.063077}}
\multiput(245.46,55.66)(-.033361,.058025){17}{\line(0,1){.058025}}
\multiput(244.89,56.65)(-.03199,.050626){19}{\line(0,1){.050626}}
\multiput(244.29,57.61)(-.032371,.046785){20}{\line(0,1){.046785}}
\multiput(243.64,58.55)(-.032662,.043232){21}{\line(0,1){.043232}}
\multiput(242.95,59.45)(-.032872,.039929){22}{\line(0,1){.039929}}
\multiput(242.23,60.33)(-.033009,.036848){23}{\line(0,1){.036848}}
\multiput(241.47,61.18)(-.03308,.033962){24}{\line(0,1){.033962}}
\multiput(240.68,61.99)(-.034469,.032552){24}{\line(-1,0){.034469}}
\multiput(239.85,62.78)(-.037353,.032436){23}{\line(-1,0){.037353}}
\multiput(238.99,63.52)(-.040432,.032252){22}{\line(-1,0){.040432}}
\multiput(238.1,64.23)(-.045917,.03359){20}{\line(-1,0){.045917}}
\multiput(237.18,64.9)(-.049767,.03331){19}{\line(-1,0){.049767}}
\multiput(236.24,65.54)(-.053954,.032938){18}{\line(-1,0){.053954}}
\multiput(235.27,66.13)(-.058533,.032461){17}{\line(-1,0){.058533}}
\multiput(234.27,66.68)(-.063577,.031864){16}{\line(-1,0){.063577}}
\multiput(233.25,67.19)(-.074115,.033352){14}{\line(-1,0){.074115}}
\multiput(232.22,67.66)(-.081246,.032556){13}{\line(-1,0){.081246}}
\multiput(231.16,68.08)(-.08941,.03157){12}{\line(-1,0){.08941}}
\multiput(230.09,68.46)(-.10878,.03337){10}{\line(-1,0){.10878}}
\multiput(229,68.79)(-.12231,.032){9}{\line(-1,0){.12231}}
\multiput(227.9,69.08)(-.13898,.03023){8}{\line(-1,0){.13898}}
\multiput(226.79,69.32)(-.18683,.03253){6}{\line(-1,0){.18683}}
\multiput(225.66,69.52)(-.22563,.02965){5}{\line(-1,0){.22563}}
\put(224.54,69.67){\line(-1,0){1.133}}
\put(223.4,69.77){\line(-1,0){1.137}}
\put(222.27,69.82){\line(-1,0){1.138}}
\put(221.13,69.83){\line(-1,0){1.137}}
\put(219.99,69.79){\line(-1,0){1.134}}
\multiput(218.86,69.7)(-.22594,-.0272){5}{\line(-1,0){.22594}}
\multiput(217.73,69.56)(-.18717,-.0305){6}{\line(-1,0){.18717}}
\multiput(216.6,69.38)(-.1592,-.03282){7}{\line(-1,0){.1592}}
\multiput(215.49,69.15)(-.12265,-.03067){9}{\line(-1,0){.12265}}
\multiput(214.39,68.87)(-.10914,-.03218){10}{\line(-1,0){.10914}}
\multiput(213.3,68.55)(-.09791,-.03337){11}{\line(-1,0){.09791}}
\multiput(212.22,68.18)(-.081594,-.031672){13}{\line(-1,0){.081594}}
\multiput(211.16,67.77)(-.074473,-.032546){14}{\line(-1,0){.074473}}
\multiput(210.11,67.32)(-.06818,-.03325){15}{\line(-1,0){.06818}}
\multiput(209.09,66.82)(-.058882,-.031824){17}{\line(-1,0){.058882}}
\multiput(208.09,66.28)(-.054308,-.03235){18}{\line(-1,0){.054308}}
\multiput(207.11,65.69)(-.050126,-.032768){19}{\line(-1,0){.050126}}
\multiput(206.16,65.07)(-.046279,-.03309){20}{\line(-1,0){.046279}}
\multiput(205.24,64.41)(-.042722,-.033326){21}{\line(-1,0){.042722}}
\multiput(204.34,63.71)(-.039417,-.033485){22}{\line(-1,0){.039417}}
\multiput(203.47,62.97)(-.036334,-.033575){23}{\line(-1,0){.036334}}
\multiput(202.64,62.2)(-.033447,-.033601){24}{\line(0,-1){.033601}}
\multiput(201.83,61.39)(-.033407,-.036488){23}{\line(0,-1){.036488}}
\multiput(201.06,60.55)(-.033304,-.03957){22}{\line(0,-1){.03957}}
\multiput(200.33,59.68)(-.033129,-.042874){21}{\line(0,-1){.042874}}
\multiput(199.64,58.78)(-.032877,-.046431){20}{\line(0,-1){.046431}}
\multiput(198.98,57.86)(-.032537,-.050276){19}{\line(0,-1){.050276}}
\multiput(198.36,56.9)(-.0321,-.054456){18}{\line(0,-1){.054456}}
\multiput(197.78,55.92)(-.033525,-.062717){16}{\line(0,-1){.062717}}
\multiput(197.25,54.92)(-.032937,-.068332){15}{\line(0,-1){.068332}}
\multiput(196.75,53.89)(-.032204,-.074622){14}{\line(0,-1){.074622}}
\multiput(196.3,52.85)(-.031297,-.081739){13}{\line(0,-1){.081739}}
\multiput(195.89,51.78)(-.03292,-.09806){11}{\line(0,-1){.09806}}
\multiput(195.53,50.71)(-.03168,-.10928){10}{\line(0,-1){.10928}}
\multiput(195.22,49.61)(-.03011,-.12279){9}{\line(0,-1){.12279}}
\multiput(194.94,48.51)(-.03209,-.15935){7}{\line(0,-1){.15935}}
\multiput(194.72,47.39)(-.02964,-.18731){6}{\line(0,-1){.18731}}
\multiput(194.54,46.27)(-.0327,-.2826){4}{\line(0,-1){.2826}}
\put(194.41,45.14){\line(0,-1){1.135}}
\put(194.33,44){\line(0,-1){1.137}}
\put(194.29,42.87){\line(0,-1){1.138}}
\put(194.3,41.73){\line(0,-1){1.136}}
\put(194.36,40.59){\line(0,-1){1.133}}
\multiput(194.47,39.46)(.03068,-.22549){5}{\line(0,-1){.22549}}
\multiput(194.62,38.33)(.03339,-.18668){6}{\line(0,-1){.18668}}
\multiput(194.82,37.21)(.03086,-.13884){8}{\line(0,-1){.13884}}
\multiput(195.07,36.1)(.03256,-.12216){9}{\line(0,-1){.12216}}
\multiput(195.36,35)(.03079,-.09875){11}{\line(0,-1){.09875}}
\multiput(195.7,33.91)(.03198,-.08927){12}{\line(0,-1){.08927}}
\multiput(196.08,32.84)(.032929,-.081096){13}{\line(0,-1){.081096}}
\multiput(196.51,31.79)(.033692,-.073961){14}{\line(0,-1){.073961}}
\multiput(196.98,30.75)(.032156,-.06343){16}{\line(0,-1){.06343}}
\multiput(197.5,29.74)(.032729,-.058384){17}{\line(0,-1){.058384}}
\multiput(198.05,28.75)(.033185,-.053802){18}{\line(0,-1){.053802}}
\multiput(198.65,27.78)(.033538,-.049614){19}{\line(0,-1){.049614}}
\multiput(199.29,26.84)(.032191,-.043583){21}{\line(0,-1){.043583}}
\multiput(199.97,25.92)(.032437,-.040284){22}{\line(0,-1){.040284}}
\multiput(200.68,25.03)(.032608,-.037204){23}{\line(0,-1){.037204}}
\multiput(201.43,24.18)(.03271,-.034319){24}{\line(0,-1){.034319}}
\multiput(202.21,23.35)(.034113,-.032924){24}{\line(1,0){.034113}}
\multiput(203.03,22.56)(.036999,-.03284){23}{\line(1,0){.036999}}
\multiput(203.88,21.81)(.04008,-.032689){22}{\line(1,0){.04008}}
\multiput(204.77,21.09)(.043381,-.032463){21}{\line(1,0){.043381}}
\multiput(205.68,20.41)(.046933,-.032156){20}{\line(1,0){.046933}}
\multiput(206.62,19.77)(.053593,-.033521){18}{\line(1,0){.053593}}
\multiput(207.58,19.16)(.058177,-.033094){17}{\line(1,0){.058177}}
\multiput(208.57,18.6)(.063227,-.032552){16}{\line(1,0){.063227}}
\multiput(209.58,18.08)(.068832,-.031878){15}{\line(1,0){.068832}}
\multiput(210.61,17.6)(.080888,-.033436){13}{\line(1,0){.080888}}
\multiput(211.66,17.17)(.08906,-.03253){12}{\line(1,0){.08906}}
\multiput(212.73,16.78)(.09856,-.03141){11}{\line(1,0){.09856}}
\multiput(213.82,16.43)(.12196,-.03332){9}{\line(1,0){.12196}}
\multiput(214.91,16.13)(.13864,-.03173){8}{\line(1,0){.13864}}
\multiput(216.02,15.88)(.15983,-.02962){7}{\line(1,0){.15983}}
\multiput(217.14,15.67)(.22529,-.0321){5}{\line(1,0){.22529}}
\put(218.27,15.51){\line(1,0){1.132}}
\put(219.4,15.39){\line(1,0){1.136}}
\put(220.54,15.33){\line(1,0){1.138}}
\put(221.68,15.31){\line(1,0){1.137}}
\put(222.81,15.34){\line(1,0){1.135}}
\multiput(223.95,15.42)(.2828,.0309){4}{\line(1,0){.2828}}
\multiput(225.08,15.54)(.18749,.02847){6}{\line(1,0){.18749}}
\multiput(226.2,15.71)(.15955,.03109){7}{\line(1,0){.15955}}
\multiput(227.32,15.93)(.13835,.033){8}{\line(1,0){.13835}}
\multiput(228.43,16.19)(.10948,.031){10}{\line(1,0){.10948}}
\multiput(229.52,16.5)(.09826,.03231){11}{\line(1,0){.09826}}
\multiput(230.6,16.86)(.08876,.03335){12}{\line(1,0){.08876}}
\multiput(231.67,17.26)(.074822,.031736){14}{\line(1,0){.074822}}
\multiput(232.72,17.7)(.068537,.032509){15}{\line(1,0){.068537}}
\multiput(233.74,18.19)(.062925,.033132){16}{\line(1,0){.062925}}
\multiput(234.75,18.72)(.057871,.033627){17}{\line(1,0){.057871}}
\multiput(235.73,19.29)(.050479,.032222){19}{\line(1,0){.050479}}
\multiput(236.69,19.9)(.046636,.032585){20}{\line(1,0){.046636}}
\multiput(237.63,20.56)(.043081,.03286){21}{\line(1,0){.043081}}
\multiput(238.53,21.25)(.039778,.033055){22}{\line(1,0){.039778}}
\multiput(239.41,21.97)(.036696,.033178){23}{\line(1,0){.036696}}
\multiput(240.25,22.74)(.033809,.033236){24}{\line(1,0){.033809}}
\multiput(241.06,23.53)(.032393,.034618){24}{\line(0,1){.034618}}
\multiput(241.84,24.36)(.033731,.039207){22}{\line(0,1){.039207}}
\multiput(242.58,25.23)(.033593,.042512){21}{\line(0,1){.042512}}
\multiput(243.29,26.12)(.033379,.046071){20}{\line(0,1){.046071}}
\multiput(243.95,27.04)(.033081,.04992){19}{\line(0,1){.04992}}
\multiput(244.58,27.99)(.03269,.054105){18}{\line(0,1){.054105}}
\multiput(245.17,28.96)(.032192,.058682){17}{\line(0,1){.058682}}
\multiput(245.72,29.96)(.033677,.06797){15}{\line(0,1){.06797}}
\multiput(246.22,30.98)(.033012,.074268){14}{\line(0,1){.074268}}
\multiput(246.69,32.02)(.032182,.081395){13}{\line(0,1){.081395}}
\multiput(247.1,33.08)(.03115,.08956){12}{\line(0,1){.08956}}
\multiput(247.48,34.15)(.03287,.10893){10}{\line(0,1){.10893}}
\multiput(247.81,35.24)(.03144,.12246){9}{\line(0,1){.12246}}
\multiput(248.09,36.34)(.02959,.13912){8}{\line(0,1){.13912}}
\multiput(248.33,37.46)(.03167,.18698){6}{\line(0,1){.18698}}
\multiput(248.52,38.58)(.02861,.22576){5}{\line(0,1){.22576}}
\put(248.66,39.71){\line(0,1){1.134}}
\put(248.76,40.84){\line(0,1){1.728}}
\put(50.73,27.88){\makebox(0,0)[cc]{$\sigma\theta$}}
\put(218.5,28.06){\makebox(0,0)[cc]{$\sigma\theta$}}
\put(127.99,35.96){\makebox(0,0)[cc]{$\breve{\q}_{\sigma
\vartheta}$}}
\put(70.64,56.8){\line(1,0){.999}}
\put(72.64,56.79){\line(1,0){.999}}
\put(74.64,56.78){\line(1,0){.999}}
\put(76.64,56.77){\line(1,0){.999}}
\put(78.64,56.76){\line(1,0){.999}}
\put(80.63,56.74){\line(1,0){.999}}
\put(82.63,56.73){\line(1,0){.999}}
\put(84.63,56.72){\line(1,0){.999}}
\put(86.63,56.71){\line(1,0){.999}}
\put(88.63,56.7){\line(1,0){.999}}
\put(90.63,56.69){\line(1,0){.999}}
\put(92.63,56.68){\line(1,0){.999}}
\put(94.63,56.67){\line(1,0){.999}}
\put(96.63,56.66){\line(1,0){.999}}
\put(98.62,56.65){\line(1,0){.999}}
\put(100.62,56.63){\line(1,0){.999}}
\put(102.62,56.62){\line(1,0){.999}}
\put(104.62,56.61){\line(1,0){.999}}
\put(106.62,56.6){\line(1,0){.999}}
\put(108.62,56.59){\line(1,0){.999}}
\put(110.62,56.58){\line(1,0){.999}}
\put(112.62,56.57){\line(1,0){.999}}
\put(114.62,56.56){\line(1,0){.999}}
\put(116.61,56.55){\line(1,0){.999}}
\put(118.61,56.54){\line(1,0){.999}}
\put(120.61,56.52){\line(1,0){.999}}
\put(122.61,56.51){\line(1,0){.999}}
\put(124.61,56.5){\line(1,0){.999}}
\put(126.61,56.49){\line(1,0){.999}}
\put(128.61,56.48){\line(1,0){.999}}
\put(130.61,56.47){\line(1,0){.999}}
\put(132.61,56.46){\line(1,0){.999}}
\put(134.6,56.45){\line(1,0){.999}}
\put(136.6,56.44){\line(1,0){.999}}
\put(138.6,56.42){\line(1,0){.999}}
\put(140.6,56.41){\line(1,0){.999}}
\put(142.6,56.4){\line(1,0){.999}}
\put(144.6,56.39){\line(1,0){.999}}
\put(146.6,56.38){\line(1,0){.999}}
\put(148.6,56.37){\line(1,0){.999}}
\put(150.6,56.36){\line(1,0){.999}}
\put(152.59,56.35){\line(1,0){.999}}
\put(154.59,56.34){\line(1,0){.999}}
\put(156.59,56.33){\line(1,0){.999}}
\put(158.59,56.31){\line(1,0){.999}}
\put(160.59,56.3){\line(1,0){.999}}
\put(162.59,56.29){\line(1,0){.999}}
\put(164.59,56.28){\line(1,0){.999}}
\put(166.59,56.27){\line(1,0){.999}}
\put(168.59,56.26){\line(1,0){.999}}
\put(170.58,56.25){\line(1,0){.999}}
\put(172.58,56.24){\line(1,0){.999}}
\put(174.58,56.23){\line(1,0){.999}}
\put(176.58,56.22){\line(1,0){.999}}
\put(178.58,56.2){\line(1,0){.999}}
\put(180.58,56.19){\line(1,0){.999}}
\put(182.58,56.18){\line(1,0){.999}}
\put(184.58,56.17){\line(1,0){.999}}
\put(186.58,56.16){\line(1,0){.999}}
\put(188.57,56.15){\line(1,0){.999}}
\put(190.57,56.14){\line(1,0){.999}}
\put(192.57,56.13){\line(1,0){.999}}
\put(194.57,56.12){\line(1,0){.999}}
\put(196.57,56.11){\line(1,0){.999}}
\put(70.46,28.52){\line(1,0){.995}}
\put(72.45,28.55){\line(1,0){.995}}
\put(74.44,28.58){\line(1,0){.995}}
\put(76.43,28.6){\line(1,0){.995}}
\put(78.42,28.63){\line(1,0){.995}}
\put(80.41,28.66){\line(1,0){.995}}
\put(82.4,28.69){\line(1,0){.995}}
\put(84.39,28.71){\line(1,0){.995}}
\put(86.38,28.74){\line(1,0){.995}}
\put(88.37,28.77){\line(1,0){.995}}
\put(90.37,28.8){\line(1,0){.995}}
\put(92.36,28.83){\line(1,0){.995}}
\put(94.35,28.85){\line(1,0){.995}}
\put(96.34,28.88){\line(1,0){.995}}
\put(98.33,28.91){\line(1,0){.995}}
\put(100.32,28.94){\line(1,0){.995}}
\put(102.31,28.97){\line(1,0){.995}}
\put(104.3,28.99){\line(1,0){.995}}
\put(106.29,29.02){\line(1,0){.995}}
\put(108.28,29.05){\line(1,0){.995}}
\put(110.27,29.08){\line(1,0){.995}}
\put(112.26,29.11){\line(1,0){.995}}
\put(114.25,29.13){\line(1,0){.995}}
\put(116.24,29.16){\line(1,0){.995}}
\put(118.23,29.19){\line(1,0){.995}}
\put(120.22,29.22){\line(1,0){.995}}
\put(122.21,29.24){\line(1,0){.995}}
\put(124.2,29.27){\line(1,0){.995}}
\put(126.2,29.3){\line(1,0){.995}}
\put(128.19,29.33){\line(1,0){.995}}
\put(130.18,29.36){\line(1,0){.995}}
\put(132.17,29.38){\line(1,0){.995}}
\put(134.16,29.41){\line(1,0){.995}}
\put(136.15,29.44){\line(1,0){.995}}
\put(138.14,29.47){\line(1,0){.995}}
\put(140.13,29.5){\line(1,0){.995}}
\put(142.12,29.52){\line(1,0){.995}}
\put(144.11,29.55){\line(1,0){.995}}
\put(146.1,29.58){\line(1,0){.995}}
\put(148.09,29.61){\line(1,0){.995}}
\put(150.08,29.63){\line(1,0){.995}}
\put(152.07,29.66){\line(1,0){.995}}
\put(154.06,29.69){\line(1,0){.995}}
\put(156.05,29.72){\line(1,0){.995}}
\put(158.04,29.75){\line(1,0){.995}}
\put(160.03,29.77){\line(1,0){.995}}
\put(162.03,29.8){\line(1,0){.995}}
\put(164.02,29.83){\line(1,0){.995}}
\put(166.01,29.86){\line(1,0){.995}}
\put(168,29.89){\line(1,0){.995}}
\put(169.99,29.91){\line(1,0){.995}}
\put(171.98,29.94){\line(1,0){.995}}
\put(173.97,29.97){\line(1,0){.995}}
\put(175.96,30){\line(1,0){.995}}
\put(177.95,30.02){\line(1,0){.995}}
\put(179.94,30.05){\line(1,0){.995}}
\put(181.93,30.08){\line(1,0){.995}}
\put(183.92,30.11){\line(1,0){.995}}
\put(185.91,30.14){\line(1,0){.995}}
\put(187.9,30.16){\line(1,0){.995}}
\put(189.89,30.19){\line(1,0){.995}}
\put(191.88,30.22){\line(1,0){.995}}
\put(193.87,30.25){\line(1,0){.995}}
\put(195.86,30.28){\line(1,0){.995}}
\put(150.26,42.91){\line(0,1){13.43}}
\put(156.45,49.8){\makebox(0,0)[cc]{$\theta$}} \thicklines
\qbezier(47.2,43.08)(68.94,64.56)(74.42,113.27)
\qbezier(74.42,113.27)(88.74,112.74)(94.58,123.52)
\qbezier(94.58,123.52)(129.58,-1.55)(171.65,123.34)
\qbezier(171.65,123.34)(179.78,109.46)(194.98,113.62)
\qbezier(194.98,113.62)(191.8,57.23)(221.85,42.55)
\put(127.1,77.2){\makebox(0,0)[cc]{$\pp\setminus\q$}}
\put(41.32,43.7){\makebox(0,0)[cc]{$x$}}
\put(227.73,42.81){\makebox(0,0)[cc]{$y$}} \thinlines
\multiput(47.2,42.74)(.033727,-.044345){17}{\line(0,-1){.044345}}
\multiput(48.35,41.23)(.033727,-.044345){17}{\line(0,-1){.044345}}
\multiput(49.49,39.73)(.033727,-.044345){17}{\line(0,-1){.044345}}
\multiput(50.64,38.22)(.033727,-.044345){17}{\line(0,-1){.044345}}
\multiput(51.79,36.71)(.033727,-.044345){17}{\line(0,-1){.044345}}
\multiput(52.93,35.2)(.033727,-.044345){17}{\line(0,-1){.044345}}
\multiput(54.08,33.69)(.033727,-.044345){17}{\line(0,-1){.044345}}
\multiput(55.23,32.19)(.033727,-.044345){17}{\line(0,-1){.044345}}
\multiput(56.37,30.68)(.033727,-.044345){17}{\line(0,-1){.044345}}
\multiput(57.52,29.17)(.033727,-.044345){17}{\line(0,-1){.044345}}
\multiput(58.67,27.66)(.033727,-.044345){17}{\line(0,-1){.044345}}
\multiput(59.81,26.16)(.033727,-.044345){17}{\line(0,-1){.044345}}
\multiput(60.96,24.65)(.033727,-.044345){17}{\line(0,-1){.044345}}
\multiput(62.11,23.14)(.033727,-.044345){17}{\line(0,-1){.044345}}
\multiput(221.12,42.74)(-.033033,-.041882){18}{\line(0,-1){.041882}}
\multiput(219.93,41.23)(-.033033,-.041882){18}{\line(0,-1){.041882}}
\multiput(218.74,39.73)(-.033033,-.041882){18}{\line(0,-1){.041882}}
\multiput(217.55,38.22)(-.033033,-.041882){18}{\line(0,-1){.041882}}
\multiput(216.36,36.71)(-.033033,-.041882){18}{\line(0,-1){.041882}}
\multiput(215.18,35.2)(-.033033,-.041882){18}{\line(0,-1){.041882}}
\multiput(213.99,33.69)(-.033033,-.041882){18}{\line(0,-1){.041882}}
\multiput(212.8,32.19)(-.033033,-.041882){18}{\line(0,-1){.041882}}
\multiput(211.61,30.68)(-.033033,-.041882){18}{\line(0,-1){.041882}}
\multiput(210.42,29.17)(-.033033,-.041882){18}{\line(0,-1){.041882}}
\multiput(209.23,27.66)(-.033033,-.041882){18}{\line(0,-1){.041882}}
\multiput(208.04,26.16)(-.033033,-.041882){18}{\line(0,-1){.041882}}
\multiput(206.85,24.65)(-.033033,-.041882){18}{\line(0,-1){.041882}}
\multiput(205.66,23.14)(-.033033,-.041882){18}{\line(0,-1){.041882}}
\end{picture}
\unitlength .32mm 
\linethickness{0.4pt}
\ifx\plotpoint\undefined\newsavebox{\plotpoint}\fi 
\begin{picture}(100.75,125.5)(20,0)
\multiput(119.75,33.75)(-.0473163842,.0337217514){1416}{\line(-1,0){.0473163842}}
\multiput(119.5,34)(.0342055485,.0337326608){1586}{\line(1,0){.0342055485}}
\thicklines \qbezier(53.25,81)(73.88,91.75)(75,124.5)
\qbezier(75,124.5)(100.88,122)(106.25,142.5)
\qbezier(106.25,142.5)(122.25,122.63)(147.25,132.25)
\qbezier(147.25,132.25)(138.13,97.88)(173.5,87)
\put(160.25,114.5){\makebox(0,0)[cc]{$\oo_x$}}
\put(121,29){\makebox(0,0)[cc]{$x$}} \thinlines
\put(148.82,34){\line(0,1){1.191}}
\put(148.79,35.19){\line(0,1){1.189}}
\multiput(148.72,36.38)(-.0302,.2963){4}{\line(0,1){.2963}}
\multiput(148.6,37.56)(-.02814,.19653){6}{\line(0,1){.19653}}
\multiput(148.43,38.74)(-.03094,.16733){7}{\line(0,1){.16733}}
\multiput(148.22,39.92)(-.033,.14519){8}{\line(0,1){.14519}}
\multiput(147.95,41.08)(-.0311,.11499){10}{\line(0,1){.11499}}
\multiput(147.64,42.23)(-.03249,.1033){11}{\line(0,1){.1033}}
\multiput(147.28,43.36)(-.03361,.0934){12}{\line(0,1){.0934}}
\multiput(146.88,44.48)(-.032033,.078824){14}{\line(0,1){.078824}}
\multiput(146.43,45.59)(-.032861,.072294){15}{\line(0,1){.072294}}
\multiput(145.94,46.67)(-.033535,.066468){16}{\line(0,1){.066468}}
\multiput(145.4,47.74)(-.032184,.057823){18}{\line(0,1){.057823}}
\multiput(144.82,48.78)(-.032691,.053496){19}{\line(0,1){.053496}}
\multiput(144.2,49.79)(-.033095,.049518){20}{\line(0,1){.049518}}
\multiput(143.54,50.78)(-.033409,.04584){21}{\line(0,1){.04584}}
\multiput(142.84,51.75)(-.033641,.042425){22}{\line(0,1){.042425}}
\multiput(142.1,52.68)(-.032392,.037605){24}{\line(0,1){.037605}}
\multiput(141.32,53.58)(-.032537,.034808){25}{\line(0,1){.034808}}
\multiput(140.51,54.45)(-.033924,.033457){25}{\line(-1,0){.033924}}
\multiput(139.66,55.29)(-.036724,.033387){24}{\line(-1,0){.036724}}
\multiput(138.78,56.09)(-.039704,.033254){23}{\line(-1,0){.039704}}
\multiput(137.86,56.85)(-.042887,.03305){22}{\line(-1,0){.042887}}
\multiput(136.92,57.58)(-.046299,.03277){21}{\line(-1,0){.046299}}
\multiput(135.95,58.27)(-.049971,.032406){20}{\line(-1,0){.049971}}
\multiput(134.95,58.92)(-.05694,.033721){18}{\line(-1,0){.05694}}
\multiput(133.92,59.52)(-.06169,.033227){17}{\line(-1,0){.06169}}
\multiput(132.88,60.09)(-.066926,.032612){16}{\line(-1,0){.066926}}
\multiput(131.8,60.61)(-.072742,.031857){15}{\line(-1,0){.072742}}
\multiput(130.71,61.09)(-.085357,.033318){13}{\line(-1,0){.085357}}
\multiput(129.6,61.52)(-.09386,.03231){12}{\line(-1,0){.09386}}
\multiput(128.48,61.91)(-.10374,.03106){11}{\line(-1,0){.10374}}
\multiput(127.34,62.25)(-.12823,.03278){9}{\line(-1,0){.12823}}
\multiput(126.18,62.55)(-.14564,.03098){8}{\line(-1,0){.14564}}
\multiput(125.02,62.79)(-.1957,.03339){6}{\line(-1,0){.1957}}
\multiput(123.84,62.99)(-.23628,.0305){5}{\line(-1,0){.23628}}
\put(122.66,63.15){\line(-1,0){1.187}}
\put(121.47,63.25){\line(-1,0){1.19}}
\put(120.29,63.31){\line(-1,0){1.191}}
\put(119.09,63.32){\line(-1,0){1.19}}
\put(117.9,63.27){\line(-1,0){1.188}}
\multiput(116.72,63.19)(-.23665,-.02744){5}{\line(-1,0){.23665}}
\multiput(115.53,63.05)(-.19612,-.03086){6}{\line(-1,0){.19612}}
\multiput(114.36,62.86)(-.16689,-.03326){7}{\line(-1,0){.16689}}
\multiput(113.19,62.63)(-.12864,-.03112){9}{\line(-1,0){.12864}}
\multiput(112.03,62.35)(-.11455,-.03268){10}{\line(-1,0){.11455}}
\multiput(110.88,62.02)(-.09427,-.03109){12}{\line(-1,0){.09427}}
\multiput(109.75,61.65)(-.08578,-.032212){13}{\line(-1,0){.08578}}
\multiput(108.64,61.23)(-.078373,-.033122){14}{\line(-1,0){.078373}}
\multiput(107.54,60.77)(-.067342,-.031743){16}{\line(-1,0){.067342}}
\multiput(106.46,60.26)(-.062115,-.032426){17}{\line(-1,0){.062115}}
\multiput(105.41,59.71)(-.057372,-.032982){18}{\line(-1,0){.057372}}
\multiput(104.37,59.12)(-.053038,-.033428){19}{\line(-1,0){.053038}}
\multiput(103.37,58.48)(-.046719,-.032169){21}{\line(-1,0){.046719}}
\multiput(102.39,57.8)(-.043311,-.032493){22}{\line(-1,0){.043311}}
\multiput(101.43,57.09)(-.040131,-.032737){23}{\line(-1,0){.040131}}
\multiput(100.51,56.34)(-.037153,-.03291){24}{\line(-1,0){.037153}}
\multiput(99.62,55.55)(-.034354,-.033016){25}{\line(-1,0){.034354}}
\multiput(98.76,54.72)(-.032984,-.034384){25}{\line(0,-1){.034384}}
\multiput(97.93,53.86)(-.032876,-.037183){24}{\line(0,-1){.037183}}
\multiput(97.15,52.97)(-.032701,-.040161){23}{\line(0,-1){.040161}}
\multiput(96.39,52.05)(-.032453,-.043341){22}{\line(0,-1){.043341}}
\multiput(95.68,51.09)(-.033732,-.049085){20}{\line(0,-1){.049085}}
\multiput(95,50.11)(-.03338,-.053069){19}{\line(0,-1){.053069}}
\multiput(94.37,49.1)(-.03293,-.057402){18}{\line(0,-1){.057402}}
\multiput(93.78,48.07)(-.032369,-.062145){17}{\line(0,-1){.062145}}
\multiput(93.23,47.01)(-.031682,-.067371){16}{\line(0,-1){.067371}}
\multiput(92.72,45.93)(-.03305,-.078403){14}{\line(0,-1){.078403}}
\multiput(92.26,44.84)(-.032133,-.08581){13}{\line(0,-1){.08581}}
\multiput(91.84,43.72)(-.03101,-.0943){12}{\line(0,-1){.0943}}
\multiput(91.47,42.59)(-.03258,-.11458){10}{\line(0,-1){.11458}}
\multiput(91.14,41.44)(-.031,-.12867){9}{\line(0,-1){.12867}}
\multiput(90.86,40.29)(-.0331,-.16692){7}{\line(0,-1){.16692}}
\multiput(90.63,39.12)(-.03068,-.19615){6}{\line(0,-1){.19615}}
\multiput(90.45,37.94)(-.02722,-.23668){5}{\line(0,-1){.23668}}
\put(90.31,36.76){\line(0,-1){1.188}}
\put(90.22,35.57){\line(0,-1){2.382}}
\put(90.19,33.19){\line(0,-1){1.19}}
\put(90.25,32){\line(0,-1){1.187}}
\multiput(90.36,30.81)(.03071,-.23625){5}{\line(0,-1){.23625}}
\multiput(90.51,29.63)(.03357,-.19567){6}{\line(0,-1){.19567}}
\multiput(90.71,28.46)(.03112,-.14561){8}{\line(0,-1){.14561}}
\multiput(90.96,27.29)(.0329,-.1282){9}{\line(0,-1){.1282}}
\multiput(91.26,26.14)(.03115,-.10371){11}{\line(0,-1){.10371}}
\multiput(91.6,25)(.03239,-.09383){12}{\line(0,-1){.09383}}
\multiput(91.99,23.87)(.033397,-.085326){13}{\line(0,-1){.085326}}
\multiput(92.42,22.76)(.031924,-.072713){15}{\line(0,-1){.072713}}
\multiput(92.9,21.67)(.032673,-.066896){16}{\line(0,-1){.066896}}
\multiput(93.42,20.6)(.033283,-.06166){17}{\line(0,-1){.06166}}
\multiput(93.99,19.55)(.031996,-.053914){19}{\line(0,-1){.053914}}
\multiput(94.6,18.53)(.032452,-.049941){20}{\line(0,-1){.049941}}
\multiput(95.25,17.53)(.032813,-.046269){21}{\line(0,-1){.046269}}
\multiput(95.93,16.56)(.03309,-.042857){22}{\line(0,-1){.042857}}
\multiput(96.66,15.61)(.03329,-.039674){23}{\line(0,-1){.039674}}
\multiput(97.43,14.7)(.033421,-.036693){24}{\line(0,-1){.036693}}
\multiput(98.23,13.82)(.033489,-.033893){25}{\line(0,-1){.033893}}
\multiput(99.07,12.97)(.034838,-.032505){25}{\line(1,0){.034838}}
\multiput(99.94,12.16)(.037635,-.032357){24}{\line(1,0){.037635}}
\multiput(100.84,11.39)(.042456,-.033602){22}{\line(1,0){.042456}}
\multiput(101.78,10.65)(.045871,-.033366){21}{\line(1,0){.045871}}
\multiput(102.74,9.95)(.049548,-.033049){20}{\line(1,0){.049548}}
\multiput(103.73,9.28)(.053526,-.032642){19}{\line(1,0){.053526}}
\multiput(104.75,8.66)(.057852,-.032131){18}{\line(1,0){.057852}}
\multiput(105.79,8.09)(.066499,-.033474){16}{\line(1,0){.066499}}
\multiput(106.85,7.55)(.072324,-.032795){15}{\line(1,0){.072324}}
\multiput(107.94,7.06)(.078853,-.031961){14}{\line(1,0){.078853}}
\multiput(109.04,6.61)(.09343,-.03352){12}{\line(1,0){.09343}}
\multiput(110.16,6.21)(.10333,-.0324){11}{\line(1,0){.10333}}
\multiput(111.3,5.85)(.11502,-.03099){10}{\line(1,0){.11502}}
\multiput(112.45,5.54)(.14522,-.03287){8}{\line(1,0){.14522}}
\multiput(113.61,5.28)(.16736,-.03079){7}{\line(1,0){.16736}}
\multiput(114.78,5.06)(.23586,-.03355){5}{\line(1,0){.23586}}
\put(115.96,4.9){\line(1,0){1.185}}
\put(117.15,4.78){\line(1,0){1.189}}
\put(118.34,4.7){\line(1,0){2.382}}
\put(120.72,4.71){\line(1,0){1.189}}
\multiput(121.91,4.78)(.2962,.0305){4}{\line(1,0){.2962}}
\multiput(123.09,4.9)(.1965,.02832){6}{\line(1,0){.1965}}
\multiput(124.27,5.07)(.1673,.03109){7}{\line(1,0){.1673}}
\multiput(125.44,5.29)(.14516,.03313){8}{\line(1,0){.14516}}
\multiput(126.6,5.56)(.11496,.0312){10}{\line(1,0){.11496}}
\multiput(127.75,5.87)(.10327,.03259){11}{\line(1,0){.10327}}
\multiput(128.89,6.23)(.09337,.03369){12}{\line(1,0){.09337}}
\multiput(130.01,6.63)(.078794,.032105){14}{\line(1,0){.078794}}
\multiput(131.11,7.08)(.072263,.032928){15}{\line(1,0){.072263}}
\multiput(132.2,7.57)(.066437,.033596){16}{\line(1,0){.066437}}
\multiput(133.26,8.11)(.057793,.032237){18}{\line(1,0){.057793}}
\multiput(134.3,8.69)(.053466,.03274){19}{\line(1,0){.053466}}
\multiput(135.32,9.31)(.049487,.03314){20}{\line(1,0){.049487}}
\multiput(136.31,9.98)(.04581,.033451){21}{\line(1,0){.04581}}
\multiput(137.27,10.68)(.042394,.03368){22}{\line(1,0){.042394}}
\multiput(138.2,11.42)(.037575,.032426){24}{\line(1,0){.037575}}
\multiput(139.1,12.2)(.034778,.032569){25}{\line(1,0){.034778}}
\multiput(139.97,13.01)(.033426,.033955){25}{\line(0,1){.033955}}
\multiput(140.81,13.86)(.033354,.036755){24}{\line(0,1){.036755}}
\multiput(141.61,14.74)(.033217,.039735){23}{\line(0,1){.039735}}
\multiput(142.37,15.66)(.033011,.042917){22}{\line(0,1){.042917}}
\multiput(143.1,16.6)(.032728,.046329){21}{\line(0,1){.046329}}
\multiput(143.78,17.57)(.03236,.050001){20}{\line(0,1){.050001}}
\multiput(144.43,18.57)(.033669,.056971){18}{\line(0,1){.056971}}
\multiput(145.04,19.6)(.03317,.061721){17}{\line(0,1){.061721}}
\multiput(145.6,20.65)(.03255,.066956){16}{\line(0,1){.066956}}
\multiput(146.12,21.72)(.03179,.072771){15}{\line(0,1){.072771}}
\multiput(146.6,22.81)(.03324,.085387){13}{\line(0,1){.085387}}
\multiput(147.03,23.92)(.03222,.09389){12}{\line(0,1){.09389}}
\multiput(147.42,25.05)(.03096,.10377){11}{\line(0,1){.10377}}
\multiput(147.76,26.19)(.03266,.12826){9}{\line(0,1){.12826}}
\multiput(148.05,27.34)(.03085,.14567){8}{\line(0,1){.14567}}
\multiput(148.3,28.51)(.03321,.19573){6}{\line(0,1){.19573}}
\multiput(148.5,29.68)(.03028,.2363){5}{\line(0,1){.2363}}
\put(148.65,30.87){\line(0,1){1.187}}
\put(148.75,32.05){\line(0,1){1.948}}
\put(119.18,34.18){\line(1,0){.991}}
\put(121.16,34.2){\line(1,0){.991}}
\put(123.15,34.21){\line(1,0){.991}}
\put(125.13,34.23){\line(1,0){.991}}
\put(127.11,34.25){\line(1,0){.991}}
\put(129.09,34.27){\line(1,0){.991}}
\put(131.08,34.28){\line(1,0){.991}}
\put(133.06,34.3){\line(1,0){.991}}
\put(135.04,34.32){\line(1,0){.991}}
\put(137.02,34.33){\line(1,0){.991}}
\put(139.01,34.35){\line(1,0){.991}}
\put(140.99,34.37){\line(1,0){.991}}
\put(142.97,34.39){\line(1,0){.991}}
\put(144.96,34.4){\line(1,0){.991}}
\put(146.94,34.42){\line(1,0){.991}}
\put(135.75,28){\makebox(0,0)[cc]{$\nu\vartheta$}}
\end{picture}

\begin{figure}[!ht]
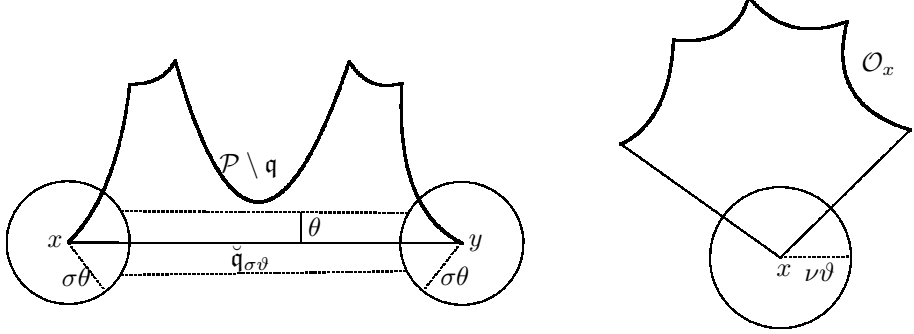

\centering \caption{Properties ($F_1$) and ($F_2$).} \label{fig1}
\end{figure}

\begin{defn}[fat polygons]\label{deffat}
Let $\vartheta>0$, $\sigma \geq 1$ and $\nu \geq 4\sigma$. We call
a $k$-gon $P$ with quasi-geodesic edges \textit{$(\vartheta ,
\sigma, \nu )$--fat} if the following properties hold:
\begin{enumerate}
  \item[$(F_1)$] for every edge $\q $ we have, with the notation
\ref{breve}, that
$$\dist \left( \breve{\q}_{\sigma \vartheta}\, ,\, P \setminus
\q \right) \geq \vartheta;$$

  \item[$(F_2)$] for
every vertex $x$ we have $$\dist(x, \oo_x)\geq \nu \vartheta.$$
\end{enumerate}

When $\sigma =2$ we say that $P$ is \textit{$(\vartheta , \nu
)$--fat}.
\end{defn}

\begin{thm}[\cite{DrutuSapir:TreeGraded},
\cite{Drutu:RelHyp}]\label{tgi}
Let $(X,\dist)$ be a geodesic metric space
and let $\aaa$ be a collection of subsets of $X$. The metric space
$X$ is asymptotically tree-graded with respect to $\aaa$ if and
only if the following properties are satisfied:

\medskip

\begin{itemize}
  \item[$(\alpha_1)$] For every $\delta >0$ the diameters of the
intersections $\nbhd_{\delta}(A)\cap \nbhd_{\delta}(A')$ are
uniformly bounded for distinct pairs of $A,A'\in\aaa$.

\medskip

\item[$(\alpha_2)$] There exists $\varepsilon$ in $\left( 0, \frac{1}{2}
\right)$ and $M>0$ such that for every geodesic $\fg$ of length
$\ell$ and every $A\in \aaa$ with $\fg(0),\fg(\ell)\in
\nn_{\varepsilon \ell}(A)$ we have $\fg([0, \ell ])\cap
\nn_{M}(A)\neq \emptyset$.

\medskip

\item[$(\beta_3)$] There exists $\vartheta >0$, $\nu \geq 8$ and $\chi>0$
        such that any  $(\vartheta , \nu)$--fat geodesic hexagon is
    contained in $\nn_\chi (A)$, for some $A\in\aaa$.
        \end{itemize}
\end{thm}

\begin{rmk}\label{rtgi}
In Theorem \ref{tgi}, property $(\alpha_2)$ can be replaced by the
following stronger property:
\begin{itemize}
    \item[$(\beta_2)$] There exists $\epsilon >0$ and $M\geq 0$
        such that for any geodesic $\fg$ of length $\ell$ and any
        $A\in \aaa$ satisfying $\fg(0),\fg(\ell) \in \nn_{\epsilon\ell}
        (A)$, the middle third
        $\fg \left( \left[ \frac{\ell}{3}\, ,\, \frac{2\ell}{3} \right] \right)$
         is contained in $\nn_M (A)$.
\end{itemize}
\end{rmk}

\medskip

The notion of asymptotically tree-graded space relates to the
standard definition of (strong) relative hyperbolicity by the
following.

\begin{thm}[Dru\c{t}u--Osin--Sapir \cite{DrutuSapir:TreeGraded}]
    \label{DrutuSapir:TreeGraded:relhypclassification}
A finitely generated group $G$ is hyperbolic relative to a
finite collection of finitely generated subgroups $\H$ if and only
if $G$ is asymptotically tree-graded with respect to $\lll\H$.
\end{thm}

The converse statement of the above theorem can be strengthened as
follows.

\begin{thm}[Dru\c{t}u \cite{Drutu:RelHyp}]\label{t3}
If $G$ is a finitely generated group which is asymptotically tree-graded with
respect to a collection $\aaa$ of subsets, then $G$ is either
hyperbolic or it is relatively hyperbolic with respect to a finite
family of finitely generated subgroups $\{H_1,...,H_m\}$ such that
every $H_i$ is contained in $\nn_\varkappa (A_i)$ for some $A_i\in
\aaa$, where $\varkappa$ is the maximum between the constant $M$
in $(\beta_2)$ and the constant $\chi$ in $(\beta_3)$.
\end{thm}

A consequence of this is the following result:

\begin{thm}[relative hyperbolicity is rigid,
    Dru\c{t}u \cite{Drutu:RelHyp}]\label{t4}
If a group $G'$ is quasi-isometric to a relatively hyperbolic
group $G$ then $G'$ is also relatively hyperbolic.
\end{thm}

Note that formulating a relation between the peripheral
subgroups of $G$ and of $G'$ is, in general, nontrivial. This can be seen for
instance when $G=G'=A*B*C$, since $G$ is hyperbolic relative to $\{
A,B,C \}$, and also hyperbolic relative to $\{ A*B,C \}$.

\section{Unconstricted and constricted metric spaces}\label{sunc}

\begin{defn}\label{ums}
A metric space $B$ is \emph{unconstricted} if the
following two properties hold:
\begin{itemize}
\item[(1)] there exists an ultrafilter $\omega$ and a
sequence $d$ such that for every sequence of observation points
$b$, $\con (B,b,d)$ does not have cut-points;
\item[(2)] for some constant $c$, every point in $B$
is at distance at most $c$ from a bi-infinite geodesic in $B$.
\end{itemize}
\end{defn}

When $B$ is an infinite finitely generated group, being
unconstricted means simply that \emph{at least one} of its asymptotic cones
does not have cut-points. Opposite to it, a \textit{constricted}
group is a group with cut-points in every asymptotic cone. See the
list following Definition~\ref{duunc} for examples of
unconstricted groups.

\begin{rmk}
Theorem \ref{DrutuSapir:TreeGraded:relhypclassification} implies
that relatively hyperbolic groups are constricted. Thus,
unconstricted groups are particular cases of NRH groups. They play
an essential part in the notion we introduce, of thick group.
\end{rmk}

Note that the definition above slightly differs from the one in
\cite{DrutuSapir:TreeGraded} in that property (2) has been added.
We incorporate this condition into the definition as it is a
required hypothesis for all the quasi-isometry rigidity results we
obtain. Since, up to bi-Lipschitz homeomorphism, the set of
asymptotic cones is a quasi-isometry invariant of a metric space
$B$, it follows that constrictedness and unconstrictedness are
quasi-isometry invariants.

The property of being constricted is related to the divergence of
geodesics \cite{Gersten:divergence}. Let $X$ be a geodesic metric
space. Given a geodesic segment $\fc\colon[-R,R] \to X$, its
\textit{divergence} is a function $div_\fc\colon (0, R] \to \R_+$,
where for every $r>0$ we define $div_\fg (r)$ as the distance
between $\fc (-r)$ and $\fc (r)$ in $X\setminus B\left( \fc \left(
0 \right) , r \right)$ endowed with the length metric (with the
assumption that $\fc (-r)$ and $\fc (r)$ can be joined in
$X\setminus B\left( \fc \left( 0 \right) , r \right)$ by a path of
finite length). To a complete minimizing geodesic $\fg :\R \to X$
is associated a function $div_\fg$ defined similarly on $\R_+$. By
a slight abuse of terminology, it is standard to refer to the
growth rate of the function $div_\fg$ as \textit{the divergence
of}~$\fg$.

A geodesic in a metric space $X$ is called \textit{periodic} if
its stabilizer in the group of isometries of $X$ is co-bounded. By
combining Proposition~4.2 of \cite{KKL:QI} with
Lemma~\ref{cutting}, we obtain:
\begin{lem}\label{suplin}
Let $\fg :\R \to X $ be a periodic geodesic. If $\fg$ has
superlinear divergence, then in any asymptotic cone, $\cone(X)$,
for which the limit of $\fg$ is nonempty there exists a collection
of proper subsets of $\cone(X)$ with respect to which it is
tree-graded. Furthermore, in this case one has that the limit of
$\fg$ is a transversal geodesic. \qed\end{lem}

\begin{defn}\label{duunc}
A collection of metric spaces, $\B$, is \emph{uniformly
unconstricted} if:
\begin{itemize}
\item[(1)] for some constant $c$, every point in every
space $B\in \B$ is at distance at most $c$ from a bi-infinite
geodesic in $B$;
\item[(2)] for every sequence of spaces $(B_{i}, \dist_i)$ in
$\B$, there exists an ultrafilter $\omega$ and a
sequence of scaling constants $d$ so that for every sequence
of basepoints $b=(b_i)$ with $b_i\in B_i$,
$\ulim (B_{i}, b_i, 1/d_i\, \dist_i)$ does not have cut-points.
\end{itemize}
\end{defn}

\medskip

Recall that a group is \emph{elementary} if it is virtually
cyclic.

\medskip

\textit{Examples of uniformly unconstricted collections of spaces:}
\begin{enumerate}
\item The collection of all cartesian products of
geodesic metric spaces of infinite diameter. This follows from the
fact that every ultralimit of a sequence of such spaces appears as
cartesian product of two non-trivial geodesic metric spaces. Such
a cartesian product cannot have a global cut-point, because
Euclidean rectangles do not have cut-points.
\item The collection of finitely generated non-elementary
groups
with a central element of infinite order is uniformly
unconstricted \cite[Theorem~6.7]{DrutuSapir:TreeGraded}.
\item The collection of finitely generated non-elementary
groups satisfying the same identity is uniformly unconstricted
\cite[Theorem~6.12]{DrutuSapir:TreeGraded}. Recall that a group
$G$ is said to \emph{satisfy an identity (a law)}
    if there exists a word $w(x_1,..., x_n)$ in $n$ letters
    $x_1,...,x_n,$ and their inverses, such that if $x_i$ are replaced
    by arbitrary elements in $G$ then the word $w$ becomes~$1$.

In particular this applies to the collection of all solvable
groups of class at most $m\in \N$, and to the collection of
Burnside groups with a uniform bound on the order of elements.
\item The collection of uniform (or cocompact) lattices in semisimple
groups of rank at least $2$ and at most $m\in \N$ is uniformly
unconstricted \cite{KleinerLeeb:buildings}.
\item Every finite collection of unconstricted metric spaces is
uniformly unconstricted, as is, more generally, every collection
of unconstricted metric spaces containing only finitely many
isometry classes.
\end{enumerate}

\medskip

\begin{rmk}
Uniform unconstrictedness is a quasi-isometry invariant in the
following sense. Consider two collections of metric spaces
$\B,\B'$ which are uniformly quasi-isometric, meaning that there
are constants $L \ge 1$ and $C \ge 0$ and a bijection between
$\B,\B'$ such that spaces that correspond under this bijection are
$(L,C)$-quasi-isometric. It follows that $\B$ is uniformly
unconstricted if and only if $\B'$ is uniformly unconstricted.
\end{rmk}

One of the main interests in (uniformly) unconstricted metric
spaces resides in their rigid behavior with respect to
quasi-isometric embeddings into ATG metric spaces.

\begin{thm}[Dru\c{t}u--Sapir \cite{DrutuSapir:TreeGraded}]\label{cutp1}
Let $X$ be ATG with respect to a collection of subsets $\aaa$. Let
$\B$ be a collection of uniformly unconstricted metric spaces. For
every $(L,C)$ there exists $M$ depending only on $L$, $C$, $X$,
$\aaa$ and $\B$, such that for every $(L,C)$--quasi-isometric
embedding $\q $ of a metric space $B$ from $\B$ into $X$, $\q (B)$
is contained in an $M$-neighborhood of a peripheral subset $A\in
\aaa$.
\end{thm}

\section{Non-relative hyperbolicity and quasi-isometric
rigidity}\label{snrh}

In the particular case when all the metric spaces in $\B$ are
finitely generated groups endowed with word metrics, Theorem
\ref{cutp1} can be greatly improved: its conclusion holds when
$\B$ is the collection of \emph{all} NRH groups.

\begin{thm}\label{qinonrelhyp}
Let $(X,\dist_X)$ be ATG with respect to a collection $\aaa$ of
subsets. For every $L\geq 1$ and $C\geq 0$ there exists $R=R(L,C,
X,\aaa)$ such that the following holds. If $(G, \dist)$ is an NRH group
 endowed with a word metric, and
 $\q\co (G,\dist ) \to (X ,\dist_X)$ is an $(L,C)$-quasi-isometric
embedding, then $\q (G)$ is contained in $\nn_R (A)$ for some
$A\in \aaa$.
\end{thm}

\begin{rmk}
    The first result of this kind appeared in Schwartz's proof of the
    classification of non-uniform lattices in rank one semisimple
    Lie groups \cite{Schwartz:RankOne}. In that case,
    one of the key technical steps is showing that
    any quasi-isometry of a neutered space coarsely preserves the
    collection of boundary horospheres. To do this he proved the
    ``Quasi-flat Lemma'' which, reformulated in the language of this
    paper, states that the quasi-isometric image of an
    unconstricted metric space into a neutered space must stay in a
    uniformly bounded neighborhood of a single boundary horosphere.

    This theorem was later generalized by Dru\c{t}u--Sapir
    \cite{DrutuSapir:TreeGraded} who
    kept the unconstricted hypothesis on the domain, but
    replaced the hypothesis that the image is in a
    neutered space by only assuming relative hyperbolicity of
    the target space.
\end{rmk}

\begin{rmk} Theorem~\ref{qinonrelhyp} also holds in the case that
    $G$ is replaced by a metric space which is not ATG. In this
    case though, the constant $R$ will additionally depend on the
    choice of metric space and the choice of quasi-isometry.
\end{rmk}

\begin{rmk}\label{oneend}
By Stallings' Ends Theorem \cite{Stallings:ends} a finitely
generated group has more than one end if and only if it splits
nontrivially as an amalgamated product or HNN-extension with
finite amalgamation. A group which splits in this manner is
obviously hyperbolic relative to its vertex subgroups.
Consequently if a group is NRH then it is one-ended.
\end{rmk}

\begin{rmk}
 A result similar to Theorem \ref{qinonrelhyp} has been
    obtained in \cite[$\S 3$]{PapasogluWhyte:ends}, for $G$ a
    one-ended group and $X$ the fundamental group of a graph of
    groups with finite edge groups.
    Although NRH groups are one-ended,
    the hypothesis in Theorem
    \ref{qinonrelhyp} cannot be weakened to ``$G$ a one-ended group,'' as
    illustrated by the case when $G=X$ and $G$ is the
    fundamental group of a finite volume real hyperbolic manifold.
\end{rmk}

Before proving Theorem~\ref{qinonrelhyp}, we state some
consequences of it, and give a list of examples of NRH groups.

\begin{cor}\label{c1}
Let $G$ be an infinite group which admits an
$(L,C)$-quasi-isometric embedding into a geodesic metric space $X$
which is \atgrs a collection of subsets $\aaa$. Then either $\q
(G)$ is contained in $\nn_R(A)$ for some $A\in \aaa$ and
$R=R(L,C,X,\aaa)$ or $G$ is relatively hyperbolic.
\end{cor}

Another consequence is a new proof of the following which was first
established in \cite[Theorem 1.8]{DrutuSapir:TreeGraded}.

\begin{cor}[see also \cite{DrutuSapir:TreeGraded}, Theorem 1.8]\label{c2}
Let $G$ be a finitely generated group hyperbolic relative to
$\H=\{ H_1, ..., H_m\}$. Let $H$ be an undistorted finitely
generated subgroup of $G$. Then either $H$ is contained in a
conjugate of $H_i$, $i\in \{1,2,...,m\}$, or $H$ is relatively
hyperbolic.
\end{cor}

Perhaps the most important consequence of Theorem
\ref{qinonrelhyp} is the following  quasi-isometric rigidity
theorem for groups hyperbolic relative to NRH subgroups.

\begin{thm}\label{sqir}
Let $G$ be a finitely generated group which is hyperbolic relative
to a finite family of finitely generated subgroups $\H$ such that
each $H\in\H$ is not relatively hyperbolic. If a group $G'$ is
quasi-isometric to $G$ then $G'$ is hyperbolic relative to
$\H'=\{H_1,...,H_m\}$, where each $H_i$ is quasi-isometric to some
$H\in\H$.
\end{thm}
\proof The proof is almost identical to the proof of Theorem~5.13
in \cite{DrutuSapir:TreeGraded}. Indeed let $X=G$ and let $\aaa=
\{gH\; :\; g\in G/H \mbox{ and } H\in\H \}$. The pair $(X,\aaa )$
satisfies all the hypotheses of Theorem~5.13 in
\cite{DrutuSapir:TreeGraded}, except (1). Still, hypothesis (1) is
used in that proof only to ensure that for every quasi-isometry
constants $L\geq 1$ and $C\geq 0$ there exists a constant
$M=M(L,C, X,\aaa)$ such that for every $A\in\aaa$ and for every
$(L,C)$--quasi-isometric embedding $\q \co A\to X$ there exists
$B\in\aaa$ for which $\q (A)\subset \nn_{M}(B)$. In our case, each
$H\in\H$ is known to be undistorted since it is a peripheral
subgroup (see, for instance, \cite{DrutuSapir:TreeGraded} for
details). Thus, the hypothesis that all $H\in\H$ are NRH implies
via Theorem~\ref{qinonrelhyp} that for every $L\geq 1$ and $C\geq
0$ there exists a constant $M$ as above depending only on $L$,
$C$, and the undistorsion constants of each $H$ in~$G$.\endproof

\medskip

In view of Theorems \ref{qinonrelhyp} and \ref{sqir}, it becomes
interesting to consider examples of NRH groups. We do this below.
In Section \ref{section:network} we give a procedure allowing one to
build NRH groups from smaller NRH groups (see Proposition
\ref{conjb}).

\medskip

\textit{Examples of NRH groups:}
\begin{enumerate}
    \item[\textbf{(I)}] Non-elementary groups without free
    non-abelian subgroups.
This follows from the fact that non-elementary relatively
hyperbolic groups contain a free non-abelian subgroup.

The class of groups without free non-abelian subgroups
contains the non-elementary amenable groups,
but it is strictly larger than that class; indeed, a well known
question attributed to J. von Neumann \cite{vonNeumann} is whether
these two classes
coincide (this is known as \textit{the von Neumann problem}). The
first examples of non-amenable groups without free non-abelian
subgroups were given in \cite{Olshanskii:vonNeumann}. Other
examples were later given in \cite{Adyan:randomwalks} and in
\cite{OlshanskiiSapir:vonNeumann}.
    \item[\textbf{(II)}] Non-elementary groups with infinite
    center.
    Indeed, if $G$ is hyperbolic then its center is finite. Assume that $G$
    is relatively hyperbolic with respect to
    $H_1,...,H_m$ and at least one $H_i$ is infinite (otherwise
    $G$ would be hyperbolic). Since $G\neq H_i$ there exists a
    left coset $gH_i\neq H_i$. For every $z\in Z(G)$, $H_i$ and
    $zH_i=H_iz$ are at Hausdorff distance at most $\dist (1,z)$.
    This and Theorem \ref{tgi}, $(\alpha_1)$, imply that $zH_i = H_i$, thus
    $Z(G) \subset H_i$. Similarly it can be proved that $Z(G) \subset gH_ig\inv$. If
    follows that $Z(G)\subset gH_ig\inv\cap H_i$, hence that it
     is finite (see for instance \cite[Lemma 4.20]{DrutuSapir:Splitting}).
    \item[\textbf{(III)}] Unconstricted groups.
    \item[\textbf{(IV)}] Inductive limits of small cancellation
    groups (see Section \ref{exnrh}).
\end{enumerate}

\bigskip

The remainder of this section provides the proof of
Theorem~\ref{qinonrelhyp}, thus we let $(X,\dist_X), \aaa, L$, $C,
\q$ and $G$ be as in the statement of the theorem.  We will
proceed by using the quasi-isometric embedding $\q$ to construct
an \atgs structure on $(G,\dist)$.

In order to produce an \atgs structure on $G$ we first search for
a constant $\tau$ such that the following set is non-empty:
\begin{equation}\label{atau}
\aaa_\tau = \{ A\in \aaa \; ;\; \nn_{\tau }(A)\cap \q (G)\neq
\emptyset\}\, .
\end{equation}
Then, for every $A\in \aaa_\tau$ we consider the pre-image
$B_A=\q\inv \left(\nn_{\tau}(A)\right)\, $ and the set
\begin{equation}\label{btau}
\bb_\tau =\left\{ B_A \; ;\; A\in \aaa_\tau \right\}\, .
\end{equation}
For an appropriate choice of $\tau$, we will show that the
collection $\bb_\tau$ defines an \atgs structure on $(G,\dist )$.
We begin with the following lemmas which will allow us to choose
the constant $\tau$.

\begin{lem}[\cite{DrutuSapir:TreeGraded}, Theorem 4.1 and Remark
    4.2, (2)]\hfill
    \begin{enumerate}
    \item[\textbf{(a)}] \label{m1}
    There exists $M'>0$ such that for every
    $(L,C)$--quasi-geodesic $\p\co [0,\ell ] \to X$ and every $A\in
    \aaa$ satisfying $\p (0)\, ,\, \p(\ell )\in
    \nn_{\ell/3L}(A)$, the tubular neighborhood $\nn_{M'}(A)$
    intersects $\p \left( \left[ 0,\ell \right] \right)$.
    \item[\textbf{(b)}] \label{m2}
    For every $\sigma \geq 1$ and $\nu\geq 4\sigma$ there exists
    $\vartheta_0$ satisfying the following: for every
    $\vartheta \geq \vartheta_0$
    there exists $\chi$ with the property that every hexagon with
    $(L,C)$-quasi-geodesic edges which is $(\vartheta , \sigma ,
    \nu)$--fat is contained in $\nn_\chi (A)$ for some $A\in \aaa$.
    \end{enumerate}
\end{lem}

\begin{lem}\label{qfat}
Let $\p\co Y\to X$ be an $(L,C)$--quasi-isometric embedding. Let
$\sigma=4L^2+L\geq 1$, $\nu=4\sigma$ and $\vartheta\geq C$. If $P$
is a $(2L\vartheta , \nu+1)$--fat geodesic hexagon, then $\p (P)$
is a hexagon with $(L,C)$--quasi-geodesic edges which is $(
\vartheta , \sigma , \nu)$--fat.
\end{lem}
\proof $\mathbf{(F_1)}$ Let $\gfrak$ be an edge of $P$, of
endpoints $x,y$. Let $x\in \p (\gfrak)\setminus \nn_{\sigma
\vartheta }\left( \left\{ \p(x), \p(y) \right\} \right)$. Then
$x=\p(t)$ with $t\in \gfrak$ at distance at most
$\frac{1}{L}\sigma \vartheta -C$ from $x$ and $y$. Since
$\frac{1}{L}\sigma \vartheta -C = (4L+1)\vartheta -C \geq
4L\vartheta$, property $(F_1)$ for $P$ implies that $t$ is at
distance at least $2L \vartheta$ from any edge $\p\neq \gfrak$ of
$P$. Then $\p(x)$ is at distance at least $\frac{1}{L} 2L\vartheta
-C=2\vartheta-C\geq \vartheta$ from $\q (\p)$.

\medskip

$\mathbf{(F_2)}$  Let $v$ be an arbitrary vertex of $P$. Property
$(F_2)$ for $P$ grants that $\dist \left( v , \oo_v \right(P))\geq
(\nu+1)(2L\vartheta)$, hence $\dist \left( \p(v) , \oo_{\p(v)}
\right(\p (P))) \geq \frac{1}{L}(\nu+1)(2L\vartheta) -C =
2(\nu+1)\vartheta -C \geq \nu \vartheta$.
\endproof

For the remainder of the proof, we fix the following constants:
\begin{itemize}
    \item  $\sigma$ and $\nu$ as in Lemma \ref{qfat};
    \item if $\vartheta_0$ is the constant provided by Lemma \ref{m2}
    for $\sigma$ and $\nu$ above, it is no loss of generality to assume further that
    $\vartheta_0\geq C$;
    \item let $\vartheta = 2L \vartheta_0$;
    \item let $\chi$ the constant given by Lemma \ref{m2} for
     $\vartheta_0$;
    \item  $\tau = \max \left( \chi , M'\right)$, where $M'$ is
the constant from Lemma \ref{m1}.
\end{itemize}

\medskip

If $G$ does not contain a $(\vartheta, \nu+1)$--fat geodesic
hexagon  or if all such hexagons have uniformly bounded diameter,
then $G$ is hyperbolic by Corollary~4.20 in \cite{Drutu:RelHyp}.
This contradicts our hypothesis on $G$. We may thus henceforth
assume that for every $\eta >0$, the space $G$ contains a
$(\vartheta ,\nu+1)$--fat  geodesic hexagon of diameter at
least~$\eta$. For every such hexagon $P$, Lemma \ref{qfat} and the
above choice of constants imply that $\q(P)\subset
\nn_{\chi}(A)\subset \nn_{\tau }(A)$ for some $A\in \aaa$. In
particular the set $\aaa_\tau$ is non-empty.

\begin{lem}\label{ltgr}
The metric space $(G,\dist )$ is asymptotically tree-graded with
respect to the set $\bb_\tau$ defined in~(\ref{btau}).
\end{lem}
\proof We start with the simple remark that if $x\in \nn_t \left(
B_A\right)$ then $\q (x)\in \nn_{Lt+C+\tau }(A)$.

\medskip

According to Theorem \ref{tgi} it suffices to verify conditions
$\mathbf{(\alpha_1)}$, $\mathbf{(\alpha_2)}$,
$\mathbf{(\beta_3)}$.

\medskip

We first establish $\mathbf{(\alpha_1)}$. Let $A,A'\in \aaa_\tau
,\, A\neq A'$, and let $x,y\in \nn_\delta \left( B_A\right)\cap
\nn_\delta \left( B_{A'}\right)$. Then $\q(x)$ and $\q(y)$ are in
$\nn_{L\delta +C+\tau} \left( A\right)\cap \nn_{L\delta +C+\tau }
\left( A'\right)$. Since $X$ is asymptotically tree-graded,
$\diam(\nn_{L\delta +C+\tau} \left( A\right)\cap \nn_{L\delta
+C+\tau } \left( A'\right))=D$ is uniformly bounded. Thus
$$
\dist (x,y)\leq L\left[ \dist_X(\q(x), \q(y))+C \right]\leq L
\left( D +C \right)\, .
$$

\medskip

We prove $\mathbf{(\alpha_2)}$ for $\varepsilon = \frac{1}{6L^2}$
and $M=\left( 3L + \frac{1}{L} \right)(C+\tau)$. Let $\gfrak\co
[0,\ell ]\to G$ be a geodesic with endpoints in $\nn_{\varepsilon
\ell} \left( B_A \right)$ for some $A\in \aaa_\tau$. Then $\q
\circ \gfrak$ is an $(L,C)$-quasi-geodesic with endpoints in
$\nn_{L \varepsilon \ell +C+\tau } \left( A \right)$. If $C+\tau
\geq  \frac{\ell}{6L}$, that is $\ell \leq 6L (C+\tau )$ then
$\gfrak \subset \overline{\nn} _{3L(C+\tau)}\left( \left\{
\gfrak(0),\gfrak(\ell) \right\} \right)\subset \nn_{\left( 3L +
1/L \right)(C+\tau)} \left( B_A \right)$. If $C+\tau <
\frac{\ell}{6L}$ then Lemma \ref{m1} implies that $\q \circ \gfrak
\left( \left[0,\ell \right] \right)$ intersects $\nn_{M'}(A)$. It
follows that $\gfrak \left( \left[0,\ell \right] \right)$
intersects $B_A$.

\medskip

We prove $\mathbf{(\beta_3)}$ for $(\vartheta , \nu +1)$ as above
and for $\chi =0$. Let $P$ be a $(\vartheta , \nu +1)$--fat
geodesic hexagon in $G$. Then by Lemma \ref{qfat}, $\q (P)$ is a
$(\vartheta_0 , \sigma , \nu)$--fat hexagon with
$(L,C)$--quasi-geodesic edges. Lemma~\ref{m2} implies that $\q
(P)$ is contained in $\nn_{\chi}(A)$. It follows that $A\in
\aaa_\tau$ and that $P\subset B_A$.\endproof

\medskip

Let $M$ be the maximum between the constant from $(\beta_2)$ and
the constant $\chi$ from $(\beta_3)$, for $(G, \bb_\tau)$. Note
that the constants in $(\beta_2)$ and $(\beta_3)$ for $(G,
\bb_\tau)$ can be obtained from the constants in the same
properties for $(X, \aaa)$, as well as from $\tau$, $L$ and $C$.
Consequently $M=M(X,\aaa , L,C)$.

Lemma~\ref{ltgr}, Theorem \ref{t3} and the hypothesis that $G$ is
NRH imply that $G\subset \nn_M \left( B_A \right)$ for some $A\in
\aaa_\tau$. Hence $\q (G)\subset \nn_{L M +C +\tau }(A)$,
completing the proof of Theorem~\ref{qinonrelhyp}.

\section{Networks of subspaces}\label{section:network}

We begin by defining the notions of networks of subspaces and of
subgroups.

\begin{defn}\textbf{(network of subspaces).}

Let $X$ be a metric space and $\lll$ a collection of subsets of
$X$. Given $\tau \ge 0$ we say that $X$ is a
\textit{$\tau$--network with respect to the collection~$\lll$} if
the following conditions are satisfied:

\smallskip

\begin{enumerate}
\item[$(\bn_1)$] $X = \bigcup_{L \in \lll} \nbhd_\tau(L)$;

\smallskip

\item[$(\bn_2)$] Any two elements $L,L'$ in $\lll$ can be
\textit{thickly connected in} $\lll$: there exists a sequence,
$L_{1}=L,L_{2}, \ldots , L_{n-1},L_{n}=L'$, with $L_{i}\in\lll$
and with $\diam(\nbhd_\tau(L_{i})\cap \nbhd_\tau(L_{i+1}))=\infty$
for all $1\leq i<n$.
\end{enumerate}
\end{defn}

We now define a version of the above notion in the context of
finitely generated groups with word metrics. Recall that a
finitely generated subgroup $H$ of a finitely generated group $G$
is {\textit{undistorted}} if any word metric of $H$ is
bi-Lipschitz equivalent to a word metric of $G$ restricted to $H$.

\begin{defn}[\textbf{algebraic network of subgroups}]\label{dgnetwork}
Let $G$ be a finitely generated group, let $\H$ be a finite
collection of subgroups of $G$ and let $M>0$. The group $G$ is an
$M$--\emph{algebraic network with respect to $\H$} if:

\medskip

\begin{enumerate}
\item[$(\ban_0)$] All subgroups in $\H$ are finitely generated and
undistorted in $G$.

\medskip

\item[$(\ban_1)$] There is a
finite index subgroup $G_1$ of $G$ such that $G\subset \nbhd_M
(G_1)$ and such that a finite generating set of $G_1$ is contained
in $\bigcup_{H\in \H} H$.

\medskip

\item[$(\ban_2)$]  Any two subgroups $H,H'$ in $\H$
can be \textit{thickly connected in }$\H$: there exists a finite
sequence $H=H_{1}, \ldots,H_{n}=H'$ of subgroups in $\H$ such that
for all $1\leq i<n$, $H_{i}\cap H_{i+1}$ is infinite.
\end{enumerate}
\end{defn}

\bigskip

\begin{prop}\label{inetw}
If a finitely generated group $G$ is an $M$--algebraic network
with respect to $\H$ then it is an $M$--network with respect to
the collection of left cosets
$$\lll=\left\{\, gH : g\in G_1\, ,\, H\in \H\, \right\}\, .$$
\end{prop}

\begin{pf} Property $(\bn_1)$ is trivial. We prove property $(\bn_2)$.
Since it is equivariant with respect to the action of $G$ it
suffices to prove it for $L=H$ and $L'=gH'$, $H,H'\in \H$ and
$g\in G_1$. Fix a finite generating set $S$ of the finite index
subgroup $G_1$ of $G$ so that $S \subset \bigcup_{H\in \H } H$;
all lengths in $G_1$ will be measured with respect to this
generating set. We argue by induction on $|g|=|g|_S$. If $|g|=1$,
then $g\in S$. By hypothesis, $g$ is contained in a subgroup
$\widetilde{H}$ in $\H$. We take a sequence $H=H_1,H_2,...,
H_k=\widetilde{H}$ as in $(\ban_2)$, and a similar sequence
$\widetilde{H}=\overline{H}_1,\overline{H}_2,...,
\overline{H}_m=H'$. Then the sequence $$H=H_1,H_2,...,
H_k=\widetilde{H}=g\widetilde{H}=g\overline{H}_1,
g\overline{H}_2,..., g\overline{H}_m=gH'$$ satisfies the
properties in $(\bn_2)$. We now assume the inductive hypothesis
that for all $g\in G_1$ with $|g|\leq n$ and all $H,H'\in \H$ the
cosets $H$ and $gH'$ can be connected by a sequence satisfying
$(\bn_2)$ with $\tau = M$. Take $g\in G_1$ such that $|g|=n+1$;
thus $g=\hat{g} s$, where $s\in S$ and $\hat{g}\in G_1$,
$|\hat{g}|=n$. By hypothesis there exists some $\widetilde{H}\in
\H$ containing $s$. Take arbitrary $H,H'\in \H$. In order to show
that $H$ and $gH'$ can be connected by a good sequence it suffices
to show, by the inductive hypothesis, that $\hat{g}\widetilde{H}
=g\widetilde{H}$ and $gH'$ can be connected by a good sequence.
This holds because $\widetilde{H}$ and $H'$ can be so connected,
according to $(\ban_2)$.
\end{pf}

One of the reasons for which one can be interested in the notion
of network of groups is that it represents a way of building up
NRH groups. More precisely the following holds:

\begin{prop}\label{conjb}
Let $G$ be a finitely generated group which is an $M$--algebraic
network with respect to $\H$, such that each of the subgroups in
$\H$ is not relatively hyperbolic.

If $G$ is an undistorted subgroup of a group $\Gamma$ hyperbolic
relative to $\widetilde{H}_1,.., \widetilde{H}_m$, then $G$ is
contained in a conjugate of some subgroup $\widetilde{H}_i, i\in
\{1,...,m\}$.

In particular $G$ is not relatively hyperbolic.
\end{prop}

\begin{pf} According to Corollary \ref{c2}, any subgroup
$H\in\H$ is contained in the conjugate of some $\widetilde{H}_i,
i\in \{1,...,m\}$. Since distinct conjugates of subgroups
$\widetilde{H}_i$ have finite intersections, it follows from
$(\ban_2)$ that all subgroups in $\H$ are in the same conjugate
$\gamma \widetilde{H}_i \gamma\inv$. Hence, condition $(\ban_1)$
implies that $G$ has a finite index subgroup $G_1$ which is
completely contained in the same conjugate $\gamma \widetilde{H}_i
\gamma\inv$. Given $M$ the constant in $(\ban_1)$, for any $g\in
G$, $gG_1 g\inv \subset \nn_M (G_1)\subset \nn_M (\gamma
\widetilde{H}_i \gamma\inv)$. It follows that $g(\gamma
\widetilde{H}_i\ \gamma\inv)g\inv \cap \nn_M (\gamma
\widetilde{H}_i \gamma\inv)$ has infinite diameter. From this it
can be deduced, by \cite[Lemma~2.2]{MSW:QTTwo}, that $g(\gamma
\widetilde{H}_i \gamma\inv)g\inv \cap \gamma \widetilde{H}_i
\gamma\inv$ is also infinite. This implies that the two conjugates
coincide and thus $g\in \gamma \widetilde{H}_i \gamma\inv$. We
have thereby shown that $G<\gamma \widetilde{H}_{i}
\gamma\inv$.\end{pf}

In Proposition \ref{conjb} the hypotheses of undistortedness (of
$G$ in $\Gamma$ and of  every subgroup $H\in \H$ in $G$) can be
removed, if the hypothesis ``all subgroups in $\H$ are NRH'' is
strengthened to ``all subgroups in $\H$ are non-elementary and
without free non-Abelian subgroups''. The latter condition implies
the former but they are not equivalent: for instance uniform
lattices in semisimple groups of rank at least two are
unconstricted hence NRH and they have many non-Abelian free
subgroups.

Thus, the following statement, generalizing the main result of
\cite{AASh:RelHypMCG}, holds:

\begin{prop}\label{conjb2}
Let $G$ be a finitely generated group with a finite collection
$\H$ of finitely generated subgroups satisfying $(\ban_1)$ and
$(\ban_2)$. Assume moreover that all $H\in\H$ are non-elementary
and do not contain free non-Abelian subgroups.

If $G$ is a subgroup of a group $\Gamma$ hyperbolic relative to
$\widetilde{H}_1,.., \widetilde{H}_m$, then $G$ is contained in a
conjugate of some subgroup $\widetilde{H}_i, i\in \{1,...,m\}$.

In particular $G$ is not relatively hyperbolic.
\end{prop}

\begin{pf} We use the Tits alternative in relatively hyperbolic
groups: a subgroup in $\Gamma$ is either virtually cyclic,
parabolic (i.e.\ contained in a conjugate of some subgroup
$\widetilde{H}_i$), or it contains a free non-Abelian subgroup;
the proof follows from \cite{Tukia:convergence} and
\cite{Bowditch:RelHyp}. Hence, with our hypotheses, any subgroup
$H\in\H$ is parabolic. The rest of the proof is identical to the
proof of Proposition \ref{conjb}. \end{pf}

\medskip


\section{Relative hyperbolicity and Dunwoody's inaccessible
group}\label{section:Dunwoody}

Having a quasi-isometric rigidity theorem for relatively
hyperbolic groups whose peripheral subgroups are not relatively
hyperbolic, one might think to ask:

\begin{qn}\label{qacc}
Given a finitely generated relatively hyperbolic group $G$, is $G$
hyperbolic relative to some finite collection of subgroups
\emph{none of which are relatively hyperbolic}?
\end{qn}

\medskip

\begin{rmk}\label{rterminal}
Note that if $G$ is hyperbolic relative to $\{ H_i\; ;\;
i=1,2,..,m\}$ and if each $H_i$ is hyperbolic relative to $\{
H_i^j\; ;\; j=1,2,..,n_i\}$ then $G$ is hyperbolic relative to $\{
H_i^j\; ;\; j=1,2,..,n_i\, ,\, i=1,2,..,m\}$. Examples where such
process never terminates are easily found (for instance when $G$
is a free non-Abelian group with $H_i$ finitely generated
non-Abelian subgroups and $H_i^j$ finitely generated non-Abelian
subgroups of $H_i$.). Still, one might ask if in every
relatively hyperbolic group there exists a terminal point
for the process above (like $H=\{1\}$ in the case of a free
group). This is the meaning of Question \ref{qacc}.
\end{rmk}

\medskip

We answer this question in the negative, using Dunwoody's example
$J$ of an inaccessible group \cite{Dunwoody:Notaccessible}:

\begin{prop}
\label{PropDunwoody} Dunwoody's group $J$ is relatively
hyperbolic. If $J$ is hyperbolic relative to a finite collection
of subgroups $A_1,\ldots,A_I$, at least one of the subgroups
$A_1,\ldots,A_I$ is relatively hyperbolic.
\end{prop}

This proposition shows that $J$ satisfies a kind of ``relatively
hyperbolic inaccessibility'': whenever $J$ is written as a
relatively hyperbolic group, one of the peripheral subgroups $A$
is also relatively hyperbolic and so $A$ can be replaced by its
list of peripheral subgroups, giving a new relatively hyperbolic
description of $J$; this operation can be repeated forever, giving
an infinite sequence of finer and finer relatively hyperbolic
descriptions of $J$.

First we review Dunwoody's construction of $J$. Let $H$ be the
group of permutations of $\Z$ generated by the transposition $t =
(0,1)$ and the shift map $s(i)=i+1$. Each element $\sigma \in H$
agrees outside a finite set with a unique power $s^p$, and the map
$\sigma \mapsto p$ defines a homomorphism $\pi \from H \mapsto \Z$
whose kernel denoted $H_\omega$ is the group of finitely supported
permutations of $\Z$. Let $H_i \subset H_\omega$ be the group of
permutations supported on
$[-i,i]=\{-i,-i+1,\ldots,0,1,\ldots,i-1,i\}$, so $H_\omega =
\union_{i=0}^\infinity H_i$. Let $V$ be the group of all functions
from $\Z$ to $\Z_2=\{\pm 1\}$ with finite support and the usual
group law. Let $V_i$ be the subgroup of all such maps with support
$[-i,i]$. Let $z_i \in V_i$ be the map defined by $z_i(n)=-1$ if
and only if $n \in [-i,i]$. The group $H_i$ acts on the left of
$V_i$ by $^h v(j) = v(h^{-1} (j))$, and so we can form the
semidirect product $G'_i = V_i \semidirect H_i$, each of whose
elements can be written uniquely as $vh$ with $v \in V_i$ and $h
\in H_i$, and the group law is $(v_0h_0) \cdot (v_1 h_1) = (v_0 \,
{}^{h_0}v_1) (h_0h_1)$. The element $z_i$ is central in $G'_i$ and
so we have a direct product subgroup $K_i = \<z_i\> \cross H_i
\approx \Z/2 \cross H_i \subgroup G'_i$. For $i=1,2,\ldots$ choose
$G_i$ to be an isomorphic copy of $G'_i$, with the $G_i$ pairwise
disjoint. The group $K_i$ being a subgroup of $G'_i$ and of
$G'_{i+1}$, we may identify $K_i$ with its images in $G_i$ and
$G_{i+1}$, which defines the following graph of groups whose
fundamental group is denoted $P$:
\begin{equation}
\xymatrix{ G_1 \ar@{-}[r]^{K_1} & G_2 \ar@{-}[r]^{K_2} & G_3
\ar@{-}[r]^{K_3} & G_4 \ar@{-}[r]^{K_4} & \ldots } \label{Pgog}
\end{equation}
We shall need below the following equation which can be regarded
as taking place within $G_{i+1}$:
\begin{equation}
K_i \intersect K_{i+1} = H_i \label{HiIntersection}
\end{equation}
Collapsing all edges in (\ref{Pgog}) to the right of the one
labeled $K_n$ produces another decomposition of $P$ as the
fundamental group of the graph of groups
\begin{equation}\xymatrix{
G_1 \ar@{-}[r]^{K_1} & G_2 \ar@{-}[r]^{K_2} & G_3 \ar@{.}[r] &
G_{n-1} \ar@{-}[r]^{K_{n-1}} & G_n \ar@{-}[r]^{K_n} & Q_n
}\label{QnDef}
\end{equation}
and then collapsing all edges except the one labeled $K_n$ we get
a decomposition $P = P_n *_{K_n} Q_n$. Noting that $P$ contains
$H_1 < H_2 < H_3 < \cdots < H_\omega$, we can form the amalgamated
product
\begin{equation}J = P *_{H_\omega} H
\end{equation}
Since $H_\omega \subset Q_n$, the group $J$ also has the
decomposition
\begin{equation}J = P_n *_{K_n}
    \underbrace{\left( Q_n * _{H_\omega} H \right)}_{J_n}
\label{JnRelHyp}
\end{equation}
Applying~(\ref{QnDef}) and the definition of $P_n$ we obtain a
decomposition of $J$ as the fundamental group of the graph of
groups
\begin{equation}\xymatrix{
G_1 \ar@{-}[r]^{K_1} & G_2 \ar@{-}[r]^{K_2} & G_3 \ar@{.}[r] &
G_{n-1} \ar@{-}[r]^{K_{n-1}} & G_n \ar@{-}[r]^{K_n} & J_n }
\label{JnRelHypTwo}
\end{equation}

From both (\ref{JnRelHyp}) and (\ref{JnRelHypTwo}) we see that $J$
is relatively hyperbolic: in either of these graph of groups
presentations, each edge group is finite and includes properly
into both adjacent vertex groups, and so $J$ is hyperbolic
relative to the vertex groups. This proves the first clause of
Proposition~\ref{PropDunwoody}.

To prepare for the rest of the proof we need some additional facts
about the group $H$.
\begin{itemize}
\item $H$ is the intersection of the nested family of subgroups
$J_1 > J_2 > J_3 > \cdots$
\end{itemize}
To prove this, since $J_n = Q_n *_{H_\omega} H$, it suffices to
show that the intersection of the nested family $Q_1 > Q_2 > Q_3 >
\cdots$ equals $H_\omega$. Consider an element $x \in P=Q_0$ that
is contained in this intersection. Since $Q_0$ is generated by its
subgroups $G_1,G_2,\ldots$, the element $x$ can be written as a
word $w_1$ whose letters are elements of $G_1,G_2,\ldots,G_n$ for
some $n$. Since $x \in Q_1$, $x$ can also be written as a product
of elements of $G_2,G_3,\ldots$. By uniqueness of normal forms in
a graph of groups, any letter of $w_1$ that is in the subgroup
$G_1$ is also in the subgroup $K_1$; each such letter can be
pulled across the $K_1$ edge into the subgroup $G_2$, and so $x$
can be written as a word $w_2$ whose letters are elements of
$G_2,\ldots,G_n$. Continuing inductively in this fashion, we see
that $x \in G_n$. Going one more step, since $x \in Q_n$, it can
be written as a product of elements of
$G_{n+1},G_{n+2},G_{n+3},\ldots$, and so by uniqueness of normal
forms we have $x \in K_n \subgroup G_{n+1}$. And going one more
step again, $x$ can be written as a product of elements of
$G_{n+2},G_{n+3},\ldots$, and so $x \in K_{n+1}$. Applying
(\ref{HiIntersection}) we have $x \in K_n \intersect K_{n+1} = H_n
\subgroup H_\omega$.

Next we need:
\begin{itemize}
\item $H$ satisfies the hypotheses of Proposition \ref{conjb2}.
\end{itemize}
To see this, let $H_\even$ be the abelian subgroup of $H$
generated by the transpositions $(2n,2n+1)$, $n \in \Z$, and let
$H_\odd$ be the abelian subgroup generated by the transpositions
$(2n+1,2n+2)$, $n \in \Z$. The squared shift map $s^2$ preserves
each of these subgroups, and so we have subgroups $H_\even
\semidirect \<s^2\>$ and $H_\odd \semidirect \<s^2\>$ of $H$, each
solvable and therefore non-elementary, without free non-Abelian
subgroups. These two subgroups generate the index~2 subgroup
$\pi\inv (2\Z)$, thus $(\ban_1)$ is satisfied. Also $\<s^2\>$ is
contained in both subgroups, whence $(\ban_2)$.

\me

Now we prove the second clause of Proposition~\ref{PropDunwoody}.
Arguing by contradiction, suppose that $J$ is relatively
hyperbolic with respect to peripheral subgroups $L_1,\ldots,L_m$
none of which is relatively hyperbolic. By Proposition
\ref{conjb2}, the group $H$ must be contained in some conjugate of
some $L_i$, so we have $H \subgroup L'_i = g L_i g\inv$ for some
$g \in J$. Since $H$ is infinite, $L'_i$ is infinite.  By
combining Corollary~\ref{c2} with the relatively hyperbolic
description (\ref{JnRelHypTwo}), the NRH subgroup $L'_i$ must be
contained in a conjugate of one of $G_1,\ldots,G_n,J_n$, but only
$J_n$ is infinite and so $L'_i \subgroup h J_n h\inv$ for some $h
\in J$. We therefore have $H \subgroup J_n \intersect h J_n
h\inv$, and so by malnormality $J_n = h J_n h\inv$. Thus $L'_i
\subgroup J_n$ for all $n$, and so $L'_i \subgroup H$. We have
therefore proved that $L'_i = H$, and so $J$ is hyperbolic
relative to a collection of subgroups that includes $H$.

Now note that $H \cap z_n H z_n^{-1}$ contains $H_n$ whose
cardinality goes to $+\infinity$ as $n \to +\infinity$. Since $H$
is a peripheral subgroup of $J$, the intersection of $H$ with its
distinct conjugates has uniformly bounded cardinal. Thus for $n$
large enough we have $H=z_nHz_n^{-1}$. In particular $H$ and
$z_nH$ are at finite Hausdorff distance $|z_n|$, which together
with the fact that $H$ is infinite and with property $(\alpha_1)$
imply that $H=z_nH$, hence that $z_n\in H$, a contradiction.

\section{Thick spaces and groups}\label{section:thick}

A particular case of NRH groups are those obtained by using the
construction in Proposition \ref{conjb} inductively, with
unconstricted groups as a starting point. This particular case of
groups are the thick groups. We begin by introducing the notion of
thickness in the general metric setting.

\begin{defn}\textbf{(metric thickness and uniform thickness).}
\begin{itemize}
\item[\textbf{(M$_1$)}] A metric space is called \textit{thick of
order
zero} if it is unconstricted. A family of metric spaces is
\textit{uniformly thick of order
zero} if it is uniformly unconstricted.

\medskip

\item[\textbf{(M$_2$)}] Let $X$ be a metric space and $\lll$ a
collection of subsets of $X$. Given $\tau \ge 0$ and $n\in\N$ we
say that $X$ is \textit{$\tau$--thick of order at most~$n+1$ with
respect to the collection~$\lll$} if $X$ is a $\tau$-network with
respect to $\lll$, and moreover:

\smallskip

\begin{enumerate}
\item[$(\theta)$] when the subsets in $\lll$ are endowed with the
restricted metric on $X$, then the collection $\lll$ is uniformly
thick of order at most $n$.
\end{enumerate}

\medskip

\noindent
We say $X$ is \emph{thick of order at most~$n$} if it is
\textit{$\tau$--thick of order at most~$n$}
 with respect to some collection $\lll$ for some~$\tau$.
Further, $X$ is said to be $\tau$--\emph{thick of order~$n$ (with
respect to the collection~$\lll$)} if it is $\tau$--\emph{thick of
order at most~$n$ (with respect to the collection~$\lll$)} and for
no choices of $\tau$ and $\lll$ is it thick of order at
most~$n-1$. When the choices of $\lll$, $\tau$, and~$n$ are
irrelevant, we simply say that $X$ is \emph{thick}.

\medskip

\item[\textbf{(M$_3$)}] A family $\left\{ X_i \mid i\in I \right\}$
of metric spaces is \textit{uniformly thick of order at most
}$n+1$ if the following hold.
\begin{itemize}
\item[$(\upsilon \theta_1)$] There exists $\tau >0$ such that
every $X_i$ is $\tau$--thick of order at most $n+1$ with respect
to a collection~$\lll_i$ of subsets of it;
\item[$(\upsilon \theta_2)$] $\bigcup_{i\in I}\lll_i$ is uniformly
thick of order at most $n$, where each $L\in \lll_i$ is endowed
with the induced metric.
\end{itemize}
\end{itemize}
\end{defn}

\begin{rmk}\label{thickqiinvariant}
Thickness is a quasi-isometry invariant in the following sense.
Let $X,X'$ be metric spaces, let $\q \colon X \to X'$ be a
$(L,C)$--quasi-isometry, let $\lll$ be a collection of subsets of
$X$, let $\lll'$ be a collection of subsets of $X'$, and suppose
that there is a bijection $\q_\# \colon \lll \to \lll'$ such that
the subsets $\q(L)$ and $\q_\#(L)$ have Hausdorff distance $\le C$
in $X'$, for each $L \in \lll$. For example, one could simply take
$\lll'=\{ \q(L) \suchthat L \in \lll\}$. If we metrize each space
in $\lll$ or in $\lll'$ by restricting the ambient metric, it
follows that $\lll$ and $\lll'$ are uniformly quasi-isometric, and
so $\lll$ is uniformly unconstricted if and only if $\lll'$ is
uniformly unconstricted. This is the basis of an easy inductive
argument which shows that $X$ is $\tau$--thick of order $n$ with
respect to $\lll$ if and only if $X'$ is $\tau'$--thick of order
$n$ with respect to $\lll'$, where $\tau' = \tau'(L,C,\tau)$.
\end{rmk}

\medskip

We now define a stronger version of thickness in the context of
finitely generated groups with word metrics.

\begin{defn}[\textbf{algebraic thickness}]\label{dgthick}
Consider a finitely generated group $G$.
\begin{itemize}
\item[\textbf{(A$_1$)}] $G$ is called \textit{algebraically thick
of order zero} if it is unconstricted.

\medskip

\item[\textbf{(A$_2$)}] $G$ is called $M$--\emph{algebraically thick
of order at most~$n+1$ with respect to $\H$}, where $\H$ is a
finite collection of subgroups of $G$ and $M>0$, if:

\medskip

\begin{itemize}
    \item $G$ is an $M$-algebraic network with respect to $\H$;

    \medskip

    \item  all subgroups in $\H$ are algebraically thick of order
    at most $n$.
\end{itemize}
\end{itemize}
\end{defn}
$G$ is said to be \emph{algebraically thick
of order~$n+1$ with respect to $\H$}, when $n$ is the smallest value
for which this statement holds.

\medskip

\begin{rmk}\label{thickword} The algebraic thickness property does
not depend on the word metric on~$G$, moreover it holds for any
metric quasi-isometric to a word
metric. Hence in what follows, when mentioning this property for a
group we shall mean that the group is considered endowed with some
metric quasi-isometric to a word metric.
See Section~\ref{SectionMCG} for an example where we use a proper
finite index subgroup $G_1$ to verify thickness.
\end{rmk}

\bigskip

\noindent \textit{Examples:} Examples of groups that are
algebraically thick of order one are provided by mapping class
groups (see Section \ref{SectionMCG} and
\cite{Behrstock:asymptotic}), right angled Artin groups whose
presentation graph is a tree of diameter greater than 2
(Corollary~\ref{CorollaryCATZeroArtin} and
Proposition~\ref{raagcutpoints}), and fundamental groups of graph
manifolds (see Section \ref{Sectiongrm} and \cite{KKL:QI}). An
example of a metric space thick of order one is the Teichm\"uller
space with the Weil-Petersson metric (see Section
\ref{Section:teich} and \cite{Behrstock:asymptotic}). An example
of a group thick of order two is described in
\cite{BehrstockDrutu:thick2}, see also Remark \ref{rthick2}.

\bigskip

\begin{qn}\label{qn:algthickqi} Since the order of metric thickness is a
    quasi-isometry
    invariant (see Remark \ref{thickqiinvariant}), we ask whether
    the order of algebraic thickness is also a quasi-isometry invariant.
\end{qn}

\begin{prop}\label{ithick}\quad
\begin{itemize}
\item[(a)] If a finitely generated group $G$ is $M$--algebraically
thick of order at most $n$ then it is $M$--metrically thick of
order at most $n$. Moreover if $n\geq 1$ and $G$ is
$M$--algebraically thick of order at most $n$ with respect to $\H$
then it is $M$--metrically thick of order at most $n$ with respect
to the collection of left cosets
$$\lll=\left\{\, gH : g\in G_1\, ,\, H\in \H\, \right\}\, .$$
\item[(b)] Let $G_1,G_2,...,G_n$ be finitely generated groups
algebraically thick of order at most $n$. Then any family $\{
X_i\mid i\in I \}$ of metric spaces such that each $X_i$ is
isometric to $G_k$ for some $k\in \{ 1,2,..., n\}$ is uniformly
metrically thick of order at most $n$.
\end{itemize}
\end{prop}

\begin{pf} We prove the proposition inductively on $n$. The
statements
(a) and (b) are true for $n=0$. Suppose that they are true for
all $k\leq n$. We prove them for $n+1$.

\medskip

\noindent\textbf{(a)} Since all groups in $\H$ are undistorted and
algebraically thick of order at most $n$ with respect to their own
word metrics, by Remark \ref{thickword} it follows that they are
algebraically thick of order at most $n$ also when endowed with
the restriction of the metric on $G$. This and (b) for $n$ imply
that $\lll$ is uniformly metrically thick of order at most $n$,
verifying condition $(\theta)$. This and Proposition \ref{inetw}
allow to finish the argument.

\medskip

\noindent\textbf{(b)} Each group $G_i$ is $M_i$--algebraically
thick of order at most $n$ with respect to some collection $\H_i$
of subgroups, where $M_i>0$. Each $H\in \H_i$ is thick of order at
most $n-1$. Property $(\upsilon\theta_1)$ holds for $\{ X_i\mid
i\in I\}$, with the constant $\tau=\max \{ M_i \mid i\in
\{1,2,...,n\} \}$. Each metric space $X_i$, $i\in I$, is isometric
to some $G_k$, hence by (a) it is metrically thick with respect to
the family of isometric images of $\{ gH \mid g\in G_k^1,\, H\in
\H_k\}$, where $G_k^1$ is a finite index subgroup in $G_k$.
Property (b) applied to the finite family of groups
$\bigcup_{k=1}^n \H_k$ yields property $(\upsilon \theta_2)$ for
the family of metric spaces $\{ X_i\mid i\in I \}$.
\end{pf}

A consequence of Proposition \ref{ithick} is that the order of
algebraic thickness is at least the order of metric
thickness. Thus, we ask the following strengthening of
Question~\ref{qn:algthickqi}.

\begin{qn}
For a finitely generated group is the order of algebraic thickness
equal to the order of metric thickness?
\end{qn}

A motivation for the study of thickness is
that it provides a metric obstruction to relative hyperbolicity.
In particular, it gives us
examples to which one can apply Theorem~\ref{qinonrelhyp}.

In the sequel, we shall not mention the collection of
subsets/subgroups with respect to which thickness is satisfied,
when irrelevant to the problem.

\begin{thm}\label{thick:peripheral} Let $\mathcal{X}$ be a
collection of uniformly thick metric spaces, and let $Y$ be a
metric space asymptotically tree-graded with respect to a
collection $\P$ of subsets. Then there is a constant $M=M(L,C,
\mathcal{X}, Y,\P )$ such that for any $X\in \mathcal{X}$ and any
$(L,C)$--quasi-isometric embedding $\q \colon X\into Y$, the image
$\q (X)$ is contained in $\nbhd_{M}(P)$ for some $P\in\P$.
\end{thm}

\begin{pf} We prove the statement by induction on the order of
thickness. If $n=0$, then the family $\mathcal{X}$ is uniformly
unconstricted and the statement follows from Theorem~\ref{cutp1}.
Assume that the statement is true for $n$. We prove it for $n+1$.
Let $\mathcal{X}$ be a collection of metric spaces uniformly thick
of order at most $n+1$. For each $X\in \mathcal{X}$ let $\lll_X$
be the collection of subsets with respect to which $X$ is thick.
The family $\lll=\bigcup_{X\in \mathcal{X}}\lll_X$ is uniformly
thick of order at most $n$. By the inductive hypothesis, there
exists $M=M(L,C, \lll , Y, \P )$ such that for any $L\in\lll$, any
$(L,C)$--quasi-isometric embedding of $L$ into $Y$ is contained
into the radius $M$ neighborhood of a set $P\in \P$. Let $X$ be
any metric space in $\mathcal{X}$ and let $\q : X\to Y$ be an
$(L,C)$--quasi-isometric embedding. For every $L\in \lll_X$, the
subset $\q (L)$ is contained in $\nbhd_{M}(P_L)$ for some
$P_L\in\P$. Further, hypothesis $(\bn_2)$ is satisfied also by the
collection of subsets $\left\{ \q (L) \mid L\in \lll \right\}$.
Theorem \ref{tgi}, $(\alpha_1)$, implies that $P_L$ is the same
for all $L\in \lll$. It follows that $\q \left( \bigcup_{L\in \lll
} L \right)$ is contained in the $M$--neighborhood of $P$.
Properties $(\upsilon\theta_1)$ and $(\bn_1)$ together imply that
$\q (X)$ is contained in the $(M+L\tau +C)$--neighborhood of the
same $P$.
\end{pf}
Taking $Y=X$ this immediately implies:

\begin{cor}\label{crit} If $X$ is a thick metric space, then
$X$ is not asymptotically tree-graded. In particular, if $X$ is a
finitely generated group, then $X$ is not relatively hyperbolic.
\end{cor}

\subsection{NRH groups which are not thick}\label{exnrh}

Thick groups provide an important class of NRH groups. It is
therefore natural to ask whether there exist examples of NRH
groups which are not thick. A construction in
\cite{ThomasVelickovic} (of which a more elaborated version can be
found in \cite[$\S 7$]{DrutuSapir:TreeGraded}) provides an example
of a two-generated group, recursively (but not finitely)
presented, which is NRH and not metrically thick.

\begin{n}
Given an alphabet $A$ and a word $w$ in this alphabet, $|w|$
denotes the length of the word.
\end{n}

\begin{defn}[property $C^*(\lambda)$]
    Let $\free_A$ denote the set of reduced words in an alphabet $A$.
A set $\mathcal{W}\subset\free_A$ which is assumed to be
closed under cyclic permutations and taking inverses, is said to
satisfy {\it property }$C^*(\lambda )$ if the following hold:

\me

\begin{itemize}
\item[(1)] if $u$ is a subword in a word $w\in \mathcal{W}$ so
that $|u|\geq \lambda |w|$ then $u$ occurs only once in $w$;

\me

\item[(2)] if $u$ is a subword in two distinct words $w_1,w_2\in
\mathcal{W}$ then $|u|\leq \lambda \min (|w_1|,|w_2|)$.
\end{itemize}
\end{defn}

Let $A=\{ a\, ,\, b \}$ and let $k_n=2^{2^n}$. In the alphabet $A$
consider the sequence of words $w_n=\left( a^{k_n} b^{k_n}
a^{-k_n} b^{-1} \right)^{k_n}$. Note that $|w_n|=k_n (3k_n+1)$. In
what follows we denote this length by $d_n$ and the sequence
$(d_n)$ by $d$.

A standard argument gives the following result (see
\cite{ThomasVelickovic} and \cite{Bowditch:qiclasses} for versions
of it).

\begin{lem}
If $\mathcal{W}$ is the minimal collection of reduced words in
$\free_A$ containing $\{ w_n \; ;\; n\in \N \, ,\, n\geq 4 \}$,
closed with respect to cyclic permutations and taking inverses,
then the following hold:
\begin{itemize}
\item[(1)] $\mathcal{W}$ can be generated recursively;

\item[(2)] $\mathcal{W}$ satisfies $C^*(1/500)$;

\item[(3)] for every $n\in \N$, the set $\{ w\in \mathcal{W} \;
;\; |w|\geq d_n \}$ satisfies $C^*(1/k_n )$.
\end{itemize}
\end{lem}

\me

\begin{prop}(\cite{ThomasVelickovic},
\cite{DrutuSapir:TreeGraded})\label{psmallc}
The two-generated and recursively presented group $G=\la a,b\mid
w_n \, ,\, n\geq 4 \ra,$ has the following properties.
\begin{itemize}
    \item[(1)] Any asymptotic cone of $G$ is either a real tree or a tree-graded
space with pieces isometric to the same circle with the arc
distance.
    \item[(2)] The group $G$ is not relatively hyperbolic.
\end{itemize}
\end{prop}

\begin{pf} \textbf{(1)}\quad Let $\Re_n$ be the loop through $1$
in the Cayley graph of $G$, labeled by the word $w_n$ starting
from $1$.

In \cite[$\S 7$]{DrutuSapir:TreeGraded} it is proved that the
asymptotic cone $\mathrm{Cone}_\omega (G; 1, d )$ is tree-graded,
with the set of pieces composed of ultralimits of sequences of the
form $(g_n\Re_n)$ where $g_n \in G$. In our case these ultralimits
are all isometric to the unit circle. The same proof works in fact
not only for $(d_n)$ but for any scaling sequence, thus giving the
statement in (1), since for other scaling sequences the
ultralimits can be either circles, points or lines. A version of
the last part of the argument can also be found in
\cite{ThomasVelickovic}.

\medskip

\textbf{(2)}\quad Assume that the group $G$ is hyperbolic relative
to a finite family of finitely generated subgroups $\H$. Then
$\mathrm{Cone}_\omega (G; 1, d )$ is tree-graded with set of
pieces ultralimits of left cosets of subgroups in $\H$. According
to Lemma 2.15 in \cite{DrutuSapir:TreeGraded}, the subset without
cut-point $\lim_\omega (\Re_n)$ is contained in some $\lim_\omega
g_n H$ where $H\in \H$.

Let $\fp_n$ be an arbitrary sub-path in $\Re_n$, of length
$\frac{1}{6} d_n$. This sub-path is a geodesic in the Cayley graph
of $G$ \cite[$\S 7.2 $]{DrutuSapir:TreeGraded}. Let $\fp_n'$ and
$\fp_n''$ be the first and the last third of $\fp_n$. Since both
have length $\frac{1}{18} d_n$ and are contained $\omega$-almost
surely in $\nbhd_{o(d_n)} (g_n H)$, property $(\alpha_2)$ implies
that both intersect a tubular neighborhood of radius $O(1)$ of
$g_n H$. The quasi-convexity of $g_n H$ (\cite[$\S
4$]{DrutuSapir:TreeGraded}, \cite[$\S 4.3 $]{Drutu:RelHyp})
implies that $\omega$-almost surely the middle third of $\fp_n$ is
contained in $\nbhd_{M} (g_n H)$, for some uniform constant $M$.
Now the loop $\Re_n$ can be divided into $18$ sub-paths of length
$\frac{1}{18} d_n$, each of which appears as the middle third of a
larger sub-path. We may conclude that $\Re_n$ is $\omega$-almost
surely contained in $\nbhd_{M} (g_n H)$. In particular $1\in
\nbhd_{M} (g_n H)$, hence it may be assumed that $g_n\in B(1,M)$.
Since $B(1,M)$ is finite, the ultrafilter allows us to assume that
$g_{n}$ is a constant sequence.

Thus we obtained that for some $g\in B(1,M)$ and some $H\in \H$
the left coset $gH$ contains in its $M$-tubular neighborhood
$\omega$-a.s. the loop $\Re_n$. It follows that $a\Re_n \subset
\nbhd_M(agH)$ and $b\Re_n \subset \nbhd_M(bgH)$ $\omega$-a.s. The
loop $a\Re_n$ has in common with $\Re_n$ the path $a\fp_a$, where
$\fp_a $ is the path of origin $1$ and label $a^{k_n -1}$. It
follows that $\omega$-a.s $\nbhd_M(gH)$ and $\nbhd_M(agH)$
intersect in a set of diameter at least $k_n -1$. Property
$(\alpha_1)$ implies that $gH=agH$, thus $a\in gHg^{-1}$.

Likewise, the remark that $b\Re_n$ and $\Re_n$ have in common the
path $b\fp_a$, together with $(\alpha_1)$, implies that $b\in
gHg^{-1}$. It follows that $G$ coincides with $gHg^{-1}$, hence
with $H$, therefore the relative hyperbolic structure defined by
$\H$ is not proper.\end{pf}

\begin{rmk}\label{smallc-nrh}
The arguments in the proof of
statement (2), Proposition \ref{psmallc}, can be carried out for a
much more general construction of the group $G$ than the one
considered here. Thus, the techniques described in \cite[$\S
7$]{DrutuSapir:TreeGraded} (following an idea from
\cite{Olshanskii:initial} further developed in
\cite{ErschlerOsin}) allow the construction of a large class of new
examples of NRH groups.
\end{rmk}

\begin{cor}
The group $G$ does not contain any subspace $B$ which endowed with
the restriction of a word metric on $G$ becomes unconstricted.

In particular $G$ is not metrically thick.
\end{cor}

\begin{pf}
Assume that $G$ contains an unconstricted subspace $B$. Then there
exists an ultrafilter $\omega$ and a sequence $\delta$ of positive
numbers such that for every sequence of observation points $b$ in
$B$ the asymptotic cone $\mathrm{Cone}_\omega (B; b, \delta )$
does not have cut-points. Since $B$ is  endowed with the
restriction of a word metric on $G$, $\mathrm{Cone}_\omega (B; b,
\delta )$ can be seen as a subset of $\mathrm{Cone}_\omega (G; b,
\delta )$.

If $\mathrm{Cone}_\omega (G; b, \delta )$ is a real tree then all
arc-connected subsets in it have cut-points, thus it cannot
contain a subset $\mathrm{Cone}_\omega (B; b, \delta )$ as above.

Assume that $\mathrm{Cone}_\omega (G; b, \delta )$ is a
tree-graded space with pieces isometric to a circle. Lemma 2.15 in
\cite{DrutuSapir:TreeGraded} implies that $\mathrm{Cone}_\omega
(B; b, \delta )$ is contained in some piece. This is impossible
since $\mathrm{Cone}_\omega (B; b, \delta )$ is infinite diameter,
by Definition \ref{ums}, (2).
\end{pf}

\begin{rmk}
Note that the group $G$ displays a sort of generalized version of
metric thickness with respect to the collection of subspaces
$\{g\Re_n \; ;\; g\in G\, ,\, n\geq 4\}$. Indeed this collection
satisfies one of the two necessary conditions for uniform
unconstrictedness (condition (1) in Definition \ref{duunc}),
property $(\bn_1)$ of a metric network obviously holds, and a
weaker version of property $(\bn_2)$ is satisfied: the diameters
of the intersections between neighborhoods of consecutive
subspaces $L_i, L_{i+1}$ in a sequence connecting thickly are no
longer infinite, but increase with the minimum between the
diameters of the starting and the target subspaces $L$ and $L'$.
\end{rmk}

\begin{qn}
Can the construction above be adapted to give an example of a
group which is metrically thick (and thus NRH) but not
algebraically thick?
\end{qn}

\section{Mapping class groups}\label{SectionMCG}

Let $S=S_{g,p}$ denote an orientable surface of genus $g$ with $p$
punctures. We parameterize the complexity of $S$ by
$\xi(S)=3g+p-3$ which is the cardinality of any set of closed
curves subdividing $S$ into pairs of pants, that is, any maximal,
pairwise disjoint, pairwise nonhomotopic set of essential,
nonperipheral closed curves on $S$. Note that every surface with
$\xi(S)\leq 1$ either has $\MCG(S)$ finite or virtually free; in
particular, these groups are all $\delta$--hyperbolic. This
section provides our first example of an algebraically thick
group:

\begin{thm}\label{mcgthick} $\MCG(S)$ is algebraically thick of order
one when $\xi(S)\geq 2$.
\end{thm}
It is known that the mapping class group is not thick of order 0
(i.e., unconstricted) by the following:

\begin{thm}[Behrstock \cite{Behrstock:asymptotic}]
    \label{mcgcutpoints} For every surface
    $S$, every asymptotic cone of $\MCG(S)$ has cut-points.
\end{thm}
$\MCG(S)$ is not hyperbolic when $\xi(S)\geq 2$ since for any set
of curves subdividing $S$ into pairs of pants, the subgroup
generated by Dehn twisting along these curves is a free abelian
subgroup of $\MCG(S)$ of rank $\xi(S)$. Indeed, according to
\cite{BirmanLubotzkyMcCarthy}, $\xi(S)$ is the maximal rank of a
free abelian subgroup of $\MCGS$. Moreover, it has been shown that
these abelian subgroups are quasi-isometrically embedded in
$\MCGS$ (see \cite{FarbLubotzkyMinsky} and
\cite{Mosher:automatic}). Masur and Minsky showed that $\MCGS$ is
weakly relatively hyperbolic with respect to a finite collection
of stabilizers of curves \cite{MasurMinsky:complex1}. (The
subgroup stabilizing a curve $\gamma$ will be denoted
${\mathrm{stab}} (\gamma)$.) Further, it is easily verified that
$\MCGS$ is not relatively hyperbolic with respect to such a
collection of subgroups. This motivates the question of whether
there exists a collection of subgroups of $\MCGS$ for which this
group is relatively hyperbolic (see \cite{Behrstock:asymptotic}).
That no such collection exists is an immediate consequence of
Theorem~\ref{mcgthick}:

\begin{cor}\label{mcgnobcp} If $S$ is any surface with
$\xi(S)\geq 2$, then there is no finite collection of finitely
generated
proper subgroups with respect to which $\MCG(S)$ is relatively
hyperbolic.
\end{cor}

Anderson, Aramayona, and Shackleton have an alternative proof of
Corollary~\ref{mcgnobcp} using an algebraic characterization of
relative hyperbolicity due to Osin \cite{AASh:RelHypMCG}. This
result also appears in both \cite{Bowditch:3manifolds} and
\cite{KarlssonNoskov} although it is not stated as such as it
appears under the guise of a fixed point theorem for actions of
the mapping class group.
We note that the techniques of each of
\cite{AASh:RelHypMCG}, \cite{Bowditch:3manifolds}, and
\cite{KarlssonNoskov} rely in an essential way on the group
structure.

\medskip

Before giving the proof of Theorem~\ref{mcgthick} we recall some
well known results concerning mapping class groups. For closed
surfaces the mapping class group was first shown to be finitely
generated by Dehn \cite{Dehn:f.g.MCG} in a result which was later
independently rediscovered by Lickorish \cite{Lickorish:f.g.MCG};
both gave generating sets consisting of finite collections of Dehn
twists. For the mapping class group $\MCGS$ of a punctured surface
$S$, the finite index subgroup which fixes the punctures pointwise
is generated by a finite set of Dehn twists \cite{Birman:Braids};
this latter group is
also called the \textit{pure mapping class group}, and denoted by
$\pmcg(S)$. The extended mapping class group, $\mcg^{\pm}(S)$, is
the group of orientation preserving and reversing mapping classes.
This is a finite extension of the mapping class group. (See
\cite{Birman:Braids}, \cite{Ivanov:mcg},
\cite{Humphries:Generators}). Since these groups are all
quasi-isometric, Remark~\ref{thickqiinvariant} implies that if we
can show that the pure mapping class group is algebraically thick
of order one, it implies that the same holds for the mapping class
group and the extended mapping class group.

Introduced by Harvey, a useful tool in the study of $\MCGS$ is the
\emph{complex of curves} $\C(S)$ \cite{Harvey:Boundary}. When
$\xi(S) \ge 2$ the complex $\C(S)$ is a simplicial complex with
one vertex corresponding to each homotopy class of nontrivial,
nonperipheral simple closed curves in $S$, and with an $n$-simplex
spanning each collection of $n+1$ vertices whose corresponding
curves can be realized on $S$ disjointly.

For later purposes we also need to define $\C(S)$ when $\xi(S)=1$,
in which case the surface $S$ is either a once-punctured torus or
a four-punctured sphere: the vertex set of $\C(S)$ is defined as
above, with an edge attached to each pair of vertices whose
corresponding curves can be realized on $S$ with minimal
intersection number, that number being~$1$ on a once-punctured
torus and~$2$ on a four-punctured sphere.

In either case the complex $\C(S)$ is connected (see for example
\cite{MasurMinsky:complex1}). The distance
$d_{\C(S)}(\alpha,\beta)$ between two vertices $\alpha,\beta$ in
$\C(S)$ is the usual simplicial metric, defined to be the length
of the shortest edge path between $\alpha$ and $\beta$.

\begin{proof}[Proof of Theorem~\ref{mcgthick}] We start by remarking
that for any essential simple closed curve $\gamma\in\C(S)$, its stabilizer
${\mathrm{stab}} (\gamma)$ in $\MCG(S)$ is a central extension of $\mcg
(S\setminus \gamma)$ by the infinite cyclic subgroup generated by
a Dehn twist about $\gamma$. Thus, if $\xi (S) \geq 2$ then
${\mathrm{stab}} (\gamma )$ is non-elementary and has a central
infinite cyclic subgroup. Consequently ${\mathrm{stab}} (\gamma)$
is unconstricted. It is an easy consequence of the distance
estimates in \cite[Theorem~6.12]{MasurMinsky:complex2}, that
${\mathrm{stab}}(\gamma )$ is undistorted for any essential simple
closed curve $\gamma$. Select a finite collection of curves
$\Gamma_0$ such that the Dehn twists along these curves generate
$\mathcal{P}\mcg (S)$. Connectivity of the curve complex implies
that there is a finite connected subgraph of $\C(S)$ containing
the vertices in $\Gamma_0$; let $\Gamma$ denote the set of
vertices in this new graph. Since $\xi (S) \geq 2$, if
$\alpha,\beta$ are curves representing vertices at distance~1 in
$\C(S)$ then $\alpha $ and $\beta$ are disjoint, and so the
subgroup ${\mathrm{stab}} (\alpha )\cap {\mathrm{stab}}
(\beta)={\mathrm{stab}} (\alpha \cup \beta )$ is infinite. It
follows that $\mcg(S)$ is algebraically thick of order at most~$1$
with respect to $\H=\left\{ {\mathrm{stab}} (\gamma)\mid \gamma
\in \Gamma \right\}$. By Theorem~\ref{mcgcutpoints}, $\MCG(S)$ is
not unconstricted and thus it is thick of order $1$.
\end{proof}

\section{$\autfn$ and $\outfn$} \label{section:outfn}

We start by fixing a set of generators
$\{x_1 ,..., x_n\}$ for the free group $\free_{n}$. We denote the
automorphism and outer automorphism groups of $\free_{n}$ by
$\autfn$ and $\outfn = \autfn / \text{Inn}(\free_n)$,
respectively, where $\text{Inn}(\free_n)$ is the group of inner
automorphisms.
Recall that an element of $\autfn$ is a \textit{special}
automorphism if the induced automorphism of $\Z^n$ has determinant
$1$. The subgroup $\sautfn$ of special automorphisms has index two
in $\autfn$.

\medskip

\noindent{\textit{Notation:}} All indices in this section are
taken modulo $n$, where $n$ is the rank of the free group we are
considering.

\medskip

 We denote the following Dehn twists in $\autfn$:
\begin{itemize}
\item $r_{i} = \begin{cases}
x_{i+1} \mapsto x_{i+1} x_{i} & \\
x_j \mapsto x_j & \text{for} \, j\neq i+1\, ,
\end{cases}
$
\item $l_i = \begin{cases}
x_{i+1}\mapsto x_ix_{i+1} & \\
x_j\mapsto x_j & \text{for} \, j\neq i+1\, ,
\end{cases}
$
\item $n_{i} = \begin{cases}
x_{i+2}\mapsto x_{i+2}\, x_{i} & \\
x_j\mapsto x_j & \text{for} \, j\neq i+2\, .
\end{cases}
$
\end{itemize}
Culler and Vogtmann proved that the set $S$ composed of all $r_i$
and $l_i$ is a set of generators of $\sautfn$, see
\cite{CullerVogtmann:FA}. Note that all elements in $S$ have
infinite order.
The elementary argument in Example 2.4 of
\cite{Alibegovic:translation} yields the following.

\begin{lem}\label{undistortedAutflats} Let $n\geq 3$.
The $\Z^{2}$ subgroup of $\autfn$
generated by the pair $\<\phi_{i},\phi_{j}\>$
is undistorted when $\phi_i\in \{ r_i,l_i\}$, $\phi_j\in
\{r_j,l_j\}$, and $\dist (i,j)\geq 2$, where $\dist (i,j)$ is
measured in $\Z / n\Z$. The $\Z^{2}$ subgroups
$\<r_{i},l_{i}\>$, $\<n_{i},r_{i}\>$, and $\<n_{i},l_{i+1}\>$
are also undistorted for all $i$. These subgroups also
inject to undistorted subgroups of $\outfn$.
\end{lem}

\begin{thm}\label{autthick}
If $n\geq 3$, then both $\autfn$ and $\outfn$ are algebraically
thick of order at most one.
\end{thm}

\proof We consider $\H$ the set of all subgroups $\langle \,
\phi_i\, ,\, \phi_j\, \rangle$, where $\phi_i\in \{ r_i,l_i\}$,
$\phi_j\in \{ r_j,l_j\}$, and $\dist (i,j)\geq 2$, and we also
include in $\H$ the subgroups $\<n_{i},r_{i}\>$ and
$\<n_{i},l_{i+1}\>$. We may regard these as subgroups of $\autfn$
or, since they each intersect $\text{Inn}(\free_n)$ trivially, as
subgroups of $\outfn$.

We shall now prove that both $\autfn$ and $\outfn$ are
algebraically thick of order one with respect to the subgroups in
$\H$.

Each subgroup $H=\langle \, \phi\, ,\, \psi\, \rangle$ in $\H$ is
isomorphic to $\Z^2$, hence unconstricted.
Lemma~\ref{undistortedAutflats} shows that each such subgroup is
undistorted.

By \cite{CullerVogtmann:FA}, the $r_{i}$ and $l_{i}$ provide a
complete set of generators for $\sautfn$, and $\sautfn$ is a
subgroup of $\autfn$ of index two, thus we have shown that
property ($\ban_1$) is satisfied for $\sautfn$ and thus for
$\autfn$.

We verify property ($\ban_2$) in the definition of algebraic
thickness. Note that since $\< \phi,l_j\> \cap\<\phi,r_j\>\supset
\<\phi\>$, it suffices to show that the subgroups generated by
$r_i$ and $n_i$ can be thickly connected. For every $\<r_i, r_j\>$
with $\dist (i,j)\geq 2$, Lemma~\ref{undistortedAutflats} shows
that the subgroup $\<r_i,r_j\>$ thickly connects any pair of
subgroups of $\H$ where one contains $r_i$ and the other $r_j$.
Thus, to finish the verification of property ($\ban_2$) it remains
to find sequences joining a pair of subgroups, where one contains
$r_i$ and the other $r_{i+1}$. Observe that the sequence of
subgroups $\<r_{i},n_{i}\>$, $\<n_{i},l_{i+1}\>$,
$\<l_{i+1},r_{i+1}\>$ each intersects the next in an infinite
diameter subset. This shows that any subgroup containing $r_{i}$
can be thickly connected to one containing $r_{i+1}$ through a
sequence of subgroups in $\H$, thereby completing our verification
of property ($\ban_2$).

All the subgroups of $\autfn$ that are used above to prove
thickness are mapped, via the canonical epimorphism, injectively
and without distortion to $\outfn$. Thus the hypotheses of
Definition~\ref{dgthick} hold as well in $\outfn$, whence $\outfn$
is algebraically thick of order one for $n\geq 3$.
\endproof

\section{Artin groups}\label{section:artin}

An \emph{Artin group} is a group given by a presentation of the
following form:
\begin{align}\label{apres}
A &= \langle x_1,..., x_n\mid (x_i,x_j)_{m_{ij}} = (x_j,x_i)_{m_{ji}}
\rangle
\end{align}
where, for all $i \ne j $ in $\{1,\ldots,n\}$,
$$m_{ij} = m_{ji} \in \{2,3,\ldots, \infinity\} \qquad\text{and}\qquad
(x_i,x_j)_{m_{ij}} =
\begin{cases}
\text{Id} &\text{if}\quad m_{ij} = \infinity\, , \\
\underbrace{x_ix_jx_i...}_{m_{ij}\mbox{ terms }} &\text{if}\quad
m_{ij} < \infinity\, .
\end{cases}
$$
Such a group can be described by a finite (possibly disconnected)
graph $\mathcal{G}_A$, the \emph{Artin presentation graph}, where
the vertices of $\ga$ are labeled $1,\ldots,n$ in correspondence
with the generators $x_1,\ldots,x_n$, and the vertices $i$ and $j$
are joined by an edge labeled by the integer $m_{ij}$ whenever
$m_{ij} < \infinity$. When $m_{ij}=\infinity$ there is no
associated relator in the presentation (\ref{apres}), and $\G_A$
has no edge between vertices $i$ and~$j$.

A subgroup generated by a subset $S$ of $\{x_1,...,x_n\}$ is
called a \emph{special subgroup} of $A$ and it is denoted by
$A_S$. Any special subgroup $A_S$ is itself an Artin group with
presentation given by the relations in (\ref{apres}) containing
only generators in $S$, and such that $\G_{A_S}$ is the subgraph
of $\G_A$ spanned by the vertices corresponding to $S$. This has
been proved by Van der Lek in \cite[Chapter II, Theorem
4.13]{Lek}. See also \cite{Paris:Parabolic} for an elementary
proof as well as for a history of the result.

In particular the two generator special subgroup $A_{ij}$
generated by $x_i,x_j$ is an Artin group: if $m_{ij}=\infinity$
then $A_{ij}$ is free of rank~2; whereas if $m_{ij}<\infinity$
then $A_{ij}$ is defined by the single relator $(x_i,x_j)_{m_{ij}}
= (x_j,x_i)_{m_{ji}}$.

The \emph{Coxeter group $W$} associated to an Artin group $A$ has
a presentation obtained from (\ref{apres}) by adding relations
saying that each $x_i^2$ is the identity.

\medskip

\begin{exa}\label{twogeneratorartin} A two generator Artin group
$\<x,y \suchthat
(x,y)_m = (y,x)_m\>$ with $m<\infinity$ is unconstricted. This holds
since the element $(x,y)_{2m}$ is central, and it is of infinite
order since it projects to a nonzero element of $\Z$ under the
exponent sum homomorphism $A \to \Z$.
\end{exa}

In \cite{KapovichSchupp:Artin} the following has been proven.
\begin{thm}[I. Kapovich--P. Schupp]
\label{TheoremKapovichSchupp}
An Artin group $A$ defined as in (\ref{apres}) with $m_{ij}\geq 7$
for all $i\neq j$ is weakly hyperbolic relatively to the collection
of two generator special subgroups
$$\H =\{A_{ij} \suchthat m_{ij}<\infinity \}\, .$$
\end{thm}
As noted in the same paper, the above result cannot be improved to
say that $A$ is strongly hyperbolic relative to $\H$. Nevertheless
the question remained whether $A$ was strongly hyperbolic relative
to other groups, or at least metrically hyperbolic relative to
some collection of subsets. Our methods give a partial answer to
this question, with the interesting outcome that when our methods
work, $A$ turns out to be algebraically thick of order at most~1 with
respect to the exact same collection $\H$. We do not go as far as
to check thickness for all of the Artin groups in
Theorem~\ref{TheoremKapovichSchupp}, but in
Corollary~\ref{CorollaryCATZeroArtin} below we show thickness as
long as the graph $\G_A$ has no triangles.
Here are some other special classes of Artin groups.

\paragraph{Free decompositions.} The graph $\G_A$ with $n$ points and
no edges describes the group with $n$ generators and no relators,
i.e., the free group on $n$ generators. More generally, if $\G_A$
is disconnected then $A$ decomposes into a free product, one
factor for each connected component in the defining graph.
The converse is true as well: if $\G_A$ is connected then $A$ is
freely indecomposable, in fact $A$ is a one-ended group. This
follows for example from Proposition \ref{allArtin} and Remark
\ref{oneend}.
Since any nontrivial free product is relatively hyperbolic with
the free factors as peripheral subgroups, we henceforth restrict
our attention to one-ended Artin groups, those whose defining
graphs have only one connected component.

\paragraph{Right angled Artin groups and even Artin groups.} The
complete graph on $n$ vertices with each $m_{ij}=2$ describes the
group with $n$ commuting
generators, i.e., $\Z^{n}$. More generally, a \emph{right angled
Artin group} is one for which $m_{ij}\in\{2,\infty\}$ for all $i,j$.
Recently there has been interest in
the quasi-isometric classification of right angled Artin groups
(see \cite{BehrstockNeumann:qigraph} and \cite{BKS:raagi}).
Generalizing a right angled Artin group, an Artin group is
\emph{even} if each $m_{ij}$ is an even integer or infinity.

\paragraph{Finite type Artin groups.} An Artin group is of
\emph{finite type} if the associated Coxeter group $W$ is finite.
For example, the \emph{braid group} on $n$ strands is the Artin
group with $n$ generators, with $m_{i,i+1}=3$, and $m_{ij}=2$ if
$|i-j|>1$ --- in this case the associated Coxeter group $W$ is
just the symmetric group on $n$ symbols. An Artin group of finite
type is unconstricted, since it has an infinite cyclic central
subgroup of infinite index, as proven in \cite{BrieskornSaito} and
\cite{Deligne:Immeubles}.

\paragraph{Affine type Artin groups.} An Artin group $A$ is of
\emph{affine type} if the associated Coxeter group $W$ is a
Euclidean crystallographic group. For example, when $\G_A$ is a
cycle of $n+1$ edges with a~3 on each edge then $W$ is the full
group of symmetries of a tiling of $\reals^n$ by cubes, in which
case we denote $A = \wt A_n$.

\me

The reason for so many different special classes of Artin groups
seems to be a proliferation of techniques for studying various
aspects of Artin groups, and a concomitant lack of any single
technique that works on all Artin groups --- most theorems about
Artin groups carry extra hypotheses on the Artin presentation. For
example, there are various constructions in the literature of
biautomatic and/or $\CAT(0)$ structures on Artin groups (we refer
the reader to \cite{ECHLPT} for the definition of a biautomatic
structure):
\begin{itemize}
\item Every right angled Artin group is $\CAT(0)$, in fact it is the
fundamental group of a nonpositively curved cube complex
\cite{BestvinaBrady:MorseTheory}, and so it is biautomatic
\cite{NibloReeves:cubed}.
\item Braid groups are biautomatic \cite{ECHLPT}. More generally,
Artin groups of finite type are biautomatic \cite{Charney:FTABiaut}.
\item If $\G_A$ has no triangles then $A$ is $\CAT(0)$
\cite{BradyMcCammond:Artin} and biautomatic (combining
\cite{Pride:Artin} and \cite{GerstenShort:automatic}; see comments in
\cite{BradyMcCammond:Artin}).
\item $A$ is $\CAT(0)$ and biautomatic if the edges of $\G_A$ can be
oriented so that each triangle has an orientation agreeing with the
orientations of all three edges, and in each square the orientations
of the four edges do not alternate when going around the square
\cite{BradyMcCammond:Artin}.
\item Artin groups for which each $m_{ij} \ge 4$ are biautomatic
\cite{Peifer:ExtraLarge}.
\item Artin groups of affine type $\wt A_n$, also known as the
\emph{affine braid groups}, are biautomatic
\cite{CharneyPeifer:affine}.
\end{itemize}

We shall prove thickness for some of these groups. The method we use
is:

\begin{lem}\label{Lemma2Gundistorted} \quad
If the graph $\G_A$ is connected, and if each two generator
special subgroup $A_{ij}$ with $m_{ij} < \infinity$ is undistorted
in $A$, then $A$ is algebraically thick of order~$\le 1$.
\end{lem}

\begin{proof} For $i,j,k$ all distinct, the subgroup $A_{ij}
\intersect A_{ik}$ contains the infinite order element $x_i$.
Since $\G_A$ is connected, and since the two generator special
subgroups $A_{ij}$ with $m_{ij} < \infinity$ are undistorted and
unconstricted (see Example \ref{twogeneratorartin}), the lemma
follows.
\end{proof}
One can verify undistortedness of two generator special subgroups in
different cases by using a variety of methods: retractions;
nonpositive curvature methods; the Masur--Minsky distance estimates
for mapping class groups; or automatic group methods.

\paragraph{Retractions.}
Our first results on Artin groups use a simple algebraic method to
prove undistortedness:

\begin{prop}
\label{PropRetraction}
Let $A$ be an Artin group. Suppose that for each 2 generator special
subgroup $A_{ij}$ with $m_{ij} < \infinity$, there exists a
retraction $p \from A \to A_{ij}$. Then $A$ is algebraically thick of
order $\le 1$.
\end{prop}

\begin{proof} This is a consequence of Lemma~\ref{Lemma2Gundistorted}
and the observation that for any finitely generated group $G$ and any
finitely generated subgroup $H \subgroup G$, if there exists a
retraction $G \to H$ then $H$ is undistorted.
\end{proof}
In each application of this proposition, the retraction from an Artin
group $A$ generated by $S$ to a special subgroup $A'$ generated by
$S' \subset S$ will be induced by a retraction from $S \union
\{\Id\}$ to $S' \union \{\Id\}$.

\paragraph{Even Artin groups.}
Consider first the case of an Artin group $A$ presented
by~(\ref{apres}) so that each $m_{ij}$ is an even integer or
$+\infinity$.
For each generator $g\in S$ define $p(g)=g$ if $g\in S'$ and $p(g)=1$
otherwise.
This projection is well defined, since any relation
$$[x_i, x_j]_{m_{ij}} = [x_j, x_i]_{m_{ij}}
$$
projects under $p$ to either: itself if both $x_i,x_j\in S'$, or to the tautological relation $x_i^{m_{ij}} = x_i^{m_{ij}}$ if $x_j\not\in S'$.

By Proposition~\ref{PropRetraction} it follows that:

\begin{thm}\label{ThmEvenArtin} Even Artin groups are algebraically
    thick of order at most~1.
\end{thm}

\paragraph{Trees.}
Consider next the case that $\G_{A}$ is a tree. There is a unique
retraction $p \from \G_A \mapsto \G_{A'}$ so that each component of
$\G_A - \G_{A'}$ maps to the unique vertex of $\G_{A'}$ incident to
that component. This induces a map $p \from S \mapsto S'$. Extend $p$
to a map from words in $S$ to words in $S'$. Again we need only prove
that given a relator $v R_{ij} v^{-1}$ for $A$ as above, $p(v R_{ij}
v^{-1}) = w p(R_{ij}) w\inv$ defines the identity in $A'$. Consider
the edge $e$ of $\G_{A'}$ connecting $s_i$ to $s_j$. If $e \subset
A'$ then $p(R_{ij}) = R_{ij}$ and we are done. If $e$ is contained in
a component of $\G_A - \G_{A'}$ incident to a vertex $s_k$ of
$\G_{A'}$ then $p(R_{ij})$ is a word in the single generator $s_k$
with exponent sum equal to zero and so is freely equal to the
identity. By Proposition~\ref{PropRetraction}, $A$ is algebraically
thick of order~$\le 1$.

\paragraph{Other examples.}
There seem still to be numerous other examples to which
Proposition~\ref{PropRetraction} applies. For example, consider
the case that the group $\G_A$ has rank~1, meaning that it
deformation retracts onto a circular subgroup $\G_{A'}$. Suppose
furthermore that each integer that occurs as a label $m_{ij}$ on
some edge of $\G_{A'}$ occurs for at least two different edges.

For any edge of $\G_A$ not in $\G_{A'}$ there is a retraction
defined as in the example above where the graph is a tree. For any
edge $A_{ij}=e \subset \G_{A'}$, let $f \subset \G_{A'}$ be
another edge with the same integer label. Removing the interiors
of $e$ and $f$ from $\G_A$ results in two connected subgraphs
$\G_i$, $\G_j$, with notation chosen so that $x_i \in \G_i$. Let
$y_i,y_j$ denote the endpoints of $f$, with notation chosen so
that $y_i \in \G_i$. There is a retract $A \to e$ defined by
taking $\G_i$ to $x_i$, and taking $f$ to $e$ so that $y_i$ goes
to $x_i$. This retract restricts to a retraction of the generating
set $S$ onto $\{x_i,x_j\}$. This map extends to a well defined
retraction $A \to A_{ij}$ for the following reasons: for edges not
equal to $f$ the corresponding Artin relation maps to a word
freely equal to the identity; and the Artin relation for the edge
$f$ maps to the Artin relation for the edge $e$ because those two
edges are labeled by the same integer.

\bigskip

We have not investigated the full extent to which
Proposition~\ref{PropRetraction} applies, but on the other hand we
can easily construct somewhat random examples to which it seems
not to apply, for example an Artin group whose presentation graph
is the complete graph on four vertices and whose six edges are
labeled by six pairwise relatively prime integers.

\paragraph{Nonpositive curvature.} A good reference for nonpositively
curved groups is \cite{BridsonHaefliger}. A geodesic metric on a
cell complex $C$ is a \emph{polyhedral Euclidean metric} if for
each cell $c$ there is a compact, convex Euclidean polyhedron $P$
and a characteristic map $P \mapsto c$ so that the metric on $P$
pushes forward to the given metric on $c$. A polyhedral spherical
metric is similarly defined, using spheres of constant curvature
$+1$ instead of Euclidean space. The link of each vertex in a
polyhedral Euclidean metric inherits a polyhedral spherical
metric.

If $C$ comes equipped with a polyhedral Euclidean metric then we
say that $C$ is a \emph{piecewise Euclidean cell complex}.
Furthermore, if the link of each vertex $v \in C$ has no closed
geodesic of length $< 2\pi$ then we say that $C$ is
\emph{nonpositively curved}. A subcomplex $D \subset C$ is
\emph{locally convex} if for each vertex $v \in D$, the link of
$v$ in $D$ is a geodesically convex subset of the link of $v$ in
$C$.

\begin{prop}[\cite{BridsonHaefliger}]
\label{PropNPC}
If $C$ is a finite piecewise Euclidean non-positively curved cell
complex, and if $D$ is a locally convex subcomplex, then the
inclusion of universal covers $\wt D \to \wt C$ is globally
isometric. It follows that the inclusion $D \inject C$ induces an
injection $\pi_1(D) \to \pi_1(C)$ with undistorted image.
\qed\end{prop}
Although right angled Artin groups are already considered
in Theorem~\ref{ThmEvenArtin}, the following gives a different
approach.

\begin{thm}\label{TheoremCATZeroArtin}
If the Artin group $A$ is right angled, or if it satisfies Pride's
condition that $\G_A$ has no triangles, then $A$ is the
fundamental group of a piecewise Euclidean non-positively curved
cell complex $C_A$ so that each 2 generator special subgroup
$A_{ij}$ is the inclusion induced image of a locally convex
subcomplex of $C_A$.
\end{thm}
The proof is given below. Combining Theorem~\ref{TheoremCATZeroArtin}
with
Lemma~\ref{Lemma2Gundistorted} and Proposition~\ref{PropNPC} we
obtain:

\begin{cor}\label{CorollaryCATZeroArtin}
Artin groups $A$ which are right angled or for which $\G_A$ has no
triangles are algebraically thick of order~$\le 1$.
\qed\end{cor}
In one case we can compute the order to be exactly~$1$:

\begin{cor}\label{raagcutpoints}
Any right angled Artin group $A$ for which $\G_A$ is a tree of
diameter at least~$3$ has cut-points in every asymptotic cone, and so
$A$ is thick of order~$1$.
\end{cor}

\begin{proof} Once we construct a compact, non-Seifert fibered,
3-dimensional graph manifold $M$ whose fundamental group is
isomorphic to $A$, the result follows by work of
\cite{KapovichLeeb:3manifolds} and \cite{KKL:QI} (see
Lemma~\ref{suplin} and Section~\ref{Sectiongrm}). The manifold $M$
will be a ``flip manifold'' in the terminology of
Section~\ref{Sectiongrm}.

Consider first a right angled Artin group $A'$ for which $\G_{A'}$
is a \emph{star graph}, meaning a tree of diameter~$2$, with
valence~$1$ vertices $v_1,\ldots,v_k$ for $k \ge 2$, and a
valence~$k$ vertex $v_0$ called the \emph{star vertex}. The group
$A'$ is the product of a rank $k$ free group with $\Z$. We can
realize $A'$ as the fundamental group of a 3-manifold $M'$ which
is the product of a ``horizontal'' $k+1$-holed sphere crossed with
a ``vertical'' circle, so that the generators $v_1,\ldots,v_k$
correspond to the horizontal circles in $k$ of the boundary tori,
and the generator $v_0$ corresponds to the vertical circle.

Suppose now that $A$ is a right angled Artin group and $\G_A$ is a
tree of diameter~$\ge 3$. Let $v_1,\ldots,v_m$ be the vertices of
$\G_A$ of valence~$\ge 2$, and note that $m \ge 2$. Let $\G_{A_i}$
denote the maximal star subgraph of $\G_A$ with star vertex $v_i$.
The graph $\G_{A_i}$ presents a special subgroup $A_i$ which is
isomorphic to the fundamental group of a 3-manifold $M_i$ as
above, homeomorphic to the product of a sphere with holes crossed
with the circle. We have $\G_A = \G_{A_1} \union \cdots \union
\G_{A_m}$. When $i \ne j$ and $\G_{A_i}$, $\G_{A_j}$ are not
disjoint then $\G_{A_i} \intersect \G_{A_j}$ is a single edge of
$\G_A$, in which case $M_i$ and $M_j$ each have a torus boundary
whose fundamental group corresponds to the $\Z^2$ special subgroup
generated by $v_i$ and $v_j$; we now glue these two tori so that
the horizontal circle on one torus glues to the vertical circle on
the other. The result of gluing $M_1,\ldots,M_m$ in this manner is
the desired 3-manifold $M$, and $M$ is not Seifert fibered because
$m \ge 2$.
\end{proof}

\begin{proof}[Proof of Theorem \ref{TheoremCATZeroArtin}]
Suppose first that $A$ is right angled. For each subset $I$ of the
set of generator indices $\{1,\ldots,n\}$ for which the generators
$\{x_i \suchthat i \in I\}$ all commute with each other, let $T_I$ be
the Cartesian product of $\abs{I}$ copies of the unit circle. Glue
these tori together using the obvious injection $T_{I'} \inject T_I$
whenever $I' \subset I$, with base point $T_\emptyset$. The result is
a nonpositively curved piecewise Euclidean cell complex $C_A$ with
fundamental group $A$.

Consider a special subgroup $A' \subset A$ with the property that
if $e$ is an edge of $\G_A$ whose endpoints are in $\G_{A'}$ then
$e$ is in $\G_{A'}$. For example, $\G_A'$ could be a single edge
of $\G_A$. Then by construction $C_{A'}$ may be regarded as a
subcomplex of $C_A$, and clearly $C_{A'}$ is locally convex.

Suppose next that $A$ is an Artin group for which $\G_A$ has no
triangles. We use the construction of Brady--McCammond
\cite{BradyMcCammond:Artin} to produce the desired piecewise
Euclidean cell complex $C_A$, and to verify local convexity of the
appropriate subcomplexes. This verification is considerably more
delicate than for right angled Artin groups. The standard
presentation of a 2 generator Artin group
$$\< y_1,y_2 \suchthat (y_1,y_2)_m = (y_2,y_1)_m \>
$$
can be transformed into the presentation
\begin{align}\label{ABMpresentation}
\< d,y_1,y_2,\ldots,y_m \suchthat d = y_1 y_2, d = y_2 y_3, \ldots, d
= y_{m-1} y_m, d = y_m y_1\>
\end{align}
by triangulating the relator $(y_1,y_2)_m = (y_2,y_1)_m$ and in the
process introducing new generators $d,y_3,\ldots,y_m$
\cite{BradyMcCammond:Artin}. Note than when each $m \ge 3$, the
ordering $y_i,y_j$ is essential to the
description of the presentation (\ref{ABMpresentation}): the word
$y_i y_j$ is a subword of some relator, but the reversed word $y_j
y_i$ is not.

The presentation complex of (\ref{ABMpresentation}) has one
vertex, $1+m$ edges, and $m$ triangular faces. The link of the
unique vertex is given in Figure~\ref{ABMlink}. Note that the
vertices come in four layers: the \emph{first layer} $d$, the
\emph{second layer} $\{y_1,\ldots,y_n\}$, the \emph{third layer}
$\{\bar y_1,\ldots,\bar y_n\}$, and the \emph{bottom layer} $\bar
d$. Also, the edges come in three horizontal layers: the \emph{top
edges} connecting first to second layer vertices; the \emph{middle
layer} connecting second to third layer vertices; and the
\emph{bottom layer} connecting third to fourth layer vertices.
\begin{figure}[ht]
\centerline{
\relabelbox\small \epsfxsize 2.7truein
\epsfbox{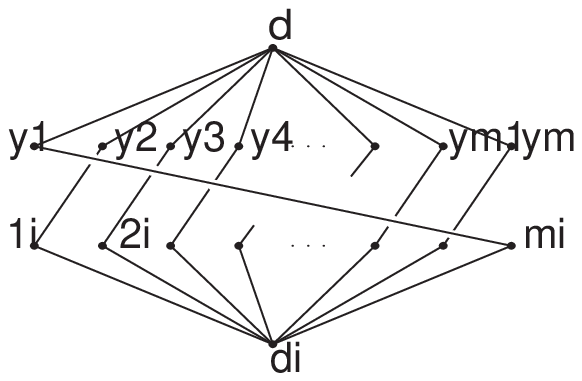}
\relabel{d}{$d$}
\relabel{y1}{$y_{1}$}
\relabel{y2}{$y_{2}$}
\relabel{y3}{$y_{3}$}
\relabel{y4}{$y_{4}$}
\relabel{ym}{$y_{m}$}
\relabel{ym1}{$y_{m-1}$}
\relabel{di}{$\bar d$}
\relabel{1i}{$\bar y_{1}$}
\relabel{2i}{$\bar y_{2}$}
\relabel{mi}{$\bar y_{m}$}
\endrelabelbox}
\caption{}
\label{ABMlink}
\end{figure}
Consider now an Artin group $A$ presented as in
(\ref{apres}). Choose an orientation on each edge of $\G_A$, which
determines an ordering of the endpoints of each edge of $\G_A$;
henceforth, when we consider the 2 generator subgroup $A_{ij} =
\<x_i,x_j \suchthat (x_i,x_j)_{m_{ij}} = (x_j,x_i)_{m_{ij}}$ we
will assume that the $ij$ edge points from $x_i$ to $x_j$. Now
rewrite the presentation (\ref{apres}) to produce the
\emph{Brady--McCammond presentation} of $A$, by triangulating each
Artin relator $(x_i,x_j)_{m_{ij}} = (x_j,x_i)_{m_{ij}}$ and
introducing new generators following the pattern of
(\ref{ABMpresentation}), where we carefully choose notation so
that new generators associated to distinct $A_{ij}$ are distinct,
as follows:
\begin{multline}\label{ABMIJ}
A_{i,j} = \langle d_{i,j}, x_i, x_j, x_{i,j,3}, x_{i,j,4} \ldots,
x_{i,j,m} \suchthat \\
d_{ij} = x_i x_j, d_{i,j} = x_j x_{i,j,3},
d_{i,j} = x_{i,j,3} x_{i,j,4} \ldots, d_{i,j} = x_{i,j,m} x_i \rangle
\end{multline}

Let $C_{ij}$ be the presentation complex for this presentation of
$A_{ij}$, and let $L_{ij}$ be the link of the unique vertex of
$C_{ij}$. The two vertex pairs $\{x_i,\bar x_i\}$ and $\{x_j,\bar
x_j\}$ in $L_{ij}$ will be called the \emph{peripheral vertex
pairs} in $L_{ij}$.

Let $C_A$ be the presentation complex for the Brady--McCammond
presentation of $A$, and note that $C_A$ is the union of its
subcomplexes $C_{ij}$ for $m_{ij} < \infinity$. Also, let $L_A$ be
the link of the unique vertex of $C_A$, and note that $L_A$ is the
union of its subcomplexes $L_{ij}$. When $m_{ij}, m_{kl} <
\infinity$ and $\{i,j\} \ne \{k,l\}$, then either $\{i,j\}
\intersect \{k,l\} = \emptyset$ in which case $C_{ij} \intersect
C_{kl}$ is the unique vertex of $C_A$ and $L_{ij} \intersect
L_{kl} = \emptyset$, or $\{i,j\} \intersect \{k,l\}$ is a
singleton, say $i=k$, in which case $C_{ij} \intersect C_{il}$ is
a single edge of $C_A$, labeled say by $x_i$, and $L_{ij}
\intersect L_{il}$ is a peripheral vertex pair, say $\{x_i,\bar
x_i\}$. It follows that the layering of vertices and edges of the
sublinks $L_{ij}$ extends to a layering of all vertices and edges
of $L_A$.

To organize $L_A$, note that there is a map $L_A$ to $\G_A$ so
that the inverse image of the $ij$ edge of $\G_A$ is $L_{ij}$, and
the inverse image of the vertex of $\G_A$ labeled $x_i$ is the
peripheral vertex pair of $L_{ij}$ labeled $x_i, \bar x_i$.

Now we use the condition that $\G_A$ has no triangles. In this
case Brady and McCammond choose a metric on $C_A$ so that each
edge labeled $d_{ij}$ has length~$\sqrt{2}$, each edge labeled
$x_i$ or $x_{ijk}$ has length~$1$, and each triangle is a
$\pi/2,\pi/4,\pi/4$ Euclidean triangle; they prove that $C_A$ is
nonpositively curved. Note that each top and bottom layered edge
in $L_A$ has spherical length $\pi/4$, and each middle layer edge
has spherical length $\pi/2$.

To verify that $C_{ij}$ is a locally convex subcomplex of $C_A$ we
must verify that for any locally injective edge path $\gamma$ in
$L_A$ with endpoints in $L_{ij}$ but with no edge in $L_{ij}$, the
spherical length of $\gamma$ is at least $\pi$.

If $\gamma$ has at least four edges then we are done. If $\gamma$
has three edges then it must connect some 2nd layer vertex to some
3rd layer vertex, and so at least one of the edges of $\gamma$ is
a middle layer vertex of length $\pi/2$, and we are done. The path
$\gamma$ cannot have one edge because $L_A$ does not have an edge
outside of $L_{ij}$ connecting a 2nd and 3rd layer vertex of
$L_{ij}$.

Suppose $\gamma$ has two edges. Since $\G_A$ has no triangles,
$\gamma$ must project to a single edge of $\G_A$ and so $\gamma$
is entirely contained in some $L_{kl}$ distinct from but
intersecting $L_{ij}$, and hence $\{i,j\} \intersect \{k,l\}$ is a
singleton. We assume that $k=i$, the other cases being handled
identically. Then $\gamma$ must connect one of the vertices
labeled $x_i,\bar x_i$ to itself. However, $L_{il}$ contains no
locally injective edge path of length two with both endpoints at
$x_i$ or both at $\bar x_i$.
\end{proof}

\paragraph{Artin groups of affine type $\wt A_n$.} Our next
verification of undistortedness uses a different method, relying
ultimately on distance estimates in mapping class groups.

\begin{thm}\label{at}
If $n \ge 3$ then the Artin group $\wt A_n$ is algebraically thick
of order at most $1$.
\end{thm}

One possible approach to proving undistortedness of special subgroups
of $\wt A_n$ is using the automatic group methods, which will be
reviewed briefly below. Charney and Peifer prove in
\cite{CharneyPeifer:affine} that $\wt A_n$ is biautomatic, and it
would suffice then to prove that the two generator special subgroups
of $\wt A_n$ are rational with respect to the Charney---Peifer
biautomatic structure. Instead we shall consider an embedding of $\wt
A_n$ into a braid group $B$, and we shall prove that all special
subgroups of $\wt A_n$ are undistorted in $B$. This trick was
suggested to us by our conversations with Ruth Charney. Our thanks to
Ruth Charney for very helpful suggestions and comments on this proof.

\begin{proof} We abbreviate $\wt A_n$ to $\wt A$, and we write its
presentation in the form
\begin{align*}
\wt A = \<x_0,x_1,\ldots,x_n &  \suchthat  x_i x_{i+1} x_i = x_{i+1}
x_i x_{i+1} \quad\text{for all}\quad i  \in \Z / (n+1)\Z, \\
& x_i x_j = x_j x_i \quad\text{for all}\quad i,j \in \Z / (n+1)\Z
\quad\text{such that}\quad j-i \not\equiv \pm 1 \>
\end{align*}
where index arithmetic takes place in $\Z / (n+1)\Z$. The cyclic
permutation of the generators $x_0,x_1,\ldots,x_n$ induces an
automorphism of $\wt A$, and this automorphism cyclically permutes
the two generator special subgroups $\<x_i,x_{i+1}\>$. It therefore
suffices to show that one of these two generator special subgroups is
undistorted.

Consider the braid group $B$ on $n+2$ strands, an Artin group with
$n+1$ generators $y_0,y_1,\ldots,y_n$ and with presentation
\begin{align*}
B = \< y_0, y_1,\ldots,y_n \quad \suchthat &\quad y_i y_j = y_j y_i
\quad\text{if}\quad 0 \le i \le j-2 \le j \le n, \\ & \quad y_i
y_{i+1} y_i = y_{i+1} y_i y_{i+1} \quad\text{if}\quad 0 \le i \le n-1
\>
\end{align*}

Let $h \from \wt A \to B$ denote the homomorphism defined on the
generators by $h(x_0)=\delta y_n \delta\inv$ and $h(x_i)=y_i$ for
$i=1,\ldots,n$, where $\delta = y_0^{-2} y_1\inv \cdots y_n\inv$.
By combining \cite{KentPeifer:annular} with the discussion at the
beginning of \cite{CharneyPeifer:affine}, it follows that $h$ is
injective. To obtain this expression for $\delta$, we refer to
\cite[Figure~4]{CharneyPeifer:affine}, which shows $\delta$ as an
element of the annular braid group on $n+1$ strands. As explained
in \cite{CharneyPeifer:affine}, this latter group is isomorphic to
the index $n+2$ subgroup of $B$ in which the $0^{\text{th}}$
strand does not move, and from this viewpoint
\cite[Figure~4]{CharneyPeifer:affine} can be redrawn as in
Figure~\ref{Figuredelta} below, which gives the desired expression
for $\delta$.

\begin{figure}[ht]
\centerline{\epsfxsize 3.5truein\epsfbox{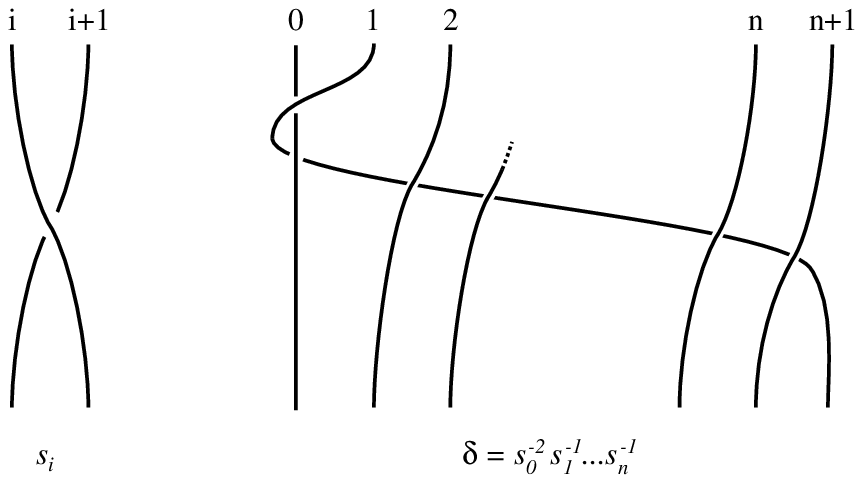}}
\caption{}
\label{Figuredelta}
\end{figure}

Clearly $h$ maps the special subgroup of $\wt A$ generated by
$x_1,x_2$ isomorphically to the special subgroup of $B$ generated
by $y_1,y_2$. It therefore suffices to show that special subgroups
in $B$ are undistorted, because of the following trick: given any
finitely generated groups $K \subgroup H \subgroup G$, if $K$ is
undistorted in $G$ then $K$ is undistorted in $H$.

The group $B$ is the mapping class group of a punctured disc $D$,
and any special subgroup of $B$ is the subgroup of mapping classes
supported on a subsurface $F \subset D$ whose boundary is a
collection of essential simple closed curves in $D$. But the fact
that the inclusion of the mapping class group of $F$ into the
mapping class group of $D$ is a quasi-isometric embedding is an
immediate consequence of
\cite[Theorem~6.12]{MasurMinsky:complex2}.
\end{proof}

\paragraph{Biautomatic groups methods.} We close this section with a
discussion of the following:

\begin{qn} \label{QnArtinThick}
Are all Artin groups algebraically thick?
\end{qn}

One of the most important problems about Artin groups is the
following:

\begin{qn} \label{QnArtinBiautomatic}
Are all Artin groups biautomatic?
\end{qn}

In an automatic or biautomatic group, a ``rational subgroup'' is a
subgroup with a particularly simple relation to the (bi)automatic
structure. See \cite{ECHLPT} for the definition; see also
\cite{GerstenShort:rational} for a detailed study of rational
subgroups of biautomatic groups. The key fact that we propose
using is:

\begin{thm}\label{ru} \cite[Theorems 3.3.4 and 8.3.1]{ECHLPT}
\label{ThmBiautUndist}
If $G$ is an automatic group and $H$ is a rational subgroup then $H$
is undistorted in $G$. \qed\end{thm}

By using Theorem~\ref{ru}, an affirmative answer to
Question~\ref{QnArtinThick} reduces to an affirmative answer to
the following refinement of Question~\ref{QnArtinBiautomatic}:

\begin{qn}
Does every Artin group have a biautomatic structure so that every
special subgroup is rational?
\end{qn}

What makes this a reasonable question to pursue are the many special
classes of Artin
groups known to be biautomatic (see the earlier list). It would be
interesting to check some of these classes for rationality of special
subgroups. For example, we have checked that all special subgroups of
the braid group $B$ are rational with respect to the symmetric
biautomatic structure defined in
\cite[Theorem~9.3.6]{ECHLPT} --- incidentally, this would provide
another proof of Theorem~\ref{at}, but to save space we opted for
quickly quoting the Masur--Minsky distance estimates for mapping
class groups.

\section{Fundamental groups of graph manifolds}\label{Sectiongrm}

Recall that graph manifolds are compact Haken manifolds of zero
Euler characteristic, such that all their geometric components are
Seifert manifolds. In this section when referring to graph
manifolds we always assume that they are connected, and we rule
out the case of Nil and Sol manifolds. Hence all graph manifolds
we consider are obtained by gluing finitely many Seifert
components with hyperbolic base orbifolds along boundary tori or
Klein-bottles, where the gluing does not identify the fibers.

In the universal cover $\widetilde{M}$ of such a manifold $M$, a
flat projecting onto a torus or Klein bottle along which different
Seifert components are glued is called a \textit{separating flat}.
A copy of a universal cover of a Seifert component is called a
\textit{geometric component}. Note that separating flats bound and
separate geometric components.

A particular case of graph manifolds are the \textit{flip
manifolds}, in the terminology of \cite{KapovichLeeb:3manifolds}.
Each Seifert component of a flip manifold is the product of a
compact, oriented surface-with-boundary (the base) and $S^1$ (the
fiber). Wherever two Seifert components are glued along a boundary
torus the gluing interchanges the basis and the fiber directions.

Every flip manifold admits a non-positively curved metric, as
follows. For each Seifert component, put a hyperbolic metric with
geodesic boundary on the base so that each boundary component has
length~1, pick a metric on the fiber to have length~1, and use the
Cartesian product metric; each gluing map of boundary tori is then
an isometry. Note that each Seifert component is locally convex.

Not every graph manifold admits a nonpositively curved metric
\cite{Leeb:nonpos}. On the other hand, according to
\cite{KapovichLeeb:3manifolds}, the fundamental group of any graph
manifold is quasi-isometric to the fundamental group of some flip
manifold. Moreover, the induced quasi-isometry between the
universal covers of the two manifolds preserves the geometric
decomposition, namely, the image of any geometric component is a
uniformly bounded distance from a geometric component. We prove
the following.

\begin{thm}
The fundamental group $G=\pi_1(M)$ of a non-Seifert fibered graph
manifold is
algebraically thick of order~$1$ with respect to the family
$\mathcal{H}$ of fundamental groups of its geometric (Seifert)
components.\footnote{This result was suggested to us by Kleiner,
who had observed that these groups are not relatively hyperbolic
with respect to any collection of finitely generated subgroups
\cite{Kleiner:personal}.}
\end{thm}

\begin{rmk}
Note that $G$ is weakly hyperbolic relative to the family
$\mathcal{H}$, because for any finite graph of finitely generated
groups, the fundamental group is weakly hyperbolic relative to the
vertex groups.
\end{rmk}

\begin{proof} By \cite[Theorem
1.1]{KapovichLeeb:3manifolds}, $G$ is quasi-isometric to the
fundamental group of a compact non-positively curved flip manifold
with totally geodesic flat boundary, and the images under a
quasi-isometry of subgroups in $\H$ are a bounded distance from
fundamental groups of geometric components. Since Seifert components
of flip manifolds are locally convex, an application of
Proposition~\ref{PropNPC}
shows that the
subgroups in $\H$ are undistorted in $G$.

Any subgroup $H$ in $\H$ has a finite index subgroup $H_1$ which
is the fundamental group of a trivial circle bundle over an
orientable surface of genus at least two. Thus $H$ is
quasi-isometric to the direct product of $\R$ with a convex subset
in $\mathbb{H}^2$. Consequently any asymptotic cone of $H$ is
bi-Lipschitz equivalent to a direct product of an $\R$-tree with
$\R$, therefore $H$ is unconstricted.

The group $G$ decomposes as a fundamental group of a graph of
groups, with vertex groups in $\H$ and edge groups commensurable
to $\Z^2$. Any two subgroups in $\H$ can be thickly connected
using a path in this graph of groups. We conclude that $G$ is
algebraically thick of order $\le 1$ with respect to~$\H$. In
\cite{KapovichLeeb:3manifolds} it is proven that the fundamental
group of a non-Seifert fibered graph manifold manifolds has
superlinear divergence. Hence, $G$ is constricted by
Lemma~\ref{suplin}. Thus we have that $G$ is algebraically thick
of order $1$ with respect to~$\H$.
\end{proof}

\begin{rmk}\label{rthick2}
Using graph manifolds one can construct an example of a group that is
thick of order two but not of order zero or one. Indeed, consider a
manifold $N$ obtained by doubling a compact flip manifold $M$ along a
periodic geodesic $g \subset M$ that is not contained in a Seifert
component of $M$. The fundamental group of $N$ is algebraically thick
of order $2$, but not of order $0$ or $1$. (Details of this
construction and further results will be provided in
\cite{BehrstockDrutu:thick2}.)
\end{rmk}

\section{Teichm\"{u}ller space and the pants graph}
\label{Section:teich}

Let $S=S_{g,p}$ be an orientable surface of genus $g$ with $p$
punctures, with complexity $\xi=\xi(S)=3g+p-3$. We let $\pants$
denote the pants graph of $S$, defined as follows. A vertex of
$\pants$ is an isotopy class of pants decompositions of $S$. Given
a pants decomposition $T = \{\gamma_1,\ldots,\gamma_\xi\}$,
associated to each $\gamma_i$ is the unique component $R_i$ of $S
- \cup\{ \gamma_j \suchthat 1 \le j \le \xi, j \ne i \}$ which is
not a pair of pants. This subsurface $R_i$ has complexity~$1$, it
is either a once punctured torus or a 4 punctured sphere, and we
refer to $R_i$ as a \emph{complexity~$1$ piece} of $T$.  Two pairs
of pants decompositions $T,T'\in\pants$ with
$T=\{\gamma_{1},\ldots,\gamma_{\xi}\}$ and
$T'=\{\gamma'_{1},\ldots,\gamma'_{\xi}\}$ are connected by an edge
of $\pants$ if they differ by an \emph{elementary move}, which
means that $T$ and $T'$ can be reindexed so that, up to isotopy,
the following conditions are satisfied:
\begin{enumerate}
\item $\gamma_{i}=\gamma'_{i}$ for all $2\leq i\leq \xi$
\item Letting $X$ be the common complexity~$1$ piece of $T$ and $T'$
associated to both $\gamma_1$ and $\gamma'_1$, we have
$d_{\C(X)}(\gamma^{\vphantom\prime}_{1},\gamma'_{1})=1$.
\end{enumerate}

We now recall two results concerning the pants complex. The first
relates the geometry of the pants complex to that of
Teichm\"{u}ller space, see \cite{Brock:wp}.

\begin{thm}[Brock]\label{WPequalsPants}
$\pants$ is quasi-isometric
to Teichm\"{u}ller space with the Weil-Petersson metric.
\end{thm}

Bowditch \cite{Bestvina:list2004} asked whether Teichm\"{u}ller space
with the
Weil-Petersson metric was a $\delta$--hyperbolic metric space.
This was first answered by Brock and Farb in
\cite{BrockFarb:curvature} where they showed:

\begin{thm}[Brock--Farb]\label{pantsqiflats}
$\R^{n}$ can be
quasi-isometrically embedded into $\pant(S)$ for all $n\leq
\lfloor \frac{\xi(S)+1}{2}\rfloor$.
\end{thm}

Combined with Theorem~\ref{WPequalsPants}, this result showed the
answer to Bowditch's question is ``no'' when the surface is
sufficiently complex, i.e.,
satisfies $\xi(S)>2$. By contrast, Brock--Farb gave an affirmative
answer to Bowditch's question when $\xi(S)\leq 2$ (see also
\cite{Aramayona:thesis} and \cite{Behrstock:asymptotic} for
alternative proofs of $\delta$--hyperbolicity in these low
complexity cases). Brock--Farb proved hyperbolicity by showing
that these Teichm\"{u}ller spaces are strongly relatively
hyperbolic with respect to a collection of subsets which
themselves are hyperbolic; they then showed that this implies
hyperbolicity of Teichm\"{u}ller space. This raises the question
of whether the presence of relative hyperbolicity generalizes to
the higher complexity cases: we show that it does not, for
sufficiently high complexity.

\begin{thm}\label{WPnotatg} For any surface $S$ of finite type with
$\xi(S)\geq 6$, Teichm\"{u}ller space with the Weil-Petersson
metric is not asymptotically tree-graded as a metric space with
respect to any collection of subsets.
\end{thm}

This is particularly interesting in light of the following:

\begin{thm}[Behrstock \cite{Behrstock:asymptotic}]
    \label{teichtreegraded} For every surface
    $S$, every asymptotic cone of
    Teichm\"{u}ller space with the Weil-Petersson metric is
    tree-graded.
\end{thm}

Together Theorems~\ref{WPnotatg} and~\ref{teichtreegraded} say
that the pieces in the tree-graded structure of an asymptotic cone
do not merely arise by taking ultralimits of a collection of
subsets.

\medskip

Theorem~\ref{WPnotatg} will follow from Theorem~\ref{crit}, once we establish:

\begin{thm}\label{teichthick} For any surface $S$ of finite type with
$\xi(S)\geq 6$, Teichm\"{u}ller space with the Weil-Petersson
metric is metrically thick of order one.
\end{thm}

\begin{pf} Let $S$ denote a surface with $\xi(S)\geq 6$.
Brock--Farb proved Theorem~\ref{pantsqiflats} by explicitly
constructing quasiflats of the desired dimension. We shall use
these same quasiflats to prove thickness, so we now recall their
construction. Cut $S$ along a pairwise disjoint family of simple
closed curves into a collection of subsurfaces each of which is
either a thrice punctured sphere or of complexity~1. Let $\Rcal =
\{R_1, \ldots, R_k\}$ be the subcollection of complexity~1
subsurfaces, and we assume that $k \ge 2$. Note that this is
possible for a given $k$ if and only if $2 \le k \le \lfloor
\frac{\xi(S)+1}{2}\rfloor$. For each $i$, let $g_{i}$ denote a
geodesic in the curve complex $\C(R_{i})$. One obtains a pants
decomposition of $S$ by taking the union of the curves $\partial
R_{i}$ and one curve from each $g_{i}$.

Theorem~\ref{pantsqiflats} is proven by showing that, for a fixed
collection of subsurfaces and geodesics as above, the collection
of all such pants decompositions is a quasiflat of rank $k$. If
$\g$ denotes the family of geodesics $\{g_1,g_2,...,g_k\}$, the
above quasiflat is denoted by $Q_{\Rcal,\g}$. For a fixed surface
$S$, all such quasi-isometric embeddings of $\R^{k}$ have
uniformly bounded quasi-isometry constants. Let $\lll$ be the
collection of all quasiflats $Q_{\Rcal,\g}$. Note that when
$\xi(S)\geq 6$ every element of $\pants$ is contained in some
$Q_{\Rcal,\g}$, thus this collection satisfies condition
$(\bn_{1})$ of metric thickness. Further, since each
$Q_{\Rcal,\g}$ is a quasiflat of dimension at least two with
uniform quasi-isometry constants, this collection is uniformly
unconstricted.

It remains to verify condition $(\bn_2)$. We proceed in two steps:
\begin{enumerate}
\item Any pair $T,T'\in\pants$ which differ by an
elementary move lie in some quasiflat in $\lll$.
\item Any pair of quasiflats in $\lll$ which intersect
can be thickly connected in $\lll$.
\end{enumerate}

Since the pants complex is connected by elementary moves, the
first step implies that given any two pants decompositions
$T,T'\in\pants$, one can find a sequence of quasiflats in $\lll$
each intersecting the next in at least one point, such that the
first quasiflat contains $T$ and the last contains $T'$. The
second step then implies that this sequence is a subsequence of a
sequence of quasiflats in $\lll$ where each intersects its
successor in an infinite diameter set. This establishes condition
$(\bn_2)$.

\medskip

\emph{(Step 1).} Fix two pair of pants decompositions
$T,T'\in\pants$ which differ by an elementary move. This
elementary move is supported in a subsurface $R_1$ with
$\xi(R_1)=1$. Since $\xi(S)\geq 6$ there exists a curve $\alpha$
of $T$ and $T'$ disjoint from $R_1$. Let $R_2$ be the union of the
pants of $T$ and $T'$ on either side of $\alpha$, so $R_1,R_2$
have disjoint interior and $\xi(R_2)=1$. Let $g_1\subset \C(R_1)$
be an infinite geodesic extending the elementary move in $R_1$.
The product of $g_1$ with a geodesic $g_2$ supported on $R_2$ is a
two-dimensional quasiflat $Q_{\Rcal,\g}$, $\Rcal=\{R_1,R_2\}$,
$\g=\{g_1,g_2\}$, an element of $\lll$, containing both $T$ and
$T'$.

\medskip

\emph{(Step 2).} Consider $Q=Q_{\Rcal,\g}$ and  $Q'=Q'_{\Rcal',\g
'}$ an arbitrary pair of quasiflats in $\lll$, with non-empty
intersection containing $T \in \pants$.

\smallskip

\emph{(a).} Assume first that there exists $R \in \Rcal \intersect
\Rcal'$. Let $g\in\g$ and $g'\in\g'$ denote the corresponding
geodesics in $Q$ and $Q'$, respectively. Consider
the quasiflat $Q''=Q''_{\Rcal,\g''}$, where $\g''$ is obtained from
$\g$ by replacing $g$ by $g'$. Then $Q''$ has infinite diameter
intersection with both $Q$ and $Q'$ thus providing the desired
thickly connecting sequence.

\smallskip

\emph{(b).} Fixing $T \in \pants$, define a finite 1-complex as
follows. A vertex is a collection of pairwise disjoint subsurfaces
$\Rcal = \{R_1,\ldots,R_k\}$ with cardinality $k \ge 2$, such that
each $R_i$ is a complexity~1 piece of $T$. Two such collections
$\Rcal,\Rcal'$ are connected by an edge if they are not disjoint. By
part \emph{(a)} it suffices to show that this 1-complex is connected,
and we prove this using that $\xi(S) \ge 6$. Consider two vertices
$\Rcal,\Rcal'$. By removing elements of each we may assume that each
has cardinality~2. Let $\Rcal = \{R_1,R_2\}$,
$\Rcal'=\{R'_1,R'_{2}\}$. We may also assume that $\Rcal,\Rcal'$ are
disjoint. If some element of $\Rcal$ is disjoint from some element of
$\Rcal'$, say $R_1 \intersect R'_1 = \emptyset$, then both are
connected by an edge to $\{R_1,R'_1\}$; we may therefore assume that
$R_i \intersect R'_j \ne \emptyset$ for $i,j=1,2$. If some element of
$\Rcal$ or $\Rcal'$ is a once-punctured torus, say $R_1$ with
boundary curve $\alpha$, then the only possible element of $\Rcal'$
that can intersect $R_1$ is the one obtained by removing $\alpha$,
contradicting that there are two elements of $\Rcal'$ that intersect
$R_1$; we may therefore assume that each $R_i$ and each $R'_j$ is a
four punctured sphere. It follows now that $T$ has four distinct
pairs of pants $P_1,P_2,P_3,P_4$ such that $R_1 = P_1 \union P_2$,
$R'_1 = P_2 \union P_3$, $R_2 = P_3 \union P_4$, $R'_2 = P_4 \union
P_1$. Since $\xi(S) \ge 6$, there is a curve $\alpha$ of $T$ not
incident to any of $P_1,\ldots,P_4$, and letting $R''$ be the
complexity~1 piece of $T$ obtained by removal of $\alpha$, it follows
that $\Rcal = \{R_1,R_2\}$ is connected by an edge to
$\{R_1,R_2,R''\}$ which is connected to $\{R'_1,R'_2,R''\}$ which is
connected to $\{R'_1,R'_2\}$.

\bigskip

We have now shown that Teichm\"{u}ller space with the
Weil-Petersson metric is thick of order at most~1 when $\xi(S)\geq
6$. That it is thick of order exactly~1 follows from
Theorem~\ref{teichtreegraded}.
\end{pf}

\begin{rmk}
With a little more work one can prove Step~\emph{(2)(b)} under the
weaker assumption that $\xi(S) \ge 5$. Condition $(\bn_{1})$ can
also be proved when $\xi(S) \ge 5$: the proof of $(\bn_{1})$ given
here has an unnecessarily strong conclusion, namely that each
point of $\pants$ lies in the union of $\lll$. However, we do not
know how to weaken the proof of Step~\emph{(1)} for any case when
$\xi(S) < 6$.
\end{rmk}

\begin{rmk}
The surfaces with $3\leq \xi(S)\leq 5$, 
do not fall under the cases where
Teichm\"{u}ller space is hyperbolic (see \cite{BrockFarb:curvature},
or for alternate arguments,
\cite{Aramayona:thesis} or \cite{Behrstock:asymptotic}) or under
the cases of Theorem~\ref{teichthick} where Teichm\"{u}ller space is
metrically thick and hence not relatively hyperbolic. 
Accordingly, in the preprint
versions of this paper (Decemeber 2005) we explicitly asked
whether in these cases the Weil-Petersson
metric on Teichm\"{u}ller space is relatively hyperbolic. The situation
for these remaining cases has now been resolved by work of Brock and
Masur \cite{BrockMasur:WPrelhyp}.
Brock and Masur prove that when $\xi(S)=3$, then
the Weil-Petersson metric is relatively hyperbolic.  On the other hand,
when $\xi(S)=4$ or $5$ they showed that the Weil-Petersson metric is
thick. More precisely,  when  $\xi(S)=4$ or $5$,
except $S=S_{2,1}$, then this space is thick
of order 1, while for $S_{2,1}$ this metric is thick of order
at most 2.
\end{rmk}

The following interesting question remains open:

\begin{qn} For $S_{2,1}$ is the Weil-Petersson
   metric on Teichm\"{u}ller space thick of degree 1 or 2?
\end{qn}

\section{Subsets in symmetric spaces and lattices}\label{lattices}

\subsection*{Subsets in symmetric spaces}

Let $X$ be a
product of finitely many nonpositively curved symmetric spaces and Euclidean buildings of rank at
least two, and let $\dist_X$ be a product metric on it (uniquely
defined up to rescaling in the factors). Given a geodesic ray $r$
in $X$, {\it the Busemann function associated to $r$} is defined
by
$$
f_r:X\to \R \, ,\; f_r(x)=\lim_{t\to \infty}[\dist_X (x,r(t))-t]\;
.
$$

\begin{rmk}\label{diff}
The Busemann functions of two asymptotic rays in $X$ differ by a
constant \cite{BridsonHaefliger}.
\end{rmk}

The level hypersurface $H(r)=\lbrace x\in X \; ;\; f_r(x)= 0
\rbrace$ is called the {\it horosphere determined by $r$}, the
sublevel set $Hbo(r)=\lbrace x\in X \; ;\; f_r(x)<0 \rbrace$ is
the {\it open horoball determined by $r$}.

\begin{prop}\label{Thick subsets symmetric space}
Let $\Rcal$ be a family of geodesic rays in $X$, such that no ray is contained in a rank one factor of
$X$ and such that the open horoballs in the family $\{ Hbo(r)\mid
r\in \Rcal\}$ are pairwise disjoint. Then for any $M>0$, any
connected component $C$ of $\bigcup_{r\in \Rcal }\nbhd_M(H(r))$
endowed with $\dist_X$ is $M$-thick of order one with respect to
$$\lll = \{ H(r)\mid r\in \Rcal, H(r) \subset C\}\, .$$
\end{prop}

\begin{pf} The fact that $\{ Hbo(r)\mid r\in \Rcal\}$ are
pairwise disjoint implies that their basepoints $\{ r(\infty )\mid
r\in \Rcal\}$ are pairwise opposite. If the previous set has
cardinality at least three, then according to \cite[proof of
Proposition 5.5, b]{Drutu:Nondistorsion} all the rays in
$\Rcal$ are congruent under the action of the group $\mathrm{Isom} (X)$.
Hence all horospheres $H(r)$ with $r\in \Rcal$ are isometric if $\Rcal$ has
cardinality at least three. Thus
in order to have property $(\theta)$ in \textbf{(M$_2$)} it suffices to prove that one such
horosphere endowed with the restriction of $\dist_X$ is
unconstricted.

Let $H$ be such a horosphere. According to
\cite[Lemma 4.2]{Drutu:Nondistorsion}, any
asymptotic cone $H_\infty = \mathrm{Cone}_\omega (H,h,d)$ is a
horosphere in the asymptotic cone $X_\infty =\mathrm{Cone}_\omega
(X,h,d)$. The cone $X_\infty$ is a Euclidean building having the
same rank as $X$ \cite{KleinerLeeb:buildings}. Let $r_\infty$ be
the ray in $X_\infty$ such that $H_\infty = H(r_\infty)$. The
hypothesis that rays in $\Rcal$ are not contained in a rank one
factor of $X$ implies that $r_\infty$ is not contained in a rank
one factor of $X_\infty$. According to Proposition 3.1.1 and Lemma 3.3.2 in \cite{Drutu:Filling}, under this hypothesis any two points $x,y$ in $H_\infty$ can be
joined by a pair of topological arcs in $H_\infty$ intersecting
only in their endpoints. This is mainly due to three facts. Firstly, a horoball $Hb(r_\infty)$ intersects a maximal flat (or apartment) in $X_\infty$ in a convex polytope, and the horosphere $H(r_\infty)$ intersects that flat in the boundary polytopic hypersurface \cite[Proposition 3.1.1]{Drutu:Filling}. Secondly, any two points $x,y\in H(r_\infty)$ are also contained in a maximal flat $F$ in $X_\infty$, by the axioms of a Euclidean building. Thirdly, given a maximal flat $F$ in $X_\infty$ intersecting $H(r_\infty)$ and two points $x,y\in F\cap H(r_\infty)$ one may ensure, by replacing a half-flat of $F$ if necessary, that $F\cap H(r_\infty)$ is a finite polytope \cite[Lemma 3.3.2]{Drutu:Filling}. We conclude that $H_\infty$ cannot have
cut-points, hence that $H$ is unconstricted.

The fact that $C$ is connected implies that for every $H(r),
H(r')\in \lll$ there exists a finite sequence $r_1=r,
r_2,...,r_n=r'$ such that $\nbhd_M(H(r_{i}))\cap
\nbhd_M(H(r_{i+1}))$ is non-empty. Then
$\dist(H(r_{i}),H(r_{i+1}))\leq 2M$. There exists a maximal flat
$F$ in $X$ containing rays asymptotic to both $r_i$ and $r_{i+1}$.
Remark \ref{diff} implies that one may suppose that both $r_i$ and
$r_{i+1}$ are contained in $F$. Since $r_i(\infty )$ and
$r_{i+1}(\infty )$ are opposite, $F\cap H(r_{i})$ and $F\cap
H(r_{i+1})$ are two parallel hyperplanes, at distance at most
$2M$. It follows that $\nbhd_M(H(r_{i}))\cap \nbhd_M(H(r_{i+1}))$
has infinite diameter. Thus $H(r)$ and $H(r')$ are thickly
connected by the sequence $H(r_1), H(r_2),..., H(r_n)$.
\end{pf}

\subsection*{Higher rank lattices}

Particularly interesting examples of spaces $C$ satisfying the
hypotheses of Proposition~\ref{Thick subsets symmetric space} are
those on which some $\Q$--rank one lattice acts cocompactly. In
this case, the space $C$ is quasi-isometric to the lattice, and
one can prove more than just metric thickness.

We recall first some known facts about lattices. In rank one
semisimple groups, uniform lattices are hyperbolic, while
non-uniform lattices are relatively hyperbolic with respect to
their maximal unipotent subgroups (this in particular implies that
maximal unipotent subgroups are undistorted in the lattice). Thus
in both cases lattices cannot be thick.

In higher rank semisimple groups, uniform lattices have as
asymptotic cones Euclidean buildings of higher rank
\cite{KleinerLeeb:buildings} so they are unconstricted, thus in
particular they are thick of order zero.

In what follows we prove that non-uniform lattices in higher rank
semisimple groups are algebraically thick of order at most one. In
our arguments we also use unipotent subgroups. Unlike in the rank
one case, these subgroups are distorted in the ambient lattice,
therefore we have to embed them into solvable undistorted
subgroups of the lattice, in order to prove thickness. For details
on the notions and the results mentioned in this section we refer
to \cite{BorelTits:greductifs}, \cite{Borel:IntroArith},
\cite{Margulis:DiscreteSubgroups} and \cite{Witte:Arithmetic}.

Let $G$ be a connected semisimple group. Then $G$ has a unique
decomposition, up to permutation of factors, as an almost direct
product $G=\prod_{i=1}^m \bg_i (k_i)$, where $k_i$ is a locally
compact non-discrete field and $\bg_i (k_i)$ is a connected group
almost simple over $k_i$. Recall that:
\begin{enumerate}
    \item[-] An algebraic group defined over a field $k$ is called
\textit{almost simple over }$k$ if all the proper $k$-algebraic
normal subgroups of it are finite.
    \item[-] An algebraic group is called
\textit{absolutely almost simple} if any proper algebraic normal
subgroup of it is finite.
    \item[-] An algebraic group $G$ is an \textit{almost direct product} of
    its subgroups $H_1,..,H_m$ if the multiplication map
    $H_1 \times ...\times H_m \to G$
    is surjective and of finite kernel (an \textit{isogeny}).
\end{enumerate}

The \textit{rank} of $G$ is defined as $\mathrm{rank}\, G\, =\,
\sum_{i=1}^m \mathrm{rank}_{k_i}\bg_i$, where
$\mathrm{rank}_{k_i}\bg_i$ is the dimension of the maximal
$k_i$--split torus of $\bg_i$. Recall that a \emph{$k_i$--split
torus} is a subgroup of $G$ defined over $k_{i}$ which is abelian, closed, connected,
with every element semisimple, and which is diagonalizable over
the field $k_{i}$.
We make the following two assumptions on $G$:
\begin{itemize}
\item[(\textbf{Hyp}$_1$)]  For every
$i=1,2,...,m$, $\mathrm{char}\, k_i=0$.

\medskip

\item[(\textbf{Hyp}$_2$)] For every $i=1,2,...,m$,
$\mathrm{rank}_{k_i}\bg_i \geq 1$, and $\mathrm{rank}\, G \geq 2$.
\end{itemize}

\medskip

\noindent \textit{Notation:} Given two functions $f,g$ defined on a set $X$
and taking real values, we write $f\ll g$ if $f(x)\leq C g(x)$ for
every $x\in X$, where $C$ is a constant uniform in $x$. We write
$f\asymp g$ if $f\ll g$ and $g\ll f$.

\medskip

The group $G$ can be endowed with a left invariant metric
$\dist_G$ with the property that for fixed embeddings of each
$\bg_i (k_i)$ into $SL(n_i, k_i)$,
\begin{equation}\label{log}
\dist_G(1,g)\asymp \sum_{i=1}^m \ln \left(1+ \|
g_i-\mathrm{Id}_i\|_{i,max} \right)\, .
\end{equation}

See for instance \cite{LMR:metrics} for details.

Let $\Gamma$ be an \emph{lattice} in $G$, that is, a discrete
subgroup of $G$ for which $G/\Gamma$ carries a $G$--left invariant
finite measure. If the projection of $\Gamma$ to any direct factor
of $G$ is dense then the lattice is called \emph{irreducible}.
Otherwise it is commensurable to a product $\Gamma_1 \times
\Gamma_2$, where $\Gamma_i$ are lattices in direct factors of $G$.
Note that in this latter case, $\Gamma $ is unconstricted,
according to the first example following Definition \ref{duunc}.
Therefore, in what follows we always assume that $\Gamma$ is
irreducible.

The lattice $\Gamma$ is called \textit{uniform} if $G/\Gamma$ is
compact.
 Throughout the rest of the section we assume that
$\Gamma$ is a \textit{non-uniform lattice}, that is $G/\Gamma$ is
not compact.

\begin{thm}[Lubotzky--Mozes--Raghunathan \cite{LMR:metrics}, Theorem
A]\label{tlmr}
The word metric on $\Gamma$ is bi-Lipschitz equivalent to
$\dist_G$ restricted to $\Gamma$.
\end{thm}

By Margulis' Arithmeticity Theorem \cite[Chapter
IX]{Margulis:DiscreteSubgroups}, the hypotheses that
$\mathrm{rank}\, G \geq 2$ and that $\Gamma$ is irreducible imply
that $\Gamma$ is an $S$--arithmetic group: there exists a global
field $F$, a simply connected absolutely almost simple algebraic
group $\bg$ defined over $F$, a finite set $\bs$ of valuations of
$F$ containing the archimedean ones and a homomorphism $\Phi :
\prod_{v\in \bs} \bg(F_v)\to G$ such that $\ker \Phi$ is compact,
$\mathrm{Im}\, \Phi$ is a closed normal subgroup of $G$ with
$G/\mathrm{Im}\, \Phi$ also compact, and $\Gamma$ is commensurable
with $\Phi(\bg (\oo_\bs))$, where $\oo_\bs$ is the ring of
$\bs$--integers in $F$, defined by $|\cdot |_v \leq 1$ for every
$v\not\in \bs$. Here $\bg (\oo_\bs)$ is $\bg (F) \cap
\mathrm{GL}(n, \oo_\bs)$ if we assume that $\bg$ is an
$F$-algebraic subgroup of $\mathrm{GL}(n)$, and we identify $\bg
(\oo_\bs)$ to its image under the diagonal embedding in
$\prod_{v\in \bs} \bg(F_v)$. The hypothesis that $\Gamma$ is
non-uniform is equivalent to the property that $\mathrm{rank}_F
\bg \geq 1$.

\begin{thm}\label{tlattices}
The lattice $\Gamma$ is algebraically thick of order at most one.
\end{thm}

\begin{pf} It suffices to prove the statement for $G=\prod_{v\in S}
\bg(F_v)$ and $\Gamma = \bg (\oo_S)$, where $\bg (\oo_S)$ is
identified to its image under the diagonal embedding in $G$.
We first recall some useful notions and results. A reductive group
is \emph{$F$--anisotropic} if it is defined over $F$ and if it
does not contain any non-trivial $F$--split torus.

\begin{lem}\label{red}
A reductive subgroup $R$ of $G$ which is defined over $F$
intersects $\Gamma$ in a lattice, that is $R/R\cap \Gamma$ has
finite measure. Moreover if $R$ is $F$--anisotropic then the
lattice is uniform, that is $R/R\cap \Gamma$ is compact.
\end{lem}

\medskip

Let $P$ be a parabolic subgroup of $G$ defined over $F$. Then we have the
following:
\begin{enumerate}
    \item The \emph{unipotent radical} $U$ of $P$ (i.e. the unique
maximal
unipotent normal subgroup of $P$) is algebraic defined over $F$.
    \item There exists an $F$--split torus $T$ such that $P=Z(T) U$,
where $Z(T)$ is the centralizer of $T$ in $G$.
    \item In its turn, $Z(T)$ is an almost direct product of a torus
$T'$ with $M=[Z(T), Z(T)]$, and the group $M$ is semisimple.
    \item The torus $T'$ is an almost direct product of $T$ with an
$F$--anisotropic torus~$C$.
    \item The groups $Z(T),M, T'$ are algebraic defined over $F$.
    \item  The
group $M$ contains countably many $F$--anisotropic tori $D$ each of
which
is maximal in $M$ and contains for each $v\in S$ a maximal
torus of $M$ defined over $F_v$ \cite[Lemma 3.10]{LMR:metrics}.
For each such torus $D$, the product $TCD$ is a maximal torus in
$G$.
     \item When the $F$-parabolic group $P$ is minimal, $T$ is a
maximal
$F$-split torus, $U$ is a maximal $F$-unipotent subgroup of $G$,
and $U\cap \Gamma$ is a maximal unipotent subgroup in $\Gamma$.
    \item Conversely, given an $F$--split torus $T$ there exist
finitely
many $F$--parabolic subgroups $P$ that can be written as $P=Z(T)
U$, and such that all the above decompositions and properties
hold. These parabolic subgroups correspond to finitely many faces of Weyl chambers composing $T$.
\end{enumerate}

Let $T$ be a maximal $F$-split torus in $G$, and let $T'=TC$, $D<
M$ and $\widetilde{T}=TCD$ be tori associated to $T$ as above. Let
$\Lambda$ be the system of $F$-roots of $G$ with respect to $T$,
and let $\widetilde{\Lambda}$ be the system of roots of $G$ with
respect to $\widetilde{T}$. For every $\widetilde{\lambda} \in
\widetilde{\Lambda}$, if its restriction to $T$,
$\widetilde{\lambda}|_{T}$, is not constant equal to $1$ then it
is in $\Lambda$.

In what follows all bases of roots and lexicographic orders on
roots will be considered as chosen simultaneously on both
$\Lambda$ and $\widetilde{\Lambda}$ so that they are compatible
with respect to the restriction from $\widetilde{T}$ to $T$.

Let $\fg$ be the Lie algebra of $G$. For every
$\widetilde{\lambda}$ in $\widetilde{\Lambda} $ denote by
$\fg_{\widetilde{\lambda}}$ the one dimensional eigenspace
$\left\{ v\in \fg \; ;\; Ad (t) (v)= \widetilde{\lambda} (t) v\,
,\, \forall t\in \widetilde{T} \right\}$. Here $Ad (t)$ is the
differential at the identity of the conjugacy by $t$. There is a
unique one-parameter unipotent subgroup $U_{\widetilde{\lambda}}$
in $G$ tangent to the Lie algebra $\fg_{\widetilde{\lambda}}$. Let
$\lambda \in \Lambda$. A multiple of it $n\lambda$ with $n\in \N$
can be in $\Lambda$ for $n\in \{ 1,2\}$. Consider the Lie
subalgebra ${\mathfrak{u}}_\lambda =
\bigoplus_{\widetilde{\lambda}|_{T}=\lambda, 2\lambda}
\fg_{\widetilde{\lambda}}$ and let $U_\lambda$ denote the unique
$T$-stable $F$-unipotent subgroup of $G$ tangent to the Lie
algebra ${\mathfrak{u}}_\lambda$. Let $\Delta$ be a basis for
$\Lambda$ (or, in another terminology, a fundamental system of
roots).
Every root $\lambda$ in $\Lambda$ can be written as $\sum_{\alpha
\in \Delta } m_\alpha (\lambda )\, \alpha$, where $(m_\alpha
(\lambda ))_{\alpha \in \Delta}$ are integers either all
non-negative or all non-positive.

Let $P$ be an $F$-maximal parabolic subgroup. There exists a
maximal $F$-split torus $T$, a basis $\Delta$ for the system
$\Lambda$ of $F$-roots of $G$ with respect to $T$, and a root
$\alpha \in \Delta $ such that the following holds. Let
$\Lambda_\alpha^+ =\left\{ \lambda\in \Lambda \; ;\; m_\alpha
(\lambda) >0 \right\}$. The parabolic $P$ decomposes as
$P=Z(T_\alpha ) \widetilde{U}_\alpha$, where
\begin{itemize}
    \item  $T_\alpha =\{t\in T \; ;\; \beta
(t)=1\, ,\, \forall \beta \in \Delta\, ,\, \beta \neq \alpha \}$;
    \item $\widetilde{U}_\alpha$ is the unipotent subgroup tangent to the Lie
algebra $\widetilde{{\mathfrak{u}}}_\alpha =
\bigoplus_{\widetilde{\lambda}|_T \in
\Lambda_\alpha^+}\fg_{\widetilde{\lambda}}$. Note that
$\widetilde{{\mathfrak{u}}}_\alpha = \bigoplus_{\lambda \in
\Lambda_\alpha} {\mathfrak{u}}_\lambda $, where $\Lambda_\alpha$
is such that any root in $\Lambda_\alpha^+$ is either contained in
$\Lambda_\alpha$ or is the double of a root in $\Lambda_\alpha$.
In particular, the above implies that $U_\alpha$ is a subgroup of
$\widetilde{U}_\alpha$.
\end{itemize}

The following result is proved in \cite{Raghunathan:generators}
for $\mathrm{rank}_F \bg\geq 2$ and in
\cite{Venkataramana:generators} for $\mathrm{rank}_F \bg=1$.

\begin{thm}[Raghunathan; Venkataramana]
Let $T$ be a maximal $F$-split torus in $G$ and let $\Lambda$ be
the system of $F$-roots of $G$ with respect to $T$. Then the group
generated by the subgroups $U_\lambda \cap \Gamma$, $\lambda \in
\Lambda$, has finite index in $\Gamma$.
\end{thm}

Note that when $\mathrm{rank}_F \bg=1$, the family of unipotent
subgroups $\left\{ U_\lambda \; ;\; \lambda \in \Lambda \right\}$
contains only two maximal $F$--unipotent subgroups, which are
\textit{opposite} (i.e. with trivial intersection).

Each of the subgroups $U_\lambda\cap \Gamma $ is finitely
generated. Thus in order to prove thickness it suffices to
construct a family $\H$ of unconstricted subgroups of $\Gamma$
satisfying properties $(\ban_0 )$ and $(\ban_2)$, and such that
each subgroup $U_\lambda \cap \Gamma$ is contained in a subgroup
in $\H$.

The parabolic groups defined over $F$ compose a spherical building
$\Sigma$ of rank $r=\mathrm{rank}_F \bg$. Minimal parabolic groups
correspond to chambers in this building, while larger parabolic
groups correspond to panels and faces in the building. Maximal
parabolic groups correspond to vertices.

In what follows we fix a maximal $F$-split torus $T$ in $G$ and
the system of $F$--roots $\Lambda$ of $G$ with respect to $T$. Let
$\pp$ be the finite collection of all the maximal $F$--parabolic
subgroups in $G$ containing $T$. They correspond to the vertices
of an apartment in $\Sigma$. Each $P\in \pp$ decomposes as
$P=Z(T_\alpha)\widetilde{U}_\alpha$ for some $\alpha \in \Lambda$.
Let $M_P=[Z(T_\alpha),Z(T_\alpha)]$, and let $C_P$ be the
$F$--anisotropic torus such that $Z(T_\alpha)$ is an almost direct
product of $T_\alpha$ with $C_P$ and with $M_P$. Also let $D_P$ be
a maximal $F$-anisotropic torus in $M_P$. We make the choice of
$D_P$ so that if $P,P'\in \pp$ correspond to opposite vertices in
the building $\Sigma$ (in which case the corresponding unipotent
radicals have trivial intersection, while the corresponding tori
$T_\alpha$ and $T_{\alpha'}$ coincide, therefore also
$M_P=M_{P'}$) then $D_P=D_{P'}$.

Consider the solvable group $S_P = C_P D_P \widetilde{U}_\alpha$.
Since $C_P$ is an $F$-anisotropic torus, by Lemma \ref{red} the
intersection $C_P\cap \Gamma$ is a uniform lattice in $C_P$,
likewise for $D_P \cap \Gamma$ in $D_P$. Also
$\widetilde{U}_\alpha\cap \Gamma$ is a (uniform) lattice in
$\widetilde{U}_\alpha$. Consequently the semidirect product
$\Gamma_P=(C_P\cap \Gamma )(D_P\cap \Gamma)(\widetilde{U}_\alpha
\cap \Gamma )$ is a uniform lattice in $S_P$. Note that $(C_P\cap
\Gamma )(D_P\cap \Gamma)$ is never trivial. This is due on one
hand to the fact that $T_\alpha C_P D_P$ is a maximal torus in
$G$, so it has dimension at least two, and since $T_\alpha$ has
dimension one it follows that $C_P D_P$ is of dimension at least
one. On the other hand $(C_P\cap \Gamma )(D_P\cap \Gamma)$ is a
uniform lattice in $C_P D_P$.

\medskip

We show that $\Gamma$ is algebraically thick of order at most~$1$
with respect to $\H =\{ \Gamma_P \mid P\in \mathcal{P} \}$. Note
that for any $\lambda \in \Lambda$ the subgroup $U_\lambda$ is
contained in the unipotent radical $\widetilde{U}_\lambda$ of some
$P\in \pp$. In particular each $U_\lambda \cap \Gamma$ is
contained in at least one $\Gamma_P$.

Each group $\Gamma_P\, ,\, P\in \pp\, ,$ is finitely generated and
solvable, hence it is unconstricted (\cite[$\S
6.2$]{DrutuSapir:TreeGraded} , see also Section \ref{sunc},
Example 3).

Therefore it only remains to prove properties $(\ban_0)$ and
$(\ban_2)$.

\medskip

$(\ban_0 )$\quad We prove that $\Gamma_P$ is undistorted in
$\Gamma$.

\medskip

\noindent \textit{Notation:} In what follows, given a
finitely generated group $H$ we write $\dist_{H}$ to denote a word
metric on $H$. Given a Lie group $L$ we denote by $\dist_{L}$ a
 metric on $L$ defined by a left-invariant Riemannian structure.

\medskip

It suffices to prove that $\dist_{\Gamma_P} (1,g)\ll \dist_\Gamma
(1,g)$ for every $g \in \Gamma_P$. An element $g$ in $\Gamma_P$
decomposes as $g= tu$, where $t\in (C_P\cap \Gamma )(D_P\cap
\Gamma)$ and $u\in U_P \cap \Gamma$. We have that
$$\dist_{\Gamma_P} (1,tu)\leq
\dist_{\Gamma_P} (1,u)+\dist_{\Gamma_P} (1,t)\, .$$
Note that $\dist_{\Gamma_P} (1,t)\leq \dist_{(C_P\cap \Gamma
)(D_P\cap \Gamma)} (1,t)\ll \dist_G (1,t)$. The last inequality
follows from the fact that a word metric on $(C_P\cap \Gamma
)(D_P\cap \Gamma)$ is bi-Lipschitz equivalent to the restriction
of a metric from $C_P D_P$, and from the fact that $C_PD_P$ is
undistorted in $G$.

\begin{lem}\label{undistu}
For every $u\in U_P \cap \Gamma$, $\dist_{\Gamma_P} (1,u)\ll
\dist_G (1,u)$.
\end{lem}

\proof The group $U_P$ is a group of type $\widetilde{U}_\alpha$
for some $\alpha \in \Lambda$. That is, if $\Lambda_\alpha^+
=\left\{ \lambda\in \Lambda \; ;\; m_\alpha (\lambda) >0
\right\}$, then the group $U_P$ has Lie algebra
$\widetilde{{\mathfrak{u}}}_\alpha =
\bigoplus_{\widetilde{\lambda}|_T \in
\Lambda_\alpha^+}\fg_{\widetilde{\lambda}}$. In particular
$\widetilde{{\mathfrak{u}}}_\alpha = \bigoplus_{\lambda \in
\Lambda_\alpha} {\mathfrak{u}}_\lambda $, where $\Lambda_\alpha$
is such that any root in $\Lambda_\alpha^+$ is either contained in
or is the double of a root in $\Lambda_\alpha$.

\begin{lem}[Lubotzky--Mozes--Raghunathan
    \cite
    {LMR:metrics}, $\S 3$]\label{dec}
Let $\lambda_1,\lambda_2,...., \lambda_N$ be the enumeration of
the roots in $\Lambda_\alpha$ in increasing order. The order here
is the lexicographic order with respect to some basis $\Delta$ of
$\Lambda$ having $\alpha$ as first root.
\begin{itemize}
\item[(1)] There exist morphisms $f_i : U_P \to U_{\lambda_i}\, ,\,
1\leq i\leq N\, ,\,
$ such that for every $u\in U_P\, ,\, u=f_1(u)\cdot f_2(u)\cdot
\ldots \cdot
f_N(u)$ and
$$
\dist_G(1,u)\leq \sum_{i=1}^N \dist_G(1,f_i(u))\ll \dist_G(1,u)\,
.
$$
\item[(2)] If $\Gamma_1$ is a suitable congruence subgroup of
$\Gamma = \mathbf{G}(\oo_S)$, then for every $u\in \Gamma_1
\cap U_P$ the components $f_i (u)$ are in $U_{\lambda_i}\cap \Gamma$
for all
$i=1,2,...,N$.
\end{itemize}
\end{lem}

Let $\Gamma_1$ be as in Lemma \ref{dec}. It has finite index in
$\Gamma$, therefore it suffices to prove Lemma \ref{undistu} for
$u\in U_P \cap \Gamma_1$. In this case $f_i (u)\in
U_{\lambda_i}\cap \Gamma$ for all $i=1,2,...,N$, and
$\dist_{\Gamma_P} (1,u)\leq \sum_{i=1}^N
\dist_{\Gamma_P}(1,f_i(u))$. By Lemma \ref{dec}, (1), it will then
suffice to prove Lemma \ref{undistu} for each $f_i(u)$. Hence we
may assume in what follows that $u\in U_\lambda \cap \Gamma$, for
some $\lambda \in \Lambda_\alpha$.

Consider the solvable subgroup $S_\lambda=C_PD_P U_\lambda$ of
$S_P$, and its uniform lattice $\Gamma_\lambda = (C_P \cap \Gamma
) (D_P \cap \Gamma) (U_\lambda \cap \Gamma )$, which is a subgroup
of $\Gamma_P$. It will suffice to prove that
$\dist_{\Gamma_\lambda}(1,u)\ll \dist_G (1,u)$, which is
equivalent to proving that $\dist_{S_\lambda}(1,u)\ll\dist_G
(1,u)$. With a decomposition similar to the one in Lemma
\ref{dec}, (1), we can reduce the problem to the case when $u$ is
in the uniparametric unipotent subgroup $U_{\widetilde{\lambda}}$
for some root $\widetilde{\lambda}$ such that its restriction to
$T$ is $\lambda$. The torus $C_PD_P$ is orthogonal to the
one-dimensional $F$-split torus $T_\alpha$ associated to $P$, in a
maximal torus containing both. If the root $\widetilde{\lambda}$
would be constant equal to $1$ on $C_PD_P$, then $\lambda$ would
be constant equal to $1$ on $C_PD_P\cap T$, hence the same would
be true for $\alpha$. Consequently $\alpha$ would be equal to $1$
on the orthogonal of $T_\alpha$ in $T$. This implies that the
$F$-structure on $G$ is reducible (see for instance
\cite{KleinerLeeb:buildings} for a geometric argument), which
implies that $\Gamma$ is reducible, contradicting the hypothesis.
Thus there exists at least one uniparametric semisimple subgroup
$T_1$ in $C_PD_P$ on which $\widetilde{\lambda}$ takes all
positive values. An argument as in \cite[$\S
3.D$]{Gromov:Asymptotic} then implies that $\dist_{S_\lambda}(1,u)
\leq \dist_{T_1U_{\widetilde{\lambda}}}(1,u)\ll \ln (1+ \| u - I
\|)\ll \dist_G (1,u)$.\endproof

Lemma \ref{undistu} together with the considerations preceding it
imply that $$\dist_{\Gamma_P} (1,tu)\ll\dist_G
(1,t)+\dist_G(1,u)\, .$$ On the other hand $\dist_G
(1,t)+\dist_G(1,u) \ll \dist_G(1,tu)$. This follows from the well
known fact that $\dist_G (1,t)\leq \dist_G(1,tu)$ and from the
triangular inequality $\dist_G(1,u) \leq \dist_G
(1,t)+\dist_G(1,tu)$. Then $\dist_{\Gamma_P} (1,tu)\ll \dist_{G}
(1,tu)\ll \dist_{\Gamma} (1,tu)$, where the latter estimate
follows from Theorem \ref{tlmr}. This completes the proof of
$(\ban_0)$.

\medskip

$(\ban_2)$\quad First we suppose that $\mathrm{rank}_F \bg=1$.
Then $\pp$ has only two elements, $P$ and $P'$, which are
opposite. Consequently $\Gamma_P \cap \Gamma_{P'}$ contains
$(C_P\cap \Gamma)(D_P\cap \Gamma )$, which is a lattice in a torus
of dimension at least one, hence it is infinite.

\medskip

Suppose now that $\mathrm{rank}_F \bg\geq 2$. This implies that
the building $\Sigma$ composed of $F$-parabolics has rank at least
two, therefore it is connected.
Let $P,P'\in \pp$. The groups $P$ and $P'$ seen as vertices in the
same apartment in $\Sigma$ can be connected by a finite gallery of
chambers in the same apartment. This gallery is represented by a
sequence of minimal $F$-parabolic subgroups $ B_1, B_2,..., B_k\,
, $ with $B_1< P$ and $B_k < P'$. For each $i=1,2,..., k-1$ there
exists $P_i\in \pp$ such that both $B_i$ and $B_{i+1}$ are
contained in $P_i$. In the spherical building $\Sigma$ the group
$P_i$ represents a vertex of the panel that the chamber $B_i$ and
the chamber $B_{i+1}$ have in common. Thus one obtains a sequence
of maximal parabolics in $\pp$, $P_0=P, P_1,P_2,..., P_{k-1},
P_k=P'$. For each $i=0,1,2,..., k-1\, ,$ the intersection of the
respective unipotent radicals $U_{P_i}$ and $U_{P_{i+1}}$ of $P_i$
and $P_{i+1}$ contains the center $U_{i,i+1}$ of the unipotent
radical of $B_i$. The group $U_{i,i+1}$ can be written as
$U_\alpha$ with $\alpha$ the maximal positive root in the basis
corresponding to the chamber $B_i$, in particular it is defined
over $F$. Hence $\Gamma_{P_i} \cap \Gamma_{P_{i+1}}$ contains
$U_{i,i+1}\cap \Gamma$, which is a lattice in $U_{i,i+1}$. We
conclude that $\Gamma_P$ and $\Gamma_{P'}$ are thickly connected
by the sequence $\Gamma_{P_0}=\Gamma_P\, ,\, \Gamma_{P_1}\, ,\,
\Gamma_{P_2}\, ,\, ....,\, \Gamma_{P_k}=\Gamma_{P'}$.
\end{pf}

\begin{qn} Are non-uniform higher rank lattices unconstricted?
\end{qn}



\end{document}